\documentclass[10pt]{article}
\usepackage{amsfonts,amsmath, amssymb,latexsym,mathrsfs,upgreek}
\usepackage{multirow}
\usepackage{verbatim}
\usepackage{leftindex}
\usepackage{colortbl}
\usepackage{longtable}
\usepackage{tabularx,ragged2e}

\usepackage{hyperref}
\usepackage{eucal}
\usepackage{tocloft}
\usepackage{cancel}
\usepackage[british]{babel}
\usepackage[new]{old-arrows} 
\usepackage{float}
\usepackage{bibspacing,mathtools}
\usepackage{tikz}
\usetikzlibrary{matrix,arrows,decorations.pathreplacing}

\usepackage{todonotes,diagbox}

\usepackage{amsthm}
\usepackage{mathtools}
\usepackage{manfnt}
\usepackage{enumitem}
\usepackage{setspace}
\usepackage{calc,footnote}
\usepackage{graphicx}
\usepackage{tikz}
\usepackage{tikz-cd}
\usepackage{comment}

\newtheorem{thm}{Theorem}

\newtheorem{theorem}{Theorem}[section]
\newtheorem{proposition}[theorem]{Proposition}
\newtheorem{lemma}[theorem]{Lemma}
\newtheorem{corollary}[theorem]{Corollary}

\theoremstyle{definition}
\newtheorem{defn}[theorem]{Definition}
\newtheorem{remark}[theorem]{Remark}
\newtheorem{rem}[thm]{Remark}

\makeatletter
\providecommand{\proofNameStyle}{\bfseries}
\renewenvironment{proof}[1][\proofname]{\par\pushQED{\qed}%
  \normalfont \topsep6\p@\@plus6\p@\relax
  \trivlist
  \item[\hskip\labelsep
        \proofNameStyle
    #1\@addpunct{:}]\ignorespaces
}{%
  \popQED\endtrivlist\@endpefalse
}
\newtheorem*{theorem*}{Theorem}
\newtheorem*{proposition*}{Proposition}
\newtheorem*{lemma*}{Lemma}

\theoremstyle{remark}

\DeclareMathOperator{\Cl}{Cl}

\DeclareMathOperator{\Sym}{Sym}

\DeclareMathOperator{\Stab}{Stab}
\def\Z{{\mathbb Z}}

\def\B{{\mathcal B}}
\def\W{{\mathcal W}}
\def\Z{{\mathbb Z}}
\def\irr{{\rm irr}}

\def\bx{{\rm bal}}

\def\wei{{\rm wei}}
\def\Res{{\rm Res}}
\def\GL{{\rm GL}}
\def\SL{{\rm SL}}
\def\Cl{{\rm Cl}}

\def\Sym{{\rm Sym}}
\def\det{\operatorname{det}}

\def\Stab{{\rm Stab}}

\def\Avg{{\rm Avg}}

\def\P{{\mathbb P}}
\def\Disc{{\rm Disc}}

\def\nd{{\rm nd}}
\def\Vol{{\rm Vol}}

\def\R{{\mathbb R}}

\def\F{{\mathbb F}}

\def\FF{{\mathcal F}}
\def\cA{{\mathcal A}}

\def\Q{{\mathbb Q}}

\def\Z{{\mathbb Z}}
\def\P{{\mathbb P}}
\def\F{{\mathbb F}}
\def\Q{{\mathbb Q}}
\def\O{{\mathcal O}}

\def\cB{{\mathcal B}}

\def\Z{{\mathbb Z}}

\def\cO{{\mathcal O}}

\def\rmH{\mathrm{H}}
\def\Sel{\mathrm{Sel}}

\newcommand{\SO}{\mathrm{SO}}

\newcommand{\J}{\mathcal{J}}

\newcommand{\rhoright}{\rho}
\newcommand{\rholeft}{\lambda}

\newcommand{\nc}{\newcommand}
\nc{\on}{\operatorname}
\nc{\renc}{\renewcommand}
\nc{\wt}{\widetilde}
\nc \defeq {\vcentcolon=}
\nc{\eqdef}{=\vcentcolon}
\nc{\Spec}{\on{Spec}}
\nc{\ol}{\overline}
\renc{\d}{\partial}
\nc{\mc}{\mathcal}

\allowdisplaybreaks

\title{\vspace*{-1.5cm}Counting integral points on symmetric varieties with \\ applications to arithmetic statistics}
\author{Arul Shankar, Artane Siad, and Ashvin A.~Swaminathan}
\date{\today}

\begin{document}

\maketitle

\vspace*{-0.6cm}
\begin{abstract}
In this article, we  
combine Bhargava's geometry-of-numbers methods with the dynamical point-counting methods of Eskin--McMullen and Benoist--Oh to develop a new technique for counting integral points on symmetric varieties lying within fundamental domains for coregular representations. As applications, we study the distribution of the $2$-torsion subgroup of the class group in thin families of cubic number fields, as well as the distribution of the $2$-Selmer groups in thin families of elliptic curves over $\Q$. For example, our results suggest that the existence of a generator of the ring of integers with small norm has an increasing effect on the average size of the $2$-torsion subgroup of the class group, relative to the Cohen--Lenstra predictions.
\end{abstract}

\setcounter{tocdepth}{1}
\tableofcontents

\section{Introduction}

The foundational heuristics of Cohen and Lenstra~\cite{MR756082}, formulated in the 1980s, constitute a system of conjectures that predict the distribution of the class groups of quadratic fields. These heuristics are based on the philosophy that algebraic objects occur in nature with probability inversely proportional to the sizes of their automorphism groups.
In the years since, the Cohen--Lenstra heuristics have enjoyed an extensive series of generalisations and elaborations. Some notable examples include the works of Cohen--Martinet~\cite{MR866103,MR1037430,MR1226813}, who generalised the heuristics to extensions of number fields of arbitrary degree, the work of Malle~\cite{MR2778658}, who realised that adjustments to the Cohen--Lenstra--Martinet heuristics were needed to account for roots of unity in the base field, and the very recent work of Wood--Sawin~\cite{2301.00791}, which proposes the necessary adjustments. See also the related works of Gerth~\cite{MR759260,MR887792}, Wood~\cite{MR3858473}, Liu--Wood--Zureick-Brown~\cite{1907.05002}, Bartel--Lenstra~\cite{MR4105790}, \mbox{Bartel--Johnston--Lenstra~\cite{2005.11533}, and Wang--Wood~\cite{MR4277275}.}

Despite these significant advances on the conjectural side, proofs have been discovered for only a handful of the predictions made by the Cohen--Lenstra--Martinet heuristics. In the 1970s (predating the Cohen--Lenstra heuristics!), the seminal work of Davenport and Heilbronn \cite{MR491593} determined the average size of the $3$-torsion subgroups of the class groups of quadratic number fields ordered by discriminant. Some three decades later, the landmark work of Bhargava determined the average size of the $2$-torsion subgroups of the class groups of cubic number fields~\cite{MR2183288} ordered by discriminant. In fact, as shown in subsequent works of Bhargava--Varma~\cite{MR3369305,MR3471250}, the aforementioned averages remain stable under the imposition of general sets of \emph{local} conditions. More recently, the stunning work of Lemke Oliver--Wang--Wood~\cite{2110.07712} determined the average size of the $3$-torsion subgroups of certain relative class groups associated to $G$-extensions of $K$, where $G \subset S_{2^m}$ is any transitive permutation $2$-group containing a transposition and $K$ is any number field; their result is also expected to be insensitive to the imposition of local conditions.

On the other hand, the question of how \emph{global} conditions impact class group distributions is much less well-understood and has developed into an active area of research in recent years. In particular, an intriguing picture has emerged of the behaviour of the $2$-primary part of the class group under the condition of being {\it monogenised} --- i.e., having ring of integers generated by one element as a $\Z$-algebra. First considered in work of Bhargava--Hanke--Shankar \cite{BHSpreprint} and further studied in works of Siad \cite{Siadthesis1,Siadthesis2}, Swaminathan \cite{swacl}, and Bhargava--Shankar--Swaminathan \cite[\S4.B]{swathesis}, the upshot of this line of inquiry is as follows: when considered as a ``random variable'' valued in finite abelian $2$-groups, the distribution of the $2$-primary parts of the class groups over the thin family of monogenised degree-$n$ fields appears to be different from the distribution predicted (for the {\it full} family of degree-$n$ fields) by the Cohen--Lenstra--Martinet heuristics. 

In this paper, 
we study how the distribution of the $2$-torsion subgroups of the class groups of \emph{cubic} number fields is impacted by the condition of being \emph{unit-monogenised} --- i.e., having ring of integers generated over $\Z$ by a unit --- among other related conditions. Our work suggests that imposing these global conditions leads to class group distributions that are different from the corresponding distributions for both the full family of cubic fields and the subfamily of monogenised fields. We also prove analogous results concerning the distribution of the $2$-Selmer groups in thin families of elliptic curves (see \S\ref{sec-selrel}), which indicate that these distributions sometimes differ from those predicted by Poonen--Rains \cite{MR2833483}. 
To accomplish this, we introduce tools from dynamics into arithmetic statistics. Specifically, we combine Bhargava's geometry-of-numbers methods --- which give a systematic procedure for counting integral orbits of representations --- with the works of Eskin--McMullen~\cite{MR1230290} and Benoist--Oh~\cite{MR3025156} --- which utilise the exponential mixing of a semisimple Lie group on its finite volume quotient to count integral points on symmetric varieties.

\subsection{Results on unit-monogenised cubic fields} \label{sec-unitresult}
A {\it unit-monogenised} cubic field is a pair $(K,\alpha)$ consisting of a cubic field $K/\Q$ with a \emph{unit monogeniser} $\alpha$ of its ring of integers --- i.e., an $\alpha \in \mc{O}_K^\times$ such that $\O_K = \Z[\alpha]$. A natural way to order such pairs $(K,\alpha)$ is by the \emph{balanced height} $\rmH_{\on{bal}}(K,\alpha)$, which is defined to be the maximum of the sizes of the coefficients of the minimal polynomial of $\alpha$. 
The balanced height gives a natural ordering since it is comparable to the discriminant $\Delta$; indeed, when unit-monogenised cubic fields $K$ are ordered by balanced height, a density of $100\%$ of them satisfy $\rmH_{\on{bal}}(K,\alpha)^4/2 \leq \Delta(K) \leq 2\rmH_{\on{bal}}(K,\alpha)^4$.\footnote{This is similar to the height on monogenised cubic fields in \cite{BHSpreprint}, which was also chosen to mimic the discriminant. In fact, it is easy to see that the constants $1/2$ and $2$ can be replaced with any constants less than $1$ and greater than $1$, respectively.}

%indeed, there exist constants $c_1, c_2 > 0$ such that when unit-monogenised cubic fields $K$ are ordered by balanced height, a density of $100\%$ of them satisfy $c_1\rmH_{\on{bal}}(K,\alpha)^4 \leq \Delta(K) \leq c_2\rmH_{\on{bal}}(K,\alpha)^4$.

For unit-monogenised cubic fields ordered by balanced height, we prove the following theorem, which describes the average size of the $2$-torsion subgroups of their class groups:

\begin{thm} \label{first main} 
When totally real unit-monogenised cubic fields are ordered by balanced height, the average size of the $2$-torsion subgroups of their class groups is at most $2 + 3/56$.
\end{thm} 

\noindent Although we only prove an upper bound in Theorem \ref{first main}, we expect that the average is in fact \emph{equal} to $2 + 3/56$; see Remark \ref{rmk:unif estimate}. We restrict to totally real fields in Theorem \ref{first main} because when unit-monogenised cubic fields are ordered by balanced height, a density of $100\%$ of them are totally real (i.e., the number of complex fields has an asymptotically lower order of growth); in fact, Theorem \ref{first main} continues to hold with the ``totally real'' condition removed. To study complex unit-monogenised cubic fields, we need to order pairs $(K,\alpha)$ in such a way that the complex fields occur with positive proportion. A natural way to do this is to order pairs $(K,\alpha)$ by the ``size'' of $\alpha$ in $\R \otimes_{\Q} K$. More precisely, given a pair $(K,\alpha)$, let $f(x)=x^3+ax^2+bx\pm 1$ be the minimal polynomial of $\alpha$; then we define the {\it weighted height} of $(K,\alpha)$ by $\rmH_{\on{wei}}(K,\alpha) \defeq \max\{|a|,\sqrt{|b|}\}$. 
For every unit-monogenised cubic field $(K,\alpha)$, the weighted height satisfies $\max_v |\alpha|_v \asymp c \times \rmH_{\on{wei}}(K,\alpha)$, where $c > 0$ is a constant and $v$ ranges over archimedean places of $K$. This estimate (for polynomials of any degree) was claimed by Lagrange in~\cite[\S2.IV]{MR0439547}; see~\cite[Theorem~3]{MR3320796} for a proof. For unit-monogenised cubic fields ordered by weighted height, we prove the following analogue of Theorem~\ref{first main}:

\begin{thm} \label{first main2} 
Let $K$ run through all unit-monogenised cubic fields \mbox{ordered by weighted height. Then:}
\vspace*{-4pt}
\begin{enumerate}[itemsep=-2pt,leftmargin=17pt]
\item[$(\mathrm a)$] Over totally real $K$, the average size of the $2$-torsion subgroup of the class group of $K$ is at most $2 + 3/364$.
\item[$(\mathrm b)$] Over complex $K$, the average size of the $2$-torsion subgroup of the class group of $K$ is at most $3 + 3/14$.
\end{enumerate}
\end{thm} 

\noindent Once again, we expect that the average sizes in Theorem~\ref{first main2} are in fact {\em equal} to $2+3/364$ and $3+3/14$; see Remark~\ref{rmk:unif estimate} for more detail. In Table \ref{tab-avgs}, we list the average size of the $2$-torsion subgroups of the class groups of cubic fields having a specified real signature running through the full family (as determined by Bhargava~\cite{MR2183288}), the monogenised family (as determined by Bhargava--Hanke--Shankar~\cite{BHSpreprint}, and conditional values for the unit-monogenised family (as suggested by Theorems~\ref{first main} and~\ref{first main2}). The symbol $\dagger$ indicates that our result is for weighted height only. The star indicates an archimedean volume factor that is height-dependent; e.g., this factor is equal to $1/4$ for balanced height and to $1/26$ for weighted height.

\begin{table}[h]
\centering
\begin{tabular}{|c||*{4}{c|}{c|}}\hline
{\bf Family of cubics}
& \multicolumn{1}{|c|}{\bf Avg $2$-torsion, totally real} & {\bf Avg $2$-torsion, complex} \\\hline\hline
\vphantom{$\frac{2^{r^{r^{r^r}}}}{2_{r_{r_{r_{r}}}}}$} Full family & $1 + 1/4$ & $1 + 1/2$ \\ \hline
\vphantom{$\frac{2^{r^{r^{r^r}}}}{2_{r_{r_{r_{r}}}}}$} Monogenised family & $1 + 1/2$ & $1 + 1$ \\\hline
\vphantom{$\frac{2^{r^{r^{r^r}}}}{2_{r_{r_{r_{r}}}}}$} Unit-monogenised family & $1 + \displaystyle 1 + 3/14 \times *$ & $1 + 2 + 3/14^\dagger$ \\ \hline
\end{tabular}
\caption{Expected averages for cubic fields, monogenised cubic fields, and unit-monogenised cubic fields}
\label{tab-avgs}
\end{table}

Theorem~\ref{first main2} suggests a phenomenon that, to the authors' knowledge, has not been observed previously for families of $S_n$-number fields:  
 there is a lack of conformity between the averages in the totally real and complex cases. 
This constitutes a deviation from the Cohen--Lenstra--Martinet philosophy, which predicts that the class group of a random totally real cubic field behaves like the quotient, by a random element, of the class group of a random complex cubic field. If this philosophy were true for a family $F$ of cubic fields, then we would expect that the average number of nontrivial elements in the 2-torsion subgroups of the class groups of totally real fields in $F$ would be exactly half the corresponding average over complex fields in $F$. From Table \ref{tab-avgs}, we can see that this holds for the family $F$ of all cubic fields (by Bhargava's work \cite{MR2183288}) and the family of monogenised cubic fields (by \cite{BHSpreprint}). It also holds for the family of $n$-monogenised cubic fields (see \cite{BHSpreprint}), and, conditional on certain tail estimates, for various families of higher-degree number fields (see \cite{MR3782066,Siadthesis1,Siadthesis2,swacl}). However, our work indicates that when unit-monogenised cubic fields are ordered by weighted height, these averages are $1+3/364$ and $2+3/14$ in the totally real and complex families, respectively. Observe that these averages are not in the ratio of $1 : 2$ predicted by the Cohen-Lenstra-Martinet philosophy.
%, we conclude that, relative to the totally real case, there is ``extra'' $2$-torsion in the class groups of complex unit-monogenised fields when ordered by weighted height. 

On a related note, the average values $2+3/56$, $2+3/364$, and $3+3/14$ occurring in Theorems~\ref{first main} and~\ref{first main2} seem quite strange. By contrast, the averages for the full and monogenised families stated in Table~\ref{tab-avgs} appear more \mbox{``natural'' --- for} instance, they are all dyadic rational numbers. As we explain further on in the introduction (see \S\ref{sec-partitioning}), the family of unit-monogenised cubic fields can be partitioned into what we call ``stable families,'' on each of which the average $2$-torsion in the class group takes a natural dyadic rational value.  The averages in Theorems~\ref{first main} and~\ref{first main2} then arise as weighted averages of these natural values, where each stable family is weighted by its density. Furthermore, we will see that the deviation from the Cohen--Lenstra--Martinet philosophy described in the previous paragraph occurs across certain stable families, but not others.

\subsection{Results on more general thin families of cubic fields} \label{sec-stat}
 Our new approach also allows us to study the distribution of the $2$-torsion in the class groups of more general thin families of cubic fields. 
Before we state our results in this direction, we first generalise the notion of a unit-monogenised cubic field. For this purpose, let $U$ be the affine scheme over $\Z$ whose $R$-points consist of binary cubic forms over $R$ for any ring $R$. To each form $f \in U(R)$ is a naturally associated cubic extension $R_f/R$ (see \S\ref{sec-rep}). When $R = \Z$ or $\Z_p$ for a prime $p$, we say that $f \in U(R)$ is \emph{maximal} if $f$ has nonzero discriminant and the corresponding cubic ring $R_f$ is integrally closed in $\on{Frac}(R) \otimes_R R_f$.

Next, take $a,d \in \Z \smallsetminus \{0\}$, and let $U_{a,d} \subset U$ be the subscheme of $U$ consisting of forms with $x^3$-coefficient $a$ and $y^3$-coefficient $d$. Then the set of unit-monogenised cubic fields $(K,\alpha)$ is in bijection with the set of maximal irreducible forms $f \in U_{1,\pm1}(\Z)$ via the map that sends $\alpha$ to the homogenisation of its characteristic polynomial. Thus, Theorems~\ref{first main} and~\ref{first main2} may be interpreted as giving the average size of the $2$-torsion in the class group of $R_f$ over the family $\{f \in U_{1,\pm 1}(\Z) : \text{$f$ is maximal, irreducible}\}$. In this section, we generalise these theorems to the family $U_{a,d}(\Z)_{\max} \defeq \{f \in U_{a,d}(\Z) : \text{$f$ is maximal, irreducible}\}$, where $a,d \in \Z \smallsetminus \{0\}$ are any fixed nonzero integers. We use the following generalisations of the balanced and weighted heights defined in \S\ref{sec-unitresult}: given $f(x,y) = ax^3 + bx^2y+cxy^2+dy^3 \in U_{a,d}(\R)$, the \emph{balanced height} is $\rmH_\bx(f)  \defeq  \max \{|b|,|c|\}$ and the \emph{weighted height} is $\rmH_{\wei}(f)  \defeq  \max \{ \left|b\right|, \left|c\right|^{1/2}\}$. 

In fact, we go one step further by stating our results for subfamilies of $U_{a,d}(\Z)_{\max}$ defined by quite general infinite sets of local conditions. To this end, we pick a sign $\varepsilon \in \{\pm\}$ and define
\begin{equation*}
U_{a,d}(\R)^{\varepsilon} \defeq \{f \in U_{a,d}(\R) : \varepsilon\on{Disc}(f) > 0\}. 
\end{equation*}
For each prime prime $p$, we pick an open subset $\Sigma_p\subset U_{a,d}(\Z_p)$
whose boundary has measure $0$ with respect to Haar measure on $U_{a,d}(\mathbb{Z}_p) \simeq \Z_p^2$. To this collection $\varepsilon,\,\{\Sigma_p\}_p$, we associate the subset $\Sigma=U_{a,d}(\R)^{\varepsilon}\cap\bigcap_p\Sigma_p$ of $U(\Z)$. A set of binary cubic forms arising this way is said to be {\it acceptable} if $\Sigma_p \supset \{f \in U_{a,d}(\Z_p) : p^2 \nmid \on{Disc}(f)\}$ for each $p \gg 1$.

 Write $a = a_k  a_m^2$ and $d = d_k  d_m^2$, where $a_k, a_m, d_k, d_m \in \Z$ with $a_k, d_k$ squarefree.\footnote{Note that $a_k$ and $a_m$ need not be coprime, and the same is also true of $d_k$ and $d_m$.} A form $f \in U_{a,d}(\Z_p)_{\max}$ is said to be \emph{right-} (resp., \emph{left-}) \emph{sufficiently-ramified} at $p$ if $p \mid a_k$ (resp., $p \mid d_k$) and $f(x,1) \equiv \alpha \cdot g(x)^2 \pmod p$ (resp., $f(1,y) \equiv \alpha \cdot g(y)^2 \pmod p$) for some $\alpha \in \F_p^\times$ and $g \in \F_p[x]$ (resp., $g \in \F_p[y]$). We say that $f$ is \emph{right-} (resp., \emph{left-}) \emph{sufficiently-ramified} if $f$ is right- (resp., left-) sufficiently-ramified at $p$ for every $p \mid d_k$ (resp., $p \mid a_k$).\footnote{In particular, if $|d_k| = 1$ (resp., $|a_k| = 1$), then $f$ is automatically right- (resp., left-) sufficiently-ramified.} For an acceptable family $\Sigma \subset U_{a,d}(\Z)_{\max}$, the \emph{right-} (resp., \emph{left-}) \emph{sufficient-ramification density} of $\Sigma$ at $p$ is denoted $\rhoright_\Sigma(p)$ (resp., $\rholeft_\Sigma(p)$) and defined as the density of right- (resp., left-) sufficiently-ramified forms in $\Sigma_p$. Set 
 \begin{equation*}
\rhoright_\Sigma \defeq \prod_{p \mid a_k} \rhoright_\Sigma(p);\quad\quad
\rholeft_\Sigma \defeq \prod_{p \mid d_k} \rholeft_\Sigma(p).
 \end{equation*}
\noindent Given the above setup, we have the following theorems, which generalise the results of \S\ref{sec-unitresult} to cubic rings associated to binary cubic forms in an acceptable family in $U_{a,d}(\Z)_{\max}$:

\begin{thm}
\label{second main real}  Let $\rmH \in \{\on{H}_{\on{bal}}, \on{H}_{\on{wei}}\}$, and let $\Sigma \subset U_{a,d}(\Z)_{\max}$ be an acceptable family of forms of positive discriminant. Then we have
\begin{alignat*}{3}
&\limsup_{X \to \infty} \underset{\substack{f \in \Sigma \\ \rmH(f) < X}}{\on{Avg}}\,\, \#\on{Cl}(R_f)[2] \hspace{1ex} & \le\hspace{1ex} & \frac{5}{4}  + \frac{\rhoright_\Sigma + \rholeft_\Sigma + \chi_{a,d} \,\, \rhoright_\Sigma \rholeft_\Sigma }{4}+ \mathfrak{d}_\infty^{(0)}\mathfrak{d}_\Sigma ,
\end{alignat*}
where $\chi_{a,d} = 1$ if $\gcd(a_k,d_k) = 1$ and \mbox{$\chi_{a,d} = 0$ otherwise,} and where $\mathfrak{d}_\infty^{(0)},\, \mathfrak{d}_\Sigma$ are \mbox{the densities defined in \S\ref{sec-dist}.}
\end{thm}

\begin{thm}
\label{second main complex}  Let $\Sigma \subset U_{a,d}(\Z)_{\max}$ be an acceptable family of forms of negative discriminant. Then we have
\begin{alignat*}{3}
&\limsup_{X \to \infty} \underset{\substack{f \in \Sigma \\ \rmH_{\on{wei}}(f) < X}}{\on{Avg}}\,\, \#\on{Cl}(R_f)[2] \hspace{1ex} & \le\hspace{1ex} &  \frac{3}{2}  + \frac{\rhoright_\Sigma + \rholeft_\Sigma + \chi_{a,d} \,\, \rhoright_\Sigma \rholeft_\Sigma }{2} + \mathfrak{d}_\Sigma,
\end{alignat*}
where $\mathfrak{d}_\Sigma$ is the density defined in \S\ref{sec-dist}.
\end{thm}
 As was the case for Theorems~\ref{first main} and~\ref{first main2}, we expect that the bounds in Theorems~\ref{second main real} and~\ref{second main complex} are in fact the exact averages (i.e., the theorems should hold with ``$\limsup$'' replaced by ``$\lim$'' and ``$\leq$'' replaced by ``$=$''); see Remark \ref{rmk:unif estimate}. The formulae given in Theorems~\ref{second main real} and~\ref{second main complex} may be evaluated explicitly in various cases of interest. For example, take $\Sigma = U_{a,d}(\Z)_{\max}$ to consist of all maximal cubic orders in $U_{a,d}(\Z)$. Formulas for the associated sufficient-ramification densities $\rhoright_\Sigma$ and $\rholeft_\Sigma$ are determined in Appendix~\ref{appendix-localcomputationsat2}; see~\eqref{eq-sqdens}. From these formulas, one readily observes that the average over all $a \in \Z \smallsetminus \{0\}$ (resp., $d \in \Z \smallsetminus \{0\}$) of $\rhoright_\Sigma$ (resp., $\rholeft_\Sigma$) is equal to zero. This observation has two noteworthy consequences. The first consequence is that, on average over all $a,d \in \Z \smallsetminus \{0\}$, the formulae given in Theorems~\ref{second main real} and~\ref{second main complex} respectively converge to $5/4$ and $3/2$. As is to be expected, these are precisely the averages determined by Bhargava~\cite{MR2183288} (cf.~Ho--Shankar--Varma~\cite{MR3782066}) for the \mbox{full family of cubic fields.} 

The second consequence is that, on average over all $d \in \Z \smallsetminus \{0\}$, the formulae given in Theorems~\ref{second main real} and~\ref{second main complex} respectively converge to $5/4 + 1/4\sigma(a_k)$ and $3/2 + 1/2\sigma(a_k)$, where $\sigma$ denotes the divisor function. Unsurprisingly, these are precisely the averages determined by Bhargava--Hanke--Shankar~\cite{BHSpreprint} for the family of cubic fields corresponding to binary cubic forms with $x^3$-coefficient equal to $a$.

We conclude this subsection by noting that the number of times each cubic field arises from a form in $U_{a,d}(\Z)_{\max}$ is bounded in terms of $a$ and $d$. Indeed, it follows from Baker's work~\cite{MR228424} on the effective solution of an arbitrary Thue equation that the fibres of the map taking a form $f \in U_{a,d}(\Z)_{\max}$ to its corresponding cubic field are finite and bounded in a way that depends only on $a$ (or on $d$). Furthermore, $a = 1$ or $d = 1$, it follows from work of Bennett~\cite{Ben2001} that the fibres of this map are of size at most $10$. Even more can be said for complex fields when \mbox{$a = d = 1$ --- here,} it follows from work of Nagell~\cite{MR1544935} and Delone~\cite{MR1545095} that the fibres have size at most $5$, and in fact at most $3$ for all but finitely many cubic fields!

\subsection{Partitioning into stable families} \label{sec-partitioning}

The work of Bhargava--Varma \cite{MR3369305,MR3471250} implies that, once we partition the family of quadratic (resp., cubic) fields into subfamilies depending on their real signature, the average size of the $3$-torsion (resp., $2$-torsion) subgroup of the class group does not change when any further finite sets of local conditions are imposed. As shown in \cite{BHSpreprint}, the same result is also true if we replace the family of cubic fields with the family of $n$-monogenic cubic fields, when $n$ is a perfect square. However, the story is more complicated for the family of $n$-monogenic cubic fields, when $n$ is not a perfect square. In this situation, local specifications at $n$ can change the average $2$-torsion. However, the following is true: suppose we fix the real signature (either totally real, or complex) and further partition the family of $n$-monogenic cubic fields into two subfamilies, the family of fields that are left-sufficiently-ramified at $n$ and the family of fields that are not left-sufficiently-ramified at $n$. Then in all four of these subfamilies, the average $2$-torsion does not change when further local conditions are imposed. We call such families {\it stable families}. Then the full family of cubic fields can be partitioned into two stable subfamilies according to their infinite signature, and the corresponding average values are $5/4$ and $3/2$. Similarly, the family of $n$-monogenic cubic fields, where $n$ is a perfect square, can be partitioned into two stable subfamilies with averages $3/2$ and $2$. The family of $n$-monogenic cubic fields, where $n$ is not a perfect square, can be partitioned into four stable subfamilies, with averages $5/4$ and $3/2$ in the totally real case, and $3/2$ and $2$ in the complex case.

The families considered in this paper can also be partitioned into stable subfamilies, by imposing local conditions at infinity and at the primes dividing $2ad$. At the infinite place, we start as usual by dividing into the families of totally real and complex cubic fields. We denote these by $F_{a,d}^\pm$. Then at the primes dividing $a$ and $d$, we break up these two families into three further subfamilies depending on whether we have neither left-sufficient-ramification at $a$ nor right-sufficient-ramification at $d$, one of the two possible sufficient-ramifications, or both. We denote these families by $F_{a,d}^\pm(\lambda)$, where $\lambda\in\{0,1,2\}$ indicates the number of sufficient-ramifications. When $\lambda \in \{0,1\}$, the corresponding subfamilies are already stable. However, when $\lambda=2$, there is a further phenomenon to consider, namely, the existence or non-existence of an additional class in the class group termed the ``$\Delta$-distinguished'' class; for more on this class, see \S\ref{sec-discuss} and \S\ref{sec-dist} (to follow). The existence of this class can have two possible obstructions: at the primes dividing $2ad$ or, in the case of totally real cubic fields, at infinity. We divide $F_{a,d}^-(2)$ into the subfamilies $F_{a,d}^+(2;\delta_{\on{fin}})$ (with $\delta_{\on{fin}}\in\{0,1\}$) and we divide $F_{a,d}^+(2)$ into the subfamilies $F_{a,d}^+(2;\delta_{\on{fin}},\delta_\infty)$ (with $\delta_{\on{fin}},\delta_\infty\in\{0,1\}$), where $\delta_{\on{fin}}$ (resp., $\delta_\infty$) being $1$ implies there is no obstruction at the primes dividing $2ad$ (resp., infinity), and being $0$ implies that there is an obstruction. The exact conditions for determining the values of $\delta_{\on{fin}}$ and $\delta_\infty$ of an $(a,d)$-monogenised cubic field are given in \S\ref{sec-dist} (to follow). 
Once all these subdivisions have been made, we are left with stable families.

We note that some of these subfamilies might be empty, depending on the values of $a$ and $d$. For example, if $a$ and $d$ are both divisible by a prime with odd multiplicity, then the $\Delta$-distinguished class is automatically prohibited. To take another example, when $a$ or (resp., and) $d$ is a perfect square, the number of possible places of sufficient-ramification is at least $1$ (resp., $2$).
In total, our families of $(a,d)$-monogenised cubic fields are naturally broken up into as many as nine stable subfamilies, with average values $5/4$, $3/2$, $2$, $2$, and $3$ in the totally real case, and average values $3/2$, $2$, $3$, and $4$, in the complex case. 
More precisely, we prove:
\begin{thm} \label{main stable}
For each of the stable families $F$ defined above, the average $\Avg(F)$ of the sizes of the $2$-torsion subgroups of the class groups of cubic fields in $F$ is bounded as follows:
\begin{equation}\label{eq:SF1}
\begin{array}{ll}
&\displaystyle\Avg(F_{a,d}^+(0))\le\frac{5}{4};\quad \Avg(F_{a,d}^+(1))\le\frac{3}{2};
\quad\Avg(F_{a,d}^+(2;0,*))\le2;
\\[.15in]
&\displaystyle\Avg(F_{a,d}^-(0))\le\frac{3}{2};\quad \Avg(F_{a,d}^-(1))\le\,2;\,\quad
\Avg(F_{a,d}^-(2;0))\le3;
\end{array}
\end{equation}
\begin{equation}\label{eq:SF2}
\begin{array}{lc}
& \displaystyle\Avg(F_{a,d}^+(2;1,0))\le2,\quad \Avg(F_{a,d}^+(2;1,1))\le3;
\\[0.15in]
&\Avg(F_{a,d}^-(2;1))\le4.
\end{array}
\end{equation}
We assume that the complex families are ordered by weighted height and the totally real families are ordered either by weighted height or by balanced height. Furthermore, conditional on the tail estimate of Remark \ref{rmk:unif estimate}, the upper bounds in~\eqref{eq:SF1} and~\eqref{eq:SF2} are in fact equalities.
\end{thm}
\begin{comment}
\[
\limsup_{X \to \infty} \underset{\substack{f \in \pm R^\sigma \Delta_2^{e_2} \Delta_\infty^{e_\infty} \\ \rmH(f) < X}}{\on{Avg}}\,\, \#\on{Cl}(R_f)[2] \hspace{1ex} \le \hspace{1ex} 1 + \frac{1+\sigma}{3 \pm 1} + e_2 \cdot e_\infty
\]

% \[
% \limsup_{X \to \infty} \underset{\substack{f \in - R^\sigma \Delta_2^{e_2} \Delta_\infty^{e_\infty} \\ \rmH(f) < X}}{\on{Avg}}\,\, \#\on{Cl}(R_f)[2] \hspace{1ex} \le \hspace{1ex} \frac{3}{2} + \frac{\sigma}{2} + e_2 \cdot e_\infty
% \]

% \[
% \limsup_{X \to \infty} \underset{\substack{f \in + R^\sigma \Delta_2^{\chi_1} \Delta_\infty^{\chi_2} \\ \rmH(f) < X}}{\on{Avg}}\,\, \#\on{Cl}(R_f)[2] \hspace{1ex} \le \hspace{1ex} \frac{5}{4} + \frac{\sigma}{4} + e_2 \cdot e_\infty
% \]

where $\sigma \in \begin{cases}
\{0,1,3\} \textrm{ if $\gcd(a_k,d_k) = 1$} \\
\{0,1,2\} \textrm{ otherwise}
\end{cases},$ and $e_2,e_\infty \in \{0,1\}$.
\end{comment}

\begin{rem} We make the following remarks in respect of Theorem~\ref{main stable}: 
    \begin{enumerate}
        \item The averages occurring in our previous theorems (Theorems \ref{first main}, \ref{first main2}, \ref{second main real}, and \ref{second main complex}) are obtained by taking a weighted average of the averages for stable families occurring in Theorem~\ref{main stable}, where each stable family is weighted by its density. The discrepancy between the average values for balanced and weighted height is simply due to the difference in the probability of a pair $(b,c)\in\R^2$ with bounded (balanced or weighted) height having $\delta_\infty=1$. This discrepancy vanishes after partitioning into stable families.
        \item The deviation from the Cohen--Lenstra--Martinet philosophy pointed out at the end of \S\ref{sec-unitresult} vanishes across the pairs of stable families $F_{a,d}^{\pm}(0)$, $F_{a,d}^{\pm}(1)$, and $F_{a,d}^+(2;0,*),\, F_{a,d}^-(2; 0)$. On the other hand, this deviation becomes apparent when we partition $F_{a,d}(2;1)$ into stable families.
    \end{enumerate}
\end{rem}

In Table \ref{table:intro}, we present values of the average $2$-torsion in the class group observed in samples of randomly generated $(a,d)$-monogenised cubic fields for several representative choices of the pair $(a,d)$. The table shows that our results are in close agreement with the data.

We now explain the notation in the table. The first row indicates the pair $(\lambda,\delta_{\on{fin}})$, where $\lambda\in\{0,1,2\}$ is the amount of sufficient-ramification in the fields of the family, and $\delta_{\on{fin}}\in\{0,1\}$ indicates the obstruction at finite places to having a $\Delta$-distinguished orbit. In particular, if $\lambda\neq 2$, there is an obstruction at some prime dividing $ad$ to the existence of a $\Delta$-distinguished orbit, and so $\delta_{\on{fin}}$ is automatically $0$. When $\lambda=2$, we have two columns depending on whether the remaining possible obstruction at the primes dividing $2ad$ exists ($(\lambda,\delta_{\on{fin}})=(2,0)$) or not ($(\lambda,\delta_{\on{fin}})=(2,1)$). The second row lists the possibilities for the sign of the discriminant, dividing up our family into totally real and complex cubic fields, and the possible values for the product $\delta_{\on{fin}}\delta_\infty\in\{0,1\}$. Note that $\delta_{\on{fin}} = 0$ unless the corresponding values in the first row is $(2,1)$, and the column corresponding to $(2,1)$ gets split into three possibilities since $\delta_\infty = 1$ for complex cubic fields. Thus the only columns of the table corresponding to cubic fields with a $\Delta$-distinguished orbit are the last two ($(+,1)$ and $(-,1)$). The third row lists the predicted upper bound from Theorem \ref{main stable}. Each column corresponds to a stable family, and the predicted average does not depend on the values of $a$ and $d$, so the predictions are uniform within each column.

The remaining rows correspond to different choices of the set $\{a,d\}$. For a given choice of $\{a,d\}$, the average is taken over a set of several thousand randomly generated cubic fields with an equal number of $(a,d)$-monogenised and $(d,a)$-monogenised cubic fields. Cells are shaded grey when there are no cubic fields with those particular specifications (for example, when $a=d=1$, we always have $\lambda=2$, so the cells corresponding to $\lambda \in \{0,1\}$ are shaded grey).
Table~\ref{table:intro} is obtained by collapsing rows and columns of the more granular Table~\ref{table:app}, which can be found in Appendix \ref{section:appendixdata}, where we provide details on \mbox{how these data were generated.}

\begin{table}[h]
\centering
\begin{tabularx}{0.9\textwidth}{|c||*{2}{X|}|*{2}{X|}|*{2}{X|}|*{3}{X|}}\hline
 $(\lambda,\delta_{\on{fin}})$
& \multicolumn{2}{c||}{$(0,0)$} 
& \multicolumn{2}{c||}{$(1,0)$} 
& \multicolumn{2}{c||}{$(2,0)$} 
& \multicolumn{3}{c|}{$(2,1)$}    \\ \hline 
$({\rm sign}(\Delta),\delta_{\on{fin}}\delta_\infty)$
& $(+,0)$
& $(-,0)$
& $(+,0)$
& $(-,0)$
& $(+,0)$
& $(-,0)$
& $(+,0)$
& $(+,1)$
& $(-,1)$
\\ \hline \hline
% {\bf Granules}
% & $+B_{0,0}^0,+B_{0,1}^0$
% & $-B_{0,1}^0$
% & $+B_{0,0}^1,+B_{0,1}^1$
% & $-B_{0,1}^1$ 
% & $+B_{0,0}^2,+B_{0,1}^2$
% & $-B_{0,1}^2$
% & $+B_{1,0}^2$
% & $+B_{1,1}^2$
% & $-B_{1,1}^2$
% \\ \hline\hline
{\bf Prediction}
& $1.25$
& $1.5$
& $1.5$
& $2$
& $2$
& $3$
& $2$
& $3$
& $4$
\\ \hline\hline
$\{1,1\}$ & \cellcolor{lightgray} & \cellcolor{lightgray} & \cellcolor{lightgray} & \cellcolor{lightgray} & 1.986 & 2.900 & 1.966 & 2.971 & 4.055 \\
\hline
$\{1,2\}$ & \cellcolor{lightgray} & \cellcolor{lightgray} & 1.484 & 2.017 & 2.005 & 2.974 & 2.025 & 3.004 & 4.059 \\
\hline
$\{1,3\}$ & \cellcolor{lightgray} & \cellcolor{lightgray} & 1.502 & 1.992 & 1.992 & 2.999 & 2.002 & 2.995 & 3.996 \\
\hline
$\{1,4\}$ & \cellcolor{lightgray} & \cellcolor{lightgray} & \cellcolor{lightgray} & \cellcolor{lightgray} & 2.023 & 2.964 & 1.976 & 2.990 & 3.990 \\
\hline
$\{1,5\}$ & \cellcolor{lightgray} & \cellcolor{lightgray} & 1.495 & 2.011 & 2.005 & 2.978 & 2.008 & 2.968 & 3.952 \\
\hline
$\{2,2\}$ & 1.254 & 1.492 & 1.491 & 1.944 & \cellcolor{lightgray} & \cellcolor{lightgray} & \cellcolor{lightgray} & \cellcolor{lightgray} & \cellcolor{lightgray} \\
\hline
$\{2,3\}$ & 1.246 & 1.499 & 1.490 & 2.020 & 1.998 & 2.984 & 2.036 & 2.976 & 4.025 \\
\hline
$\{2,4\}$ & \cellcolor{lightgray} & \cellcolor{lightgray} & 1.476 & 1.958 & \cellcolor{lightgray} & \cellcolor{lightgray} & \cellcolor{lightgray} & \cellcolor{lightgray} & \cellcolor{lightgray} \\
\hline
$\{2,5\}$ & 1.246 & 1.503 & 1.502 & 2.003 & 1.990 & 2.957 & 2.034 & 2.979 & 3.936 \\
\hline
$\{3,3\}$ & 1.237 & 1.491 & 1.481 & 2.001 & \cellcolor{lightgray} & \cellcolor{lightgray} & \cellcolor{lightgray} & \cellcolor{lightgray} & \cellcolor{lightgray} \\
\hline
$\{3,4\}$ & \cellcolor{lightgray} & \cellcolor{lightgray} & 1.501 & 1.926 & 1.993 & 2.978 & 2.022 & 2.968 & 3.963 \\
\hline
$\{3,5\}$ & 1.246 & 1.503 & 1.496 & 2.007 & 2.008 & 3.040 & 1.969 & 2.932 & 3.991 \\
\hline
$\{4,4\}$ & \cellcolor{lightgray} & \cellcolor{lightgray} & \cellcolor{lightgray} & \cellcolor{lightgray} & 1.954 & 3.049 & 2.032 & 2.949 & 4.014 \\
\hline
$\{4,5\}$ & \cellcolor{lightgray} & \cellcolor{lightgray} & 1.490 & 1.989 & 2.039 & 2.961 & 1.995 & 2.981 & 4.002 \\
\hline
$\{5,5\}$ & 1.250 & 1.503 & 1.493 & 1.948 & \cellcolor{lightgray} & \cellcolor{lightgray} & \cellcolor{lightgray} & \cellcolor{lightgray} & \cellcolor{lightgray} \\
\hline
\end{tabularx}
\caption{Predicted averages for stable subfamilies of $(a,d)$-monogenised fields, along with averages observed in random samples of fields from each subfamily.}
\label{table:intro}
\end{table}

\begin{rem}
The Cohen--Lenstra--Martinet heuristics correctly predict the average $2$-torsion in the class group in the full family of cubic fields. While the deviation from these predictions for the family of monogenic cubic fields was surprising at the time of its discovery, new work of the second named author and Venkatesh \cite{SiadVenkatesh} explains the deviation for stable families of $n$-monogenic fields through a new theory of algebraic structures on the class group. Briefly, the central concept is that of a \textit{spin structure} on a number field $K$. These are $\sigma \in K^\times$ such that $\delta_K^{-1} = (\sigma) I^2$ for some fractional ideal $I$ of $K$ taken up to the equivalence $\sigma \sim \lambda^2 \sigma$, where $\delta_K$ denotes the different ideal of $K$. They are the arithmetic counterparts of spin structures in the geometry of compact Riemann surfaces and closed oriented $3$-manifolds.   
For any number field $K$, the theory gives a map
from spin structures on $K$ to quadratic refinements of the system of bilinear structures on the dual of $\on{Cl}(K)$.
The (duals) of the class groups of fields in families having canonical family-wide spin structures are then expected to behave like random abelian groups with quadratic structure, rather than random abelian groups with bilinear structure. In theory, stable families of $n$-monogenic fields fall into one of these two regimes. Monogenic fields form the typical example of a family having a family-wide spin structure --- indeed, a monogenic field $K$ whose ring of integers is given by $\Z[\alpha]$, where the monogeniser $\alpha$ has minimal polynomial $f$, has a family-wide spin structure given by $\sigma_{K} = 1/f'(\alpha)$.  
It is a natural question to ask whether these structures account for the averages observed in stable families of $(a,d)$-monogenised fields. 
\end{rem}

\subsection{Results on $2$-Selmer groups of elliptic curves} \label{sec-selrel}

From the perspective of arithmetic invariant theory, the $2$-torsion in the class group of a monogenised cubic field $(K,\alpha)$ is closely related to the $2$-Selmer group of the elliptic curve $E_f \colon y^2 = f(x,1)$, where $f$ is the homogenisation of the minimal polynomial of $\alpha$ (see \S\ref{sec-parametrisations}). In this regard, our methods for proving Theorems~\ref{first main}--\ref{second main complex} may be adapted to obtain similar theorems about the average size of the $2$-Selmer groups in families of elliptic curves of the form $y^2 = x^3 + bx^2 + cx + d$, where $d$ is fixed. Indeed, we prove the following Selmer-group analogues of Theorems~\ref{second main real} and~\ref{second main complex}:

\begin{thm}\label{third main} Let $\rmH \in \{\rmH_{\on{bal}},\rmH_{\on{wei}}\}$. Suppose $d \equiv 1 \pmod 8$, and let $\Sigma \subset U_{1,d}(\Z)_{\max}$ be an acceptable family of forms of positive discriminant such that $\Sigma_2 = \{f \in U_{1,d}(\Z_2) : f(x,y) \equiv x^3 + x^2y + y^3 \pmod8\}$. \mbox{Then we have} 
\[\limsup_{X \to \infty} \underset{\substack{f \in \Sigma \\ \rmH(f) < X}}{\on{Avg}}\,\, \#\on{Sel}_2(E_f) \hspace{1ex}  \le 
3 + 3\Pi_{d,\rmH}\rholeft_\Sigma\]
where $\Pi_{d,\rmH}$ is the real asymptotic volume given in~\eqref{eq-newpid}.
\end{thm}

\begin{thm}\label{third main2} Let $\Sigma \subset U_{1,d}(\Z)_{\max}$ be an acceptable family of forms of negative discriminant such that $\Sigma_2 = \{f \in U_{1,d}(\Z_2) : f(x,y) \equiv x^3 + x^2y + y^3 \pmod8\}$. Then we have 
\[\limsup_{X \to \infty} \underset{\substack{f \in \Sigma \\ \rmH_{\on{wei}}(f) < X}}{\on{Avg}}\,\, \#\on{Sel}_2(E_f) \hspace{1ex}  \le 3+ 3\rholeft_\Sigma.\]
\end{thm}
As was the case for Theorems~\ref{first main}--\ref{second main complex}, we expect that the bounds in Theorems~\ref{third main} and~\ref{third main2} are in fact the exact averages (i.e., the theorems should hold with ``$\limsup$'' replaced by ``$\lim$'' and ``$\leq$'' replaced by ``$=$''); see Remark \ref{rmk:unif estimate}. The condition on $\Sigma_2$ that occurs in the statements of Theorems~\ref{third main} and~\ref{third main2} may be interpreted as fixing the mod-$2$ reduction-type of the elliptic curves in our family. We impose this condition for the sake of simplicity; with further computation, one can use our method of proof to obtain similar results for any family of elliptic curves with discriminant having bounded $2$-adic valuation. We note that our methods do not apply when $d=0$, in which case the average $2$-Selmer size was shown to be infinite in work of Klagsbrun--Lemke Oliver~\cite{MR3375649}. The determination of precise asymptotics for the count of $2$-Selmer elements in the case $d = 0$ is the subject of current work-in-progress of Alp\"{o}ge and Ho (private communication).

Just like with Theorems~\ref{second main real} and~\ref{second main complex}, the formulae given in Theorems~\ref{third main} and~\ref{third main2} may be evaluated explicitly in various cases of interest. When $\Sigma = U_{a,d}(\Z)_{\max}$ the average over all $d \in \Z \smallsetminus \{0\}$ of $\rholeft_\Sigma$ is equal to zero. Thus, on average over all $d \in \Z \smallsetminus \{0\}$, the formulae given in Theorems~\ref{third main} and~\ref{third main2} both converge to $3$. This is, of course, precisely the average determined by Bhargava--Shankar~\cite{MR3272925} for the full family of elliptic curves. On the other hand, when $d$ is a perfect square, the averages given in Theorems~\ref{third main} and~\ref{third main2} are both equal to $6$. This doubling of the average size of the $2$-Selmer group makes sense, as $d = f(0,1)$ being a square forces $E_f$ to have an extra rational point; furthermore, this accords with the work of Bhargava--Ho~\cite{manjulwei}, who obtain an average of $6$ for the family of elliptic curves of the form $y^2 = x^3 + bx^2 + cx + d$, where $d$ runs through nonzero perfect squares (cf.,~the related work of Ananth Shankar~\cite{MR3968769}, which generalises this work to curves of higher genus).
     
The Poonen--Rains heuristics~\cite{MR2833483} make predictions for the distributions of Selmer groups in quite general families of elliptic curves. In particular, these heuristics predict that the average size of the $2$-Selmer group in quite general families of elliptic curves is $3$; this was verified for the universal family in the aforementioned work of Bhargava--Shankar. They also predict an average of $6$ for families with a marked point (as was shown in the aforementioned work of Bhargava--Ho). Moreover, these averages have been shown to be stable under the imposition of sets of local conditions. The averages obtained in Theorems~\ref{third main} and~\ref{third main2} seem to be consistent with these findings: the family $\Sigma$ can be partitioned into two stable subfamilies, one on which the average is $3$ (namely the subfamily of cubics above which the $\Delta$-distinguished orbit mentioned in \S\ref{sec-partitioning} does not exist), and one on which the average is $6$ (namely the subfamily of cubics above \mbox{which the $\Delta$-distinguished orbit exists).}
%On the other hand, the averages obtained in Theorem~\ref{third main} with $d < 0$ are more mysterious: in this case, it is unclear how to partition the family $\Sigma$ into subfamilies on which the average is either $3$ or $6$ (in particular, no local conditions achieve this). This suggests an example of a \emph{global} condition under which a family of elliptic curves behaves differently from the full family of elliptic curves and from the family of elliptic curves with a marked point (or marked nontrivial $2$-Selmer element).

%In the literature, examples of \emph{local} conditions have been found that cause families of elliptic curves to fail the Poonen--Rains heuristics; see, e.g., the work of Bhargava--Elkies--Shnidman~\cite{MR4072495} on a family of Mordell curves $y^2 = x^3 + k$, for which the average size of the $3$-isogeny Selmer group is sensitive to congruence conditions on $k$. To the authors' knowledge, Theorem \ref{third main} in the case $d < 0$ suggests an example of a \emph{global} condition under which a family of elliptic curves fails to satisfy the Poonen--Rains heuristics. % \todo{\small is this oversold?}

\subsection{Discussion of the proof} \label{sec-discuss}

In this subsection, we discuss the structure of the proof of Theorem \ref{first main2}, and we explain the modifications needed to obtain the other main theorems. To prove Theorem \ref{first main2}, we need an asymptotic upper bound for the number of $2$-torsion classes in unit-monogenised cubic fields of bounded weighted height. As is the case for several other results in arithmetic statistics, such asymptotics are derived by first parametrising the arithmetic objects of interest in terms of integral orbits for the action of a reductive group $G$ on a scheme $W$ over $\Z$, and by then counting these orbits, which amounts to counting integer points within a fundamental domain for the action of $G(\Z)$ on $W(\R)$.

In our setting, we use a result of Bhargava, which parametrises the dual of the $2$-torsion in the class group of a maximal cubic order $R_f$ by the set of $\SL_3(\Z)$-orbits on the set of pairs $(A,B)$ of integral ternary quadratic forms with $4\det(Ax-By)=f(x,y)$ (see \S\ref{sec-params}).\footnote{By the ``determinant'' of a ternary quadratic form, we mean the determinant of its Gram matrix.} Using this parametrisation, Bhargava determined the average size of the $2$-torsion subgroup in the class groups of all cubic fields~\cite{MR2183288}. To do so, he developed a general ``averaging method,'' which applies to any $(G,W)$ as above. This method transforms the problem of determining asymptotics for the number of $G(\Z)$-orbits on $W(\Z)$ of bounded height into one of evaluating the integral over elements $\gamma\in G(\Z)\backslash G(\R)$ of the number of integral points lying in the $\gamma$-translates of a fixed bounded region in $W(\R)$.

For our applications, applying Bhargava's averaging method requires us to count pairs of ternary quadratic forms $(A,B)$ lying within regions skewed by elements in $\FF\defeq\SL_3(\Z)\backslash\SL_3(\R)$, satisfying the additional conditions $\det A =\det B =1$. A linearisation trick (used to handle the family of monogenised cubic fields considered in \cite{BHSpreprint}) can be applied to fix the choice of $A$, but that still leaves us with the problem of counting integral points on the degree-$3$ hypersurface cut out by $\det B =1$ in the $6$-dimensional space $V(\R)$ of ternary quadratic forms over $\R$. This problem falls outside the range of applicability of the circle method, which has recently been applied to great effect in arithmetic statistics; see, e.g., the works of Ruth \cite{Ruththesis}, Alp\"{o}ge \cite{2110.03947}, Sanjaya--Wang \cite{Sanjaya2022}, and Alp\"{o}ge--Bhargava--Shnidman \cite{ABSpreprint}.
Since the set of forms $B \in V(\Z)$ with $\det B =1$ is a union of finitely many $\SL_3(\Z)$-orbits of forms $B_i \in V(\Z)$, our problem is reduced to one of counting integral points in regions of each of the finitely many symmetric $\R$-varieties $\SL_3(\R)/\SO_{B_i}(\R)$. A key difficulty that arises when counting integral points in regions of symmetric $\R$-varieties is as follows: the volume of a thickening of the boundary of such a region is of the same order of magnitude as the volume of the region itself. Consequently, some input about the equidistribution of integral points near the boundary is needed to obtain the desired asymptotic point counts.

There are at least three methods to prove these equidistribution results in the literature. The first, by Duke--Rudnick--Sarnak \cite{MR1230289}, uses automorphic forms; the second, by Eskin--McMullen \cite{MR1230290}, builds on the seminal work of Margulis in his thesis (see, e.g.,~\cite{MR2035655}) by using the mixing property of the geodesic flow on Lie groups; and the third, by Eskin--Mozes--Shah \cite{MR1381987}, uses Ratner theory.
The aforementioned work of Eskin--McMullen was subsequently effectivised by Benoist--Oh \cite{MR3025156}, whose work can be used to count the number of unit-determinant integral ternary quadratic forms in $X\cdot \mc{R}$, for sufficiently nice regions $\mc{R}\subset V(\R)$, with a power saving error term.
More precisely, if we denote the set of unit-determinant elements in $X\cdot \mc{R}$ by $\mc{R}_X$, then Benoist--Oh prove that the number of lattice points in $\mc{R}_X$ is equal to $\on{Vol}(\mc{R}_X)$, up to an error of $O(\on{Vol}(\mc{R}_X)^{1-\delta})$, for some positive constant $\delta$. 
For our application, we must additionally understand how the error term in the above asymptotic changes when the regions $\{\mc{R}_X\}_{X > 0}$ is skewed by $\gamma\in\FF$. To this end, we apply the results of Benoist--Oh to deduce that the number of lattice points in $\gamma\mc{R}_X$ is%\todo{\small change made here --- is it okay?}
\begin{equation} \label{eq-skewasymp}
\Vol(\mc{R}_X)+O\left(h(\gamma)^M\Vol(\mc{R}_X)^{1-\delta}\right),
\end{equation}
where $h(\gamma)$ denotes the maximum of the absolute values of the entries of $\gamma$ and $\gamma^{-1}$, and where $M > 0$ is a constant (see Theorem~\ref{main effective counting}). This asymptotic formula~\eqref{eq-skewasymp} serves as a viable substitute for Davenport's Lemma (used in the case when counting in affine space) in our setting as long as $h(\gamma)$
is small --- this corresponds to $\gamma$ lying in the ``main body'' of the fundamental domain $\mc{F}$.

When $h(\gamma)$ is large --- this corresponds to $\gamma$ lying in the ``cusp'' of $\mc{F}$. In this regime, the asymptotic~\eqref{eq-skewasymp} is insufficient. We find an interesting phenomenon in our analysis of the cusp. Following Bhargava's counting techniques in \cite{MR2183288}, we partition the cusp into two pieces, a shallow cusp and a deep cusp. Similar to the situation in \cite{MR2183288}, every relevant integral element in the deep cusp corresponds to the identity element in the class group of the corresponding cubic field. However, the shallow cusp behaves very differently. In particular, we find a positive proportion of non-identity elements in the shallow cusp. We refer to these ``new'' elements as \emph{$\Delta$-distinguished}, and we count them separately by hand. Interestingly, these $\Delta$-distinguished orbits only appear for complex cubic fields and never for totally real fields. In fact, these orbits explain the aforementioned deviation from the Cohen--Lenstra philosophy that was observed in Theorem~\ref{first main2}.
To control the contribution from the non-$\Delta$-distinguished elements in the shallow cusp, we prove a suitable \emph{upper bound} for the number of lattice points in skewed regions in $\gamma \mc{B}_X$ for symmetric $\R$-varieties of the form $S_B(\R)$ (see Theorem~\ref{th:skbound}). Our proof of this upper bound is completely elementary, and it exploits the particular structure of the formula for the determinant of a ternary quadratic form. 
We note that our methods are widely applicable for counting away from the cusp (i.e., in the main body); in particular, they work for any symmetric $\R$-variety, not just the variety of unit-determinant quadratic forms. It is only our handling of the cusp that is specific to our situation.

In \S\ref{sec-countonsymmetricvarieties} and \S\ref{sec-countpairsternaryfixeddet}, we combine Bhargava's averaging method with our adaptation of the Benoist--Oh technique as described above. This allows us to obtain an upper bound on the average size of the $2$-torsion in the class group of cubic fields in our families, as well as the average size of the $2$-Selmer groups of elliptic curves in our families. However, the obtained upper bound is not explicit but instead expressed as a product of local masses. 
In order to compute these masses, it is necessary to compute the joint distribution of the Hasse invariants (at each place) of the pairs of ternary quadratic forms $(A,B)$ that arise in our count. 
In \S\ref{sec-proofofmainresults} and \S\ref{sec-selproofs}, we compute this distribution for the cases of cubic fields and elliptic curves, respectively, thus obtaining our main results.

\subsection*{Acknowledgments}
\noindent Shankar was supported by an NSERC discovery grant and a Sloan fellowship. Siad was supported by a QEII/Steve Halperin GSST from the University of Toronto, a Postdoctoral Fellowship from the NSERC, Princeton University, the Institute for Advanced Study, and NSF Grant No.~DMS-1926686. Swaminathan was supported by the NSF, under the Graduate Research Fellowship as well as Award No.~2202839.

We are very grateful to Iman Setayesh for many helpful conversations, and in particular for helping us understand many of the key ideas behind the work of Eskin--McMullen. It is a pleasure to thank Manjul Bhargava for many helpful discussions and for sharing class group data for many thin families, and sharing his code used to compute the class groups and their averages; We thank Hee Oh immensely for explaining how her joint work with Benoist~\cite{MR3025156} immediately implies Theorem~\ref{main effective counting}, simplifying our previous proof. We are grateful to Manjul Bhargava, Hee Oh, Akshay Venkatesh, and Melanie Matchett Wood for providing us with numerous helpful comments and suggestions on earlier drafts of this paper, as well as for suggesting further lines of inquiry. We thank the anonymous referee for their detailed study of our paper and for providing us with insightful comments and feedback. We also have the pleasure of thanking Levent Alp\"{o}ge, Wei Ho, Pouya Honaryar, Aaron Landesman, Lillian Pierce, Peter Sarnak, James Tao, Yuri Tschinkel, and Jacob Tsimerman for useful discussions.
Finally, we thank the Program in Applied and Computational Mathematics (PACM) for the use of their CPUs to generate our data.

\section{Orbit parametrisations} \label{sec-parametrisations}

Due to work of Bhargava~\cite{MR2183288} (resp., Bhargava--Gross~\cite{MR3156850}), the $2$-torsion subgroups of the class groups of cubic fields (resp., $2$-Selmer groups of elliptic curves) may be parametrised in terms of the integral orbits of the representation $2 \otimes \on{Sym}_2(3)$ of pairs of integral ternary quadratic forms (resp., $2 \otimes \on{Sym}^2(3)$ of pairs of integer matrix ternary quadratic forms) of the group $\on{SL}_3$. In \S\ref{sec-rep}, we define these representations, and in \S\ref{sec-params}--\S\ref{sec-sel}, we recall the parametrisations of $2$-class groups and $2$-Selmer groups. Then, in~\S\ref{sec-dist}, we discuss the notions of distinguished and $\Delta$-distinguished orbits mentioned in~\S\ref{sec-discuss}.

\subsection{The representations $2 \otimes \on{Sym}_2(3)$ and $2 \otimes \on{Sym}^2(3)$ of $\on{SL}_3$} \label{sec-rep}

As in \S\ref{sec-stat}, let $U$ be the affine scheme over $\Z$ whose $R$-points are given by binary cubic forms over $R$ for any ring $R$. The algebraic group $\GL_2$ acts on $U$ via a twisted action $\gamma\cdot f(x,y) \defeq (\det(\gamma)^{-1}f((x,y)\cdot\gamma)$, and the work of Delone--Faddeev \cite{MR0160744} and Gan--Gross--Savin \cite{MR1932327} yields a bijection between $\on{GL}_2(\Z)$-orbits on $U(\Z)$ and isomorphism classes of cubic rings. Denote the cubic ring corresponding to $f\in U(\Z)$ under \mbox{this bijection by $R_f$.}

Next, let $V$ be the affine scheme over $\Z$ whose $R$-points are given by ternary quadratic forms over $R$ for any ring $R$, let $W=V \times_{\Z} V$ denote the space of pairs of ternary quadratic forms, and let $W^\vee \subset W$ be the subscheme of pairs of integer-matrix ternary quadratic forms. We represent elements in $V$ (resp., $W$) by their Gram matrices (resp., pairs of Gram matrices). The algebraic group $\SL_3$ acts on $V$ via $\gamma\cdot A = \gamma A \gamma^T$, and this action induces an action of $\on{SL}_3$ on $W$.

Let $R$ be a ring. Given $(A,B) \in W(R)$, we define the {\it resolvent} of $(A,B)$ to be the binary cubic form $4 \det(xA - yB) \in U(R)$. Let $\on{Res} \colon W \to U$ be the function that takes a pair of ternary quadratic forms to its resolvent; then one readily verifies that the map $\on{Res}$ is $\on{SL}_3$-invariant and that its coefficients generate the ring of polynomial invariants for the action of $\on{SL}_3$ on $W$.

\subsection{Parametrisation of $2$-torsion subgroups of the class groups of cubic fields} \label{sec-params}

Let $K$ be an $S_3$-cubic field (i.e., a cubic field extension of $\Q$ whose normal closure has Galois group equal to the symmetric group $S_3$ over $\Q$; equivalently, a non-Galois cubic field), and let $L$ be an unramified quadratic extension of $K$. Heilbronn~\cite{MR0280461} proved that the Galois closure of $L$ is a degree-$24$ extension of $\Q$ with Galois group $S_4$, thus yielding a quartic field up to conjugacy. Let $\Cl(K)$ denote the class group of $K$, let $\Cl(K)^*$ denote its dual. Nontrivial elements in $\Cl(K)^*[2]$, the $2$-torsion subgroup of $\Cl(K)^*$, correspond to index-$2$ subgroups in $\Cl(K)$. It follows that we may use class field theory to parametrise elements in $\Cl(K)^*[2]$ via quartic fields which are not totally complex. Let $\Cl^+(K)$ denote the narrow class group of $K$, and let $\Cl^+(K)^*$ denote its dual. We may similarly parametrise elements in $\Cl^+(K)^*[2]$ via quartic fields (including the totally complex ones). In combination with Bhargava's parametrisation of quartic fields in~\cite{MR2113024}, this gives the following parametrisation of 
elements in the $2$-torsion of the dual of the class groups and narrow class groups of cubic fields in terms of $\on{SL}_3(\Z)$ orbits on $W(\Z)$ (see \cite[Theorem~2.15]{BHSpreprint} for more details):
\begin{theorem} \label{thm-dual2torsparam}
Let $K$ be an $S_3$-cubic field and let $f \in U(\Z)$ be an $($irreducible$)$ binary cubic form corresponding to the ring of integers of $K$ under the Delone--Faddeev parametrisation. Then
\begin{itemize}
    \item[{\rm (a)}] If $\Delta(K) > 0$, then there is a canonical bijection between the elements of the dual group $\on{Cl}^+(K)[2]^*$ and $\on{SL}_3(\Z)$-orbits of pairs $(A,B) \in W(\Z)$ with $\on{Res}(A,B) = f$.
    
    \noindent Under this bijection, the elements of the dual group $\on{Cl}(K)[2]^* \subset \on{Cl}^+(K)[2]^*$ correspond to $\on{SL}_3(\Z)$-orbits of pairs $(A,B) \in W(\Z)$ such that $A$ and $B$ have a common zero over $\R$.
    \item[{\rm (b)}] If $\Delta(K) < 0$, there is a canonical bijection between the elements of the dual group $\on{Cl}^+(K)[2]^* = \on{Cl}(K)[2]^*$ and $\on{SL}_3(\Z)$-orbits of pairs $(A,B) \in W(\Z)$ with $\on{Res}(A,B)=f$.
\end{itemize}
\end{theorem}

\subsection{Parametrisation of $2$-Selmer groups of elliptic curves} \label{sec-sel}

We now recall how the orbits of $\on{SL}_3(\Z)$ on $W^\vee(\Z)$ parametrise $2$-Selmer elements (i.e., $2$-covers that have rational points over every completion of $\Q$) of elliptic curves.

Let $K = \Q$ or $\Q_v$ for a place $v$ of $\Q$, and let $f \in U(K)$ be monic with nonzero discriminant. Let $\mc{A} \in V(\Z)$ be the anti-diagonal matrix with anti-diagonal entries given by $1$, $-1$, and $1$. We start by explaining how to construct a $2$-cover of the elliptic curve $E_f$ from a pair of the form $(\mc{A},B) \in \on{Res}^{-1}(4f) \cap W^\vee(K)$. Given such a pair $(\mc{A},B)$, we may add a $K$-rational multiple of $\mc{A}$ to $B = (b_{ij})$ to arrange that $b_{13} = 0$. Then it follows from~\cite{MR2899953} and~\cite{MR3800357} that the genus $1$ curve with equation
\begin{equation} \label{eq-associate2cover}
z^2 = \frac{b_{11}}{4} x^4 + b_{12} x^3y + b_{22} x^2y^2 + 2b_{23} x y^3 + b_{33} y^4
\end{equation}
is a $2$-cover of $E_f$. If $(A',B') \in W(K)$ is $\on{SL}_3(K)$-equivalent to a pair of the form $(\cA, B)$ with $b_{13} = 0$, we define the binary quartic form on the right-hand side of~\eqref{eq-associate2cover} to be \emph{binary quartic form associated to the pair $(A',B')$}. The $2$-cover~\eqref{eq-associate2cover} is trivial if and only if the corresponding binary quartic has a root over $K$.

A pair $(\mc{A},B) \in \on{Res}^{-1}(4f) \cap W^\vee(K)$ is said to be \textit{soluble} over $K$ if the associated $2$-cover~\eqref{eq-associate2cover} has a $K$-point. If $K = \Q$, we say that $(\mc{A},B)$ is \textit{locally soluble} if the associated $2$-cover is soluble over $\Q_v$ for every place $v$ of $\Q$. The following parametrisation result is due to Bhargava--Gross (see~\cite{MR3156850}):
\begin{theorem} \label{thm-2selparam}
For a monic binary form $f \in U(\Z)$ of nonzero discriminant, consider the corresponding elliptic curve $E_f\colon y^2=f(x,1)$. Then there is a natural bijection between elements in $\Sel_2(E_f)$ and $\on{SL}_3(\Q)$-equivalence classes on the set of locally soluble pairs $(\cA,B)\in \on{Res}^{-1}(4f) \cap W^\vee(\Z)$.
\end{theorem}

 Furthermore, it follows from \cite[\S8.1]{MR3156850} that if $f$ is maximal, then the locally soluble $\SL_3(\Q)$-equivalence classes on $\Res^{-1}(4f)\cap W^\vee(\Z)$ are the same as the $\SL_3(\Z)$-orbits on $\Res^{-1}(4f)\cap W^\vee(\Z)$. More precisely, every element in $\Res^{-1}(4f)\cap W^\vee(\Z)$ is locally soluble, and no two $\SL_3(\Z)$-orbits on $\Res^{-1}(4f)\cap W^\vee(\Z)$ collapse into the same $\SL_3(\Q)$-equivalence class.

\subsection{Distinguished and $\Delta$-distinguished orbits} \label{sec-dist}

Assume without loss of generality that $a > 0$. We start by defining the notions of distinguished and $\Delta$-distinguished orbits:

\begin{defn} \label{def-disting}
Let $K$ be a field of characteristic not $2$, and let $f \in U(K)$ of nonzero discriminant. We say that an $\on{SL}_3(K)$-orbit on $\on{Res}^{-1}(f) \cap W(K)$ is: \emph{distinguished} if it contains a representative $(A,B)$ with $a_{11} = b_{11} = 0$; and $\Delta$\emph{-distinguished} if it contains a representative $(A,B)$ with $a_{11}a_{22}-a_{12}^2 = b_{11}b_{22} - b_{12}^2 = 0$. Here, we denote by $m_{ij}$ the $ij$-matrix entry of $M \in W(K)$ if $i = j$ and twice the $ij$-matrix entry if $i \neq j$.
\end{defn}

Geometrically, the notion of $\Delta$-distinguished means that the conics cut out by the ternary quadratic forms $A$ and $B$ possess a bitangent line defined over $K$. The following result establishes existence and uniqueness of distinguished and $\Delta$-distinguished orbits over $K$:

\begin{lemma}
Let $K$ be a field of characteristic not $2$. Let $f(x,y)\in U(K)$ be an element with nonzero discriminant. Then there exist unique $\SL_3(K)$-orbits on $\Res^{-1}(K)\cap W(K)$ which are distinguished and $\Delta$-distinguished. The same result holds when $W(K)$ is replaced with $W^\vee(K)$.
\end{lemma}
\begin{proof}
We prove the result for $W(K)$; the proof for $W^\vee(K)$ is identical. Denote the coefficients of $f$ by $a$, $b$, $c$, and $d$, and assume (by replacing $f$ with a $\GL_2(K)$-translate, if necessary) that $ad\neq 0$.
First, it is easy to check that there exists a distinguished orbit (resp.\ a $\Delta$-distinguished orbit) in $W(K)$. The distinguished orbit is represented by the following pair:
\begin{equation*}
(A,B) = \left(\begin{pmatrix}
0 & 0 & 1/2\\
0 & -a & 0\\
1/2 & 0 & -c
\end{pmatrix},
\begin{pmatrix}
0 & 1/2 & 0\\
1/2 & b & 0\\
0 & 0 & d
\end{pmatrix}\right).
\end{equation*}
The $\Delta$-distinguished orbit with resolvent $f$ is represented by the following pair:
\begin{equation} \label{eq-datsdeldist}
(A,B) = \left(\begin{pmatrix}
-a & 0 & 0\\
0 & 0 & 1/2\\
0 & 1/2 & b/4ad
\end{pmatrix},
\begin{pmatrix}
0 & 0 & 1/2\\
0 & d & 0\\
1/2 & 0 & -c/4ad
\end{pmatrix}\right).
\end{equation}

It only remains to prove uniqueness. For the distinguished case, let $L$ be the \'etale cubic extension of $K$ corresponding to $f$ under the Delone--Faddeev bijection. $\SL_3(K)$-orbits on $W(K)\cap\Res^{-1}(f)$ correspond bijectively to \'etale quartic extensions of $K$ with resolvent $f$ by Bhargava's parametrization. An orbit is distinguished if and only if it corresponds to $L\oplus K$, and so there is a unique distinguished orbit.

We give a much more hands-on proof of uniqueness of the $\Delta$-distinguished orbit. Let $(A,B)$ (with coefficients $a_{ij}$ and $b_{ij}$) be an element in a $\Delta$-distinguished orbit with $\Res(A,B)=f$. By the definition of being $\Delta$-distinguished, there exists a $\P^1_K\subset \P^2K$ such that the restriction of $A$ and $B$ to this $\P^1_K$ give reducible binary quadratic forms. (Note that neither of these binary forms can be $0$, otherwise we would have $\Delta(A,B)=0$.) By replacing $(A,B)$ with an $\SL_3(K)$-translate, we may assume that this $\P^1$ is the subspace $(*,*,0)$, and we may hence assume that $a_{12}=a_{22}=b_{11}=b_{12}=0$. Since $a_{11}$ and $b_{22}$ are nonzero, we may use another $\SL_3(K)$-transformation to further assume that $a_{13}=b_{23}=0$. Then we have $a/4=\det(A)=-a_{11}a_{23}^2$ and $b/4=\det(B)=-b_{22}b_{13}^2$. Thus, with another $\SL_3(K)$ transformation, we can ensure that $a_{11}=-a$, $b_{22}=-d$, and $a_{23}=b_{13}=1/2$. Finally, the coefficients $a_{33}$ and $b_{33}$ are uniquely determined by the values of $b$ and $c$. The lemma follows.
\end{proof}

In the next result, we prove that the distinguished orbit over $\Q$ with maximal integral resolvent contains a unique orbit over $\Z$.

\begin{theorem} \label{thm-howmanydist}
Let $f(x,y) = ax^3 + bx^2y + cxy^2 + dy^3 \in U(\Z)_{\max}$. Then the distinguished $\on{SL}_3(\Q)$-orbit on $\on{Res}^{-1}(f) \cap W(\Q)$ and on $\on{Res}^{-1}(4f) \cap W^\vee(\Q)$ contains a unique $\on{SL}_3(\Z)$-orbit.
\end{theorem}
\begin{proof}
Let $K$ be the cubic field corresponding to $f$ under the Delone--Faddev parametrization over fields. Over $\Z$, the form $f$ then corresponds to the maximal order $\O_K$ (since $f$ is maximal by assumption). Bhargava's parametrization over fields associates the algebra $K\oplus\Q$ to the (unique) distinguished orbit with resolvent $f$. Therefore, any $\SL_3(\Z)$-orbit in $\Res^{-1}(f)\cap W(\Z)$ corresponds to an order $R$ of $K\oplus\Q$. However, since $\Delta(K)=\Delta(f)=\Delta(A,B)=\Delta(R)$, it follows that $R$ must be $\O_K\oplus\Z$, the maximal order of $K\oplus\Q$. Since there is only one choice of order, there must be a unique $\SL_3(\Z)$-orbit inside $\Res^{-1}(f)\cap W(\Z)$ as claimed.
\end{proof}

The analogous result for $\Delta$-distinguished orbits is significantly more complicated.     In the following pair of theorems, we  determine necessary and sufficient conditions under which the $\Delta$-distinguished $\on{SL}_3(\Z)$-orbit exists over a maximal binary cubic resolvent (note in particular that this orbit does \emph{not} always exist!), and we establish that this orbit is unique if it exists. For ease of reading, we postpone the proofs to the appendix (see \S\ref{sec-distdensities}). The first theorem handles the case of $\Delta$-distinguished $\on{SL}_3(\Z)$-orbits on $W(\Z)$:

\begin{theorem}\label{th:deltadistex}
Let $f(x,y) = ax^3 + bx^2y + cxy^2 + dy^3 \in U_{a,d}(\Z)_{\max}$. Then the $\Delta$-distinguished $\on{SL}_3(\Q)$-orbit on $\on{Res}^{-1}(f) \cap W(\Q)$ contains an integral representative if the following conditions are satisfied:
\begin{itemize}
    \item $\gcd(a_k,d) = \gcd(a, d_k) = 1$;
    \item $b^2 - 4ac \equiv 0 \pmod{d_k}$ and $c^2 - 4bd \equiv 0 \pmod{a_k}$; and
    \item If $2 \nmid a_kd_k$, then the mod-$4$ residue of $(b,c)$ occurs in the following table, where the left column indicates the mod-$4$ residue of $(a_k, d_k)$, and the top row indicates the mod-$2$ residue of $(a_m, d_m):$
\begin{center}
\begin{tabular}{|l||*{4}{c|}}\hline
\backslashbox{$(a_k,d_k)$}{$(a_m,d_m)$}
&$(0,0)$&$(0,1)$&$(1,0)$
&$(1,1)$\\\hline\hline
$(1,1)$ & $(1,1)$ & $(1,3)$ & $(3,1)$ & $(0,0),\,(1,2),\,(2,1)$\\\hline
$(1,3)$ & $(3,1)$ & $(3,3)$ & $(1,1)$ & $(0,0),\,(2,1),\,(3,2)$ \\\hline
$(3,3)$ & $(3,3)$ & $(3,1)$ & $(1,3)$ & $(0,0),\,(2,3),\,(3,2)$ \\\hline
\end{tabular}
\end{center}

On the other hand, if $2 \mid d_k$, then the mod-$8$ residue of $(b,c)$ occurs in the following table, where the left column indicates the mod-$8$ residue of $(a_k, d_k)$, and the top row indicates the mod-$2$ residue of $(a_m, d_m):$
\begin{center}
\begin{tabular}{|l||*{2}{c|}}\hline
\backslashbox{$(a_k,d_k)$}{$(a_m,d_m)$}
&$(1,0)$&$(1,1)$\\\hline\hline
$(1,2)$ & $(0,1),\,(4,1)$ & $(0,0),\,(2,1),\,(2,4),\,(4,4),\,(6,0),\,(6,1)$ \\\hline
$(1,6)$ & $(0,1),\,(4,1)$ & $(0,0),\,(2,0),\,(2,1),\,(4,4),\,(6,1),\,(6,4)$ \\\hline
$(3,2)$ & $(0,3),\,(4,3)$  & $(0,0),\,(2,3),\,(2,4),\,(4,4),\,(6,0),\,(6,3)$ \\\hline
$(3,6)$ & $(0,3),\,(4,3)$ & $(0,0),\,(2,0),\,(2,3),\,(4,4),\,(6,3),\,(6,4)$ \\\hline
$(5,2)$ & $(0,5),\,(4,5)$  & $(0,0),\,(2,4),\,(2,5),\,(4,4),\,(6,0),\,(6,5)$ \\\hline
$(5,6)$ & $(0,5),\,(4,5)$ & $(0,0),\,(2,0),\,(2,5),\,(4,4),\,(6,4),\,(6,5)$ \\\hline
$(7,2)$ & $(0,7),\,(4,7)$  & $(0,0),\,(2,4),\,(2,7),\,(4,4),\,(6,0),\,(6,7)$ \\\hline
$(7,6)$ & $(0,7),\,(4,7)$ & $(0,0),\,(2,0),\,(2,7),\,(4,4),\,(6,4),\,(6,7)$ \\\hline
\end{tabular}
\end{center}

The entries corresponding to the values of $(a_k, d_k)$ not represented in the above tables may be deduced by symmetry.
\end{itemize}
If such an integral representative exists, it is unique up to $\on{SL}_3(\Z)$-action.
\end{theorem}

  Given an acceptable family $\Sigma \subset U_{a,d}(\Z)_{\max}$, we define $\mathfrak{d}_\Sigma$ to be the density of $f \in U_{a,d}(\Z)_{\max}$, ordered by either balanced or weighted height, such that there is a $\Delta$-distinguished $\on{SL}_3(\Z)$-orbit on $W(\Z)$ with resolvent $f$. The density $\mathfrak{d}_\Sigma$ can be computed in any given case of interest using the criteria given in Theorems~\ref{th:deltadistex}, and it appears in the bounds given in Theorems~\ref{second main real}--\ref{second main complex}.

  The next theorem handles the case of $\Delta$-distinguished $\on{SL}_3(\Z)$-orbits on $W^\vee(\Z)$:

\begin{theorem} \label{th:deltadistex2}
Let $f(x,y) = ax^3 + bx^2y + cxy^2 + dy^3 \in U_{a,d}(\Z)_{\max}$. Then the $\Delta$-distinguished $\on{SL}_3(\Q)$-orbit on $\on{Res}^{-1}(4f) \cap W^\vee(\Q)$ contains an integral representative if and only if the following conditions are satisfied:
\begin{itemize}
    \item $\gcd(a_k,d) = \gcd(a, d_k) = 1$; and
    \item $b^2 - 4ac \equiv 0 \pmod{d_k}$ and $c^2 - 4bd \equiv 0 \pmod{a_k}$.
    \end{itemize}
    If such an integral representative exists, it is unique up to $\on{SL}_3(\Z)$-action.
    \end{theorem}

    It is possible for an $\on{SL}_3(\Z)$-orbit on $W(\Z)$ or $W^\vee(\Z)$ to be both distinguished and $\Delta$-distinguished. However, this happens quite rarely; for instance, we prove the following result:
    \begin{proposition} \label{prop-coincidenot}
    When binary cubic forms $f \in U_{a,d}(\Z)_{\max}$ are ordered by $\on{H}_{\on{bal}}$ or $\on{H}_{\on{wei}}$, the density of $f$ such that the distinguished and $\Delta$-distinguished $\on{SL}_3(\Z)$-orbits on $\on{Res}^{-1}(f) \cap W(\Z)$ $($or $\on{Res}^{-1}(4f) \cap W^\vee(\Z)${}$)$ coincide is zero.
    \end{proposition}
    \begin{proof}
Let $f \in U_{a,d}(\Z)$ be of nonzero discriminant, and suppose that the distinguished and $\Delta$-distinguished $\on{SL}_3(\Z)$-orbits on $\on{Res}^{-1}(f) \cap W(\Z)$ (or $\on{Res}^{-1}(4f) \cap W^\vee(\Z)$) coincide. Then the canonical pair $(A,B)$ in~\eqref{eq-datsdeldist} must have a common root with coordinates in $\Q$ (the case of $W(\Z)$ is related to the case of $W^\vee(\Z)$ by acting on the pair in~\eqref{eq-datsdeldist} via the diagonal matrix with diagonal entries $(1, 1, 2)$). This property remains true for the pair $(A',B') \defeq \gamma \cdot (A,B)$, where $\gamma$ is the diagonal matrix with diagonal entries $1/\sqrt{a}$, $1$, and $1$. Note that $\on{Res}(A',B')$ is monic, so there is a binary quartic form $g$ associated to the pair $(A',B')$ (see \S\ref{sec-sel}), and this form has a root over $\Q(\sqrt{a})$. Explicitly, the form $g$ is given as follows:
\begin{equation} \label{eq-binquartdel}
    g(x,y) = dx^4 - \frac{b}{2a}x^2y^2 + \frac{1}{\sqrt{a}}xy^3 +  \left( \frac{b^2}{16a^2d} - \frac{c}{4ad} \right)y^4.
\end{equation}
We claim that, for each pair of fixed values for $a,d \in \Z \smallsetminus \{0\}$, the binary quartic form in~\eqref{eq-binquartdel} is irreducible over the function field $\Q(\sqrt{a})(b,c)$. To prove this claim, it suffices to take $b = 0$ and $y = 1$ and show that, for a suitable choice of $c$, the resulting quartic polynomial
\begin{equation} \label{eq-binquartdel0}
    dx^4 + \frac{1}{\sqrt{a}}x - \frac{c}{4ad}
\end{equation}
has no real roots. But we can do this by choosing $c$ to be sufficiently large relative to $a$ and $d$ and to have the correct sign, so as to ensure that the constant term in~\eqref{eq-binquartdel0} pushes the quartic polynomial completely above or below the horizontal axis. Having proven the claim, the proposition then follows from Cohen's effective form of Hilbert's Irreducibility Theorem; see~\cite[Theorem~2.5]{MR628276}.
    \end{proof}

We now describe $\Delta$-distinguished orbits over $\R$:

\begin{proposition} \label{prop-zeroset}
    Let $f \in U_{a,d}(\R)$ be of nonzero discriminant. Then any pair $(A,B)$ lying in the $\Delta$-distinguished $\on{SL}_3(\R)$-orbit on $\on{Res}^{-1}(f) \cap W(\R)$ $($or $\on{Res}^{-1}(4f) \cap W^\vee(\R)${}$)$ is such that the zero sets of $A$ and $B$ have:
    \begin{itemize}
        \item two $\R$-points of intersection if $\on{Disc}(f) < 0$;
        \item zero $\R$-points of intersection if $\on{Disc}(f) > 0$, and $bd < 0$ or $c < 0$; and
        \item four $\R$-points of intersection if $\on{Disc}(f) > 0$, $bd > 0$, and $c > 0$.
    \end{itemize}
\end{proposition}
\begin{proof}
The number of $\R$-points of intersection of the zero sets of $A$ and $B$ is equal to the number of real roots of the binary quartic form $g$ in~\eqref{eq-binquartdel} (note that $g$ is in fact defined over $\R$ because we chose $a > 0$). Note that the discriminant of $g$ has the same sign as $\on{Disc}(f)$. Thus, if $\on{Disc}(f) < 0$, then $g$ must have two real roots. On the other hand, if $\on{Disc}(f) > 0$, then $g$ can have zero or four real roots. A simple criterion is given in~\cite[\S1]{MR4349740} for when a binary quartic form with positive discriminant has zero or four real roots. Applying this criterion to the form $g$, we see that $g$ has zero real roots if one of $b$ or $c$ is negative.
    \end{proof}
The next result determines the conditional probabilities with which the last two outcomes itemised in Proposition~\ref{prop-zeroset} occur, with respect to our two different height functions:
\begin{proposition} \label{prop-realvols}
    We have that
    \begin{align*}
        & \mathfrak{d}_\infty^{(2)}  \defeq \lim_{X \to \infty} \frac{\on{Vol}(\{f \in U_{a,d}(\R) : \on{H}(f) < X,\, \on{Disc}(f) > 0,\, bd < 0 \text{ or } c < 0\})}{\on{Vol}(\{f \in U_{a,d}(\R) : \on{H}(f) < X,\, \on{Disc}(f) > 0\})} = \begin{cases} 3/4, & \text{if $\on{H} = \on{H}_{\on{bal}}$,} \\ 1 - \frac{1}{24a+2}, &  \text{if $\on{H} = \on{H}_{\on{wei}}$} \end{cases}\\
        & \mathfrak{d}_\infty^{(0)} \defeq \lim_{X \to \infty} \frac{\on{Vol}(\{f \in U_{a,d}(\R) : \on{H}(f) < X,\, \on{Disc}(f) > 0,\, bd > 0 \text{ and } c > 0\})}{\on{Vol}(\{f \in U_{a,d}(\R) : \on{H}(f) < X,\, \on{Disc}(f) > 0\})} = \begin{cases} 1/4, & \text{if $\on{H} = \on{H}_{\on{bal}}$,} \\ \frac{1}{24a+2}, & \text{if $\on{H} = \on{H}_{\on{wei}}$} \end{cases}
    \end{align*}
\end{proposition}
\begin{proof}
    In the case of balanced height, binary cubic forms $f \in U_{a,d}(\R)$ have positive discriminant with probability $1 - o(1)$, so the probabilities reduce to the probability that $bd < 0$ or $c < 0$ (resp., $bd > 0$ and $c > 0$), which is $3/4$ (resp., $1/4$). In the case of weighted height, the condition $\on{Disc}(f) > 0$ is equivalent, with probability $1 - o(1)$, to stipulating that $b^2 - 4ac > 0$; indeed, the two dominant terms of the discriminant are $b^2c^2 - 4ac^3$.  The condition $b^2 - 4ac > 0$ happens with probability $1/2 + 1/24a$. This happens together with the condition $bd < 0$ or $c < 0$ (resp., $bd > $ and $c > 0$) with probability $1/2 + 1/48a$ (resp., $1/48a$).
\end{proof}
    
    The final result of this section, which describes the $\R$-solubility of the $\Delta$-distinguished orbit, is a corollary of (the proof of) Proposition~\ref{prop-zeroset}:
    \begin{corollary} \label{cor-rsol}
Let $f \in U_{1,d}(\R)$. Then the $\Delta$-distinguished $\on{SL}_3(\R)$-orbit on $\on{Res}^{-1}(4f) \cap W^\vee(\R)$ is $\R$-soluble if and only if one of the following three conditions is satisfied:
\begin{itemize}
\item $\on{Disc}(f) < 0$; or 
\item $\on{Disc}(f) > 0$, and $d > 0$; or
\item $\on{Disc}(f) > 0$, $b < 0$, $c > 0$, and $d < 0$.
\end{itemize}
    \end{corollary}
    \begin{proof}
        The associated binary quartic form $g$~\eqref{eq-binquartdel} has two real roots precisely when $\on{Disc}(f) < 0$. On the other hand, when $\on{Disc}(f) > 0$ and $d > 0$, the quartic $g$ is either positive-definite or has a real root (in which case it must have four real roots). Finally, when $\on{Disc}(f) > 0$ and $d < 0$, we can ensure that the quartic $g$ is not negative-definite by taking $b< 0$ and $c > 0$.
    \end{proof}

\section{Counting integral points on symmetric varieties} \label{sec-countonsymmetricvarieties}

In this section, we start by defining the notion of a symmetric variety (see \S\ref{sec-notset}). We then obtain an asymptotic formula with effective error term for the number of lattice points in skews of families of homogeneously expanding regions in certain symmetric varieties (see \S\ref{sec-maincounter}). We also prove a version of this formula for lattice points satisfying sets of congruence conditions (see \S\ref{sec-congrconds}).

\subsection{Notation and setup} \label{sec-notset}

Let $G$ be a connected almost simple algebraic group defined over $\Z$, let $V$ be a finite-dimensional representation of $G$, defined over $\Z$, with finite kernel, and fix a basis of $V$. A sub-algebraic-group $H \subset G$ defined over $\Z$ is said to be a \emph{symmetric $\Z$-subgroup} if the identity component of $H$ is the identity component of the fixed points of some involution $\sigma$ of $G$ defined over $\Z$. Let $v\in V(\Z)$ be such that the stabiliser of $v$ in $G$ is a symmetric $\Z$-subgroup $H_v \subset G$. Then $\mc{S}_v \defeq G/H_v \simeq Gv$ is said to be a \emph{symmetric variety}; note that $\mc{S}_v$ is in fact a subscheme of $V$.

Let $S_v$ be the functor that assigns to a ring $R$ the quotient $S_v(R) \defeq G(R)/H_v(R)$. Then $S_v(R) = G(R)/H_v(R) \simeq  G(R)v$ for any $R$; when $R = K$ is a field, $S_v(K)$ is said to be a {\it symmetric $K$-variety}. Note that
the functor $S_v$ is \emph{not} to be confused with the functor of points of $\mc{S}_v$. Indeed, regard $G$, $H_v$, and $\mc{S}_v$ as sheaves on the category of $R$-algebras via their functors of points. Then the short exact sequence 
\begin{equation} \label{eq-exactgroups}
0 \to H_v \overset{\iota}\longrightarrow G \overset{\pi}\longrightarrow \mc{S}_v \to 0
\end{equation}
of sheaves induces the following short exact sequence in fppf cohomology:
\begin{equation} \label{eq-galois}
    0 \to S_v(R) \to \mc{S}_v(R) \to \ker\big(H_{\on{fppf}}^1(\on{Spec} R, H_v) \to H_{\on{fppf}}^1(\on{Spec} R, G)\big) \to 0.
\end{equation}
We see that the set $\ker\big(H_{\on{fppf}}^1(\on{Spec} R, H_v) \to H_{\on{fppf}}^1(\on{Spec} R, G)\big)$ is in bijection with the set of orbits of the action of $G(R)$ on $\mc{S}_v(R)$, which is not necessarily transitive, and the subset $S_v(R) \subset \mc{S}_v(R)$ constitutes just \emph{one} of these orbits, namely the orbit containing $v$. It is a remarkable fact, due to Borel and Harish-Chandra~\cite{MR147566}, that when $R = \R$, and even when $R = \Z$, the set of $G(R)$-orbits on $\mc{S}_v(R)$ is finite.

Let $\mathrm{d}g$ and $\mathrm{d}h$ respectively denote \emph{volume forms} (i.e., generators of the $1$-dimensional vector space of left-invariant top-degree differential forms) on $G$ and $H_v$ defined over $\Q$. This naturally yields a $G$-invariant quotient measure $\mathrm{d}s$ on $\mc{S}_v$; in what follows, we often restrict this quotient measure to the subset $S_v(R) \simeq G(R)v$, where $R = \R$ or $\Z_p$ for a prime $p$. We write $\tau_{G,\infty} \defeq \on{Vol}(G(\R)/G(\Z))$ and $\tau_{v,\infty} \defeq \on{Vol}(H_v(\R)/H_v(\Z))$, and we write $\tau_{G,p} \defeq \on{Vol}(G(\Z_p))$ and $\tau_{v,p} \defeq \on{Vol}(H_v(\Z_p))$.

\subsection{Main counting results} \label{sec-maincounter}

%A \emph{lattice} $\Gamma \subset G(\R)$ is a subgroup such that the quotient space $G(\R)/\Gamma$ has finite volume under Haar measure. 
The main theorem of this section is the following:

\begin{theorem}[\protect{\cite{MR3025156}}]\label{main effective counting}
 Let $\Gamma\subset G(\Z)$ be a finite-index subgroup, denote $\Gamma\cap H_v(\R)$ by $\Gamma_{H_v}$, and assume $\on{Vol}(H_v(\R)/\Gamma_{H_v}) < \infty$. Let $\cB\subset V(\R)$ be a nonempty open bounded set with smooth boundary. Let $X\cB$ be dilation of $\cB$ by $X > 0$, and write $\cB_X(v) = X\cB \cap G(\R)v$. Then for any $\gamma \in G(\R)$, we have %\todo{\small removed commensurate}
\begin{equation*}
\#(\Gamma v\cap \gamma \cB_X(v) )=\frac{\Vol(H_v(\R)/\Gamma_{H_v})}{\Vol( G(\R)/\Gamma)}\Vol\big(\cB_X(v)\big)+E_v(X,\gamma),
\end{equation*}
where $E_v(X,\gamma) \defeq O(\Vol(\cB_X(v))^{1-\delta}h(\gamma)^M)$ for some $\delta,M>0$. Here, $h(\gamma)$ denotes the maximum of the absolute values of the entries of $\gamma$ and $\gamma^{-1}$, considered as elements of $\GL(V)$, and the implied constant is independent of $\gamma$ and $X$.\footnote{As communicated to us by Oh, a suitable upgrade of this theorem to semisimple groups $G$ can also be deduced from her work with Benoist along similar lines and with arguments well-known to experts.}
\end{theorem} 
\begin{proof}

The theorem follows directly from the work of Benoist--Oh \cite{MR3025156}, as we now explain. 
%The theorem was proven in the case $\gamma$ by Benoist--Oh \cite{MR3025156}, who effectivised the earlier work of Eskin--McMullen \cite{MR1230290}. As communicated to us by Oh, their proof can be readily applied to handle the case of arbitrary $\gamma$, and in particular to obtain the desired effective dependence on $\gamma$ in the error term. In what follows, we briefly explain how to do this.
We first recall the necessary notation from \cite{MR3025156}. Set $\cB_X \defeq \cB_X(v)$ for the sake of brevity. Let $\varepsilon > 0$, and let $B_\varepsilon \subset G(\R)$ be the closed ball of radius $\varepsilon$ centered at $1 \in G(\R)$. Let $(\gamma\cB_X)_\varepsilon^+ \defeq B_{\varepsilon}\gamma \cB_X$ and $(\gamma\cB_X)_\varepsilon^- \defeq \bigcap_{g \in B_{\varepsilon}} g\gamma\cB_X$. %By~\cite[Lemma~10.8]{MR3025156}, %and a smooth function $\alpha_\varepsilon \colon G(\R) \to \R$ supported in $B_\varepsilon$ such that $\int_{g \in G(\R)} \alpha_\varepsilon(g) \mathrm{d}g = 1$ and $S_m(\alpha_\varepsilon) \le \varepsilon^{-p}$; here $S_m$ denotes the Sobolev norm (see \cite[\S11]{MR3025156} for the definition).
Then it is shown in \cite[Lemma 10.8 and proof of Theorem 12.2]{MR3025156} that there exists $p > 0$ such that, for all sufficiently small $\varepsilon$,  
the desired asymptotic holds up to an error bounded by 
\begin{equation} \label{eq-initbound}
\varepsilon^{-p} \int_{s \in (\gamma \mc{B}_X)_{\varepsilon}^+} \lVert s \rVert^{-k} \mathrm{d}s + \Vol\big((\gamma\cB_X)_\varepsilon^+ \smallsetminus (\gamma\cB_X)_\varepsilon^-\big),
\end{equation}
for any $k > 0$; here, $||-||$ denotes any fixed Euclidean norm on $V(\R)$.

Now, the family $\cB_X$ is effectively well-rounded (see~\cite[Definition~12.1]{MR3025156}). Indeed, effective well-roundedness is verified in~\cite{MR3025156} in the special case where $\mc{B}$ is a Euclidean ball centered at the origin; the general case can be proven by adapting this argument and applying the asymptotic expansion for $\Vol(\cB_X)$ provided by \cite[Theorem 6.4]{MR2488484}. As a consequence of effective well-roundedness, we have that
$$\int_{s \in (\gamma \cB_X)_\varepsilon^+} \lVert s \rVert^{-k} \mathrm{d}s \le \int_{s \in \gamma(\cB_X)_{\varepsilon'}^+} \lVert s \rVert^{-k} \mathrm{d}s \le \lVert \gamma^{-1} \rVert^{k} \int_{s \in (\cB_X)_{\varepsilon'}^+} \lVert s \rVert^{-k} \mathrm{d}s \ll h(\gamma)^{M'}\on{Vol}(\mc{B}_X)^{1- \delta'},$$
for some $M',\delta' > 0$; here, $\varepsilon' \ll h(\gamma)^{M''}\varepsilon$ is chosen so that $(\gamma \mc{B}_X)_{\varepsilon}^+ \subset \gamma(\mc{B}_X)_{\varepsilon'}^+$ and $(\gamma \mc{B}_X)_{\varepsilon}^- \supset \gamma(\mc{B}_X)_{\varepsilon'}^-$. 
%All that remains is to understand how the error term in their counting theorem depends on $\gamma$. This error term takes the form: 
%\[cS_m(\alpha_{\varepsilon})\int_{s \in (\cB_X)_\varepsilon^+} \lVert \gamma s \rVert^{-k} \mathrm{d}s + \Vol\big((\gamma\cB_X)_\varepsilon^+\big) - \Vol\big((\gamma\cB_X)_\varepsilon^-\big),\]
%where $(\gamma\cB_X)_\varepsilon^+ \defeq B_{\varepsilon}\gamma \cB_X$ and $(\gamma\cB_X)_\varepsilon^- \defeq \bigcap_{g \in B_{\varepsilon}} g\gamma\cB_X$. 
This gives a bound on the first term in~\eqref{eq-initbound}; as for the second term, we have \((\gamma\cB_X)_\varepsilon^+ \smallsetminus (\gamma\cB_X)_\varepsilon^- \subset \gamma \left( (\cB_X)_{\varepsilon'}^+ \smallsetminus (\cB_X)_{\varepsilon'}^- \right)\), so it follows that
$$\Vol\big((\gamma\cB_X)_\varepsilon^+ \smallsetminus (\gamma\cB_X)_\varepsilon^-\big) \ll \on{Vol}\big((\cB_X)_{\varepsilon'}^+ \smallsetminus (\cB_X)_{\varepsilon'}^-\big).$$
As another consequence of effective well-roundedness, we have for some $\kappa > 0$ that
$$\on{Vol}\big((\cB_X)_{\varepsilon'}^+ \smallsetminus (\cB_X)_{\varepsilon'}^-\big) \ll {\varepsilon'}^\kappa \on{Vol}(\cB_X).$$
Taking $\varepsilon = h(\gamma)^{-\frac{\kappa (M''-M')}{p+\kappa}}\on{Vol}(\mc{B}_X)^{-\frac{\delta'}{p + k}}$ then yields the theorem.
\end{proof}

Take $G = \SL_n$; take $V = \Sym^2(n)$ to be the space of $n \times n$ symmetric matrices with the action given by $g \cdot A' = g A' g^T$. So $V(\Z)$ consists of integral $3\times 3$ symmetric matrices; take $A \in V(\Z)$ of nonzero determinant; and take $H = \on{SO}_A$, the orthogonal group scheme associated to $A$. As is well-known (see, e.g., \cite[Example (B)]{MR3025156}), these choices fit into the setup outlined in \S\ref{sec-notset}. Indeed, $\SL_n$ is a connected almost simple algebraic group over $\Z$, and $H_A = \SO_{A}$ is the fixed points of the involution $\sigma(g) = A (g^{-1})^{T} A^{-1}$ on $\SL_n$. 

Let $a = \det A$, and denote by $V_a \subset V$ the subvariety of matrices with determinant $a$. The following lemma shows that the subscheme $\mc{S}_A \subset V$ is precisely $V_a$:

\begin{lemma}
The map of $\Z$-schemes $\on{SL}_3 \to V_a$ defined by $g \mapsto gAg^T$ induces an isomorphism of $\Z$-schemes $\mc{S} \simeq V_a$.
\end{lemma}
\begin{proof}
As $\mc{S}_A$ and $V_a$ are flat over $\Z$, it suffices to prove the desired isomorphism after basechanging to $K$, where $K = \Q$ or $\F_p$ for a prime $p$. Fix an algebraic closure $\ol{K}$ of $K$. As the map $\on{SL}_3 \to V_a$ is defined by $\on{SO}_A$-invariant polynomial functions, it induces a map $\mc{S}_A \to V_a$. We claim that this induced map is a bijection on $\ol{K}$-points. Indeed, the map $\on{SL}_3(\ol{K}) \to V_a(\ol{K})$ is surjective since $\on{SL}_3(\ol{K})$ acts transitively on $V_a(\ol{K})$, and two elements $g_1, g_2 \in \on{SL}_3(\ol{K})$ have the same image if and only if $g_1^{-1}g_2 \in  \on{SO}_A(\ol{K})$, so we obtain a bijection $S_A(\ol{K}) \to V_a(\ol{K})$. But by the exact sequence~\eqref{eq-galois}, we have $\mc{S}_A(\ol{K}) = S_A(\ol{K})$, as claimed. It now follows from Zariski's Main Theorem~\cite[Th\'{e}or\`{e}me~4.4.3]{MR163911} that the map $\mc{S}_A \to V_a$ defines an isomorphism of $\ol{K}$-varieties; then, the theory of descent~\cite[Th\'{e}or\`{e}me~2.7.1(viii)]{MR199181} implies that the map $\mc{S}_A \to V_a$ is also an isomorphism at the level of $K$-varieties.
\end{proof}

By the result of Borel and Harish-Chandra referenced in \S\ref{sec-notset}, obtaining asymptotics for the number of integral points of $\mc{S}_A(\Z) = V_a(\Z)$ that lie in skewed expanding domains of $V(\R)$ reduces to obtaining asymptotics for $\#\{ \SL_n(\Z) A\cap \gamma X\cB \}$ for arbitrary elements $A\in V_a(\Z)$, which is precisely the problem solved in Theorem~\ref{main effective counting}.
We thus obtain the following corollary to the theorem:
\begin{corollary}\label{cor:nqf}
For any $a \in \Z \smallsetminus \{0\}$, we have
\begin{equation*}
\displaystyle\#\big(V_a(\Z) \cap \gamma X\cB\big)=\displaystyle
\sum_{A\in\frac{V_a(\Z)}{\SL_n(\Z)}}
\frac{\tau_{A,\infty}}{\tau_{n,\infty}}\Vol\big(\cB_X(A)\big) + E_{A}(X,\gamma).
\end{equation*}
\end{corollary}

\subsection{Congruence conditions} \label{sec-congrconds}

Assume now that $G$ is smooth over $\Z$. We next prove a version of Theorem~\ref{main effective counting} for counting lattice points defined by finitely many congruence conditions. Given a $\Z$-scheme $S$ equipped with a nowhere vanishing top-degree differential form defined over $\Q$, a function $\phi \colon S(\Z) \to [0,1]$ is said to be \emph{defined by congruence conditions} if, for every prime $p$, there exists a corresponding function $\phi_p \colon S(\Z_p) \to [0,1]$ such that:
\begin{enumerate}
    \item For every $s \in S(\Z)$, the product $\prod_p \phi_p(s)$ converges and is equal to $\phi(s)$;
    \item For every $p$, the function $\phi_p$ is locally constant outside a closed subset of $S(\Z_p)$ of volume zero with respect to the chosen differential form.
\end{enumerate}
If $\phi_p \equiv 1$ for all but finitely many $p$, then we say that $\phi$ is \emph{defined by finitely many congruence conditions}. 
We then have the following result:

\begin{proposition} \label{prop-singcong}
Let $\phi \colon \mc{S}_v(\Z) \to [0,1]$ be a function defined by finitely many congruence conditions. Then we have
\begin{equation*}
\sum_{v' \in \mc{S}_v(\Z) \cap \gamma X\cB}\phi(v') =
\sum_{v'\in\frac{\mc{S}_v(\Z)}{G(\Z)}}
\frac{\tau_{H_{v'},\infty}}{\tau_{G,\infty}}\Vol\big(\cB_X(v')\big) \Vol\big(S_{v'}(\widehat{\Z})\big)^{-1}\int_{s \in S_{v'}(\widehat{\Z})} \phi(s) \mathrm{d}s
+E_{v'}(X,\gamma).
\end{equation*}
\end{proposition}
\begin{proof}
We assume that $\phi$ is defined modulo $m$ for some integer $m>0$, and we define $\Gamma\subset G(\Z)$ to be the congruence subgroup of matrices reducing to the identity modulo $m$. In particular, the action of $\Gamma$ on $\mc{S}_v(\Z)$ preserves $\phi$. Given $v' \in V(\Z)$ and $g \in G(\Z)$, we can choose the measures on $H_v$ and $H_{g v}$ to correspond to each other via the map $H_v \to H_{g v}$ defined by conjugation $h \mapsto g h g^{-1}$. 
This implies in particular that $\tau_{v,\nu}=\tau_{\gamma v,\nu}$ for all places $\nu$ of $\Q$. Then an application of Theorem~\ref{main effective counting} yields that
\begin{align*}
\sum_{v'\in \mc{S}_v(\Z) \cap \gamma X\cB}\phi(v')
&=
\sum_{v'\in \frac{\mc{S}_v(\Z)}{\Gamma}}\phi(v')\times \#(\Gamma v' \cap \gamma X)
\\ &= \sum_{v'\in \frac{\mc{S}_v(\Z)}{\Gamma}}\phi(v')
\frac{\Vol(H_{v'}(\R)/\Gamma_{H_{v'}})}{\Vol(G(\R)/\Gamma)}\Vol\big(\cB_X(v')\big)+E_{v'}(X,\gamma)
\\ 
&= \sum_{v'\in\frac{\mc{S}_v(\Z)}{G(\Z)}}
\frac{\tau_{H_{v'},\infty}}{\tau_{G,\infty}}\Vol\big(\cB_X(v')\big)
\sum_{g\in G(\Z)/\langle \Gamma , H_{v'}(\Z) \rangle}
\phi(g v')\frac{[H_{g v'}(\Z):\Gamma_{H_{g v'}}]}{[G(\Z):\Gamma]}
+E_{v'}(X,\gamma).
\end{align*}
Next, since $\phi$ is defined modulo $m$, we have
\begin{equation*}
\Vol\big(S_{v'}(\widehat{\Z})\big)^{-1}\int_{s\in S_{v'}(\widehat{\Z})}\phi(s)\mathrm{d}s=\big(\#S_{v'}(\Z/m\Z)\big)^{-1}\sum_{s\in S_{v'}(\Z/m\Z)}\phi(s).
\end{equation*}
Writing $S_{v'}(\Z/m\Z)$ as $G(\Z/m\Z)/H_{v'}(\Z/m\Z)$, we see that it now suffices to prove that
\begin{equation}\label{eq:tempcong}
\sum_{g \in G(\Z)/\langle \Gamma , H_{v'}(\Z) \rangle} \phi(g v')[H_{g v'}(\Z): \Gamma_{H_{g v'}}]=
\sum_{g\in G(\Z/m\Z)/H_{v'}(\Z/m\Z)}\phi(g v') \times \#H_{g v'}(\Z/m\Z).
\end{equation}
but this equality follows since $[H_{g v'}(\Z):\Gamma_{H_{g v'}}]$ and $\#H_{g v'}(\Z/m\Z)$ are independent of $g$ and since we have $G(\Z)/\langle \Gamma , H_{v'}(\Z) \rangle=G(\Z/m\Z)/(H_{v'}(\Z)/\Gamma_{H_{v'}})$.
\end{proof}

We now specialise to the case when $G=\SL_n$ and $V=\Sym^2(n)$ and prove the following analogue of Corollary~\ref{cor:nqf} for counting symmetric matrices of fixed determinant that satisfy finitely many congruence conditions:
\begin{corollary}
Let $a \in \Z \smallsetminus \{0\}$, and let $\phi \colon V_a(\Z) \to [0,1]$ be a function defined by finitely many congruence conditions. Then we have
\begin{equation*}
\sum_{A \in V_a(\Z) \cap \gamma X\cB}\phi(A) =
\sum_{A\in\frac{V_a(\Z)}{\SL_n(\Z)}}
\frac{\tau_{A,\infty}}{\tau_{\on{SL}_n,\infty}}\Vol\big(\cB_X(A)\big) \Vol\big(S_{A}(\widehat{\Z})\big)^{-1}\int_{s \in S_{A}(\widehat{\Z})} \phi(s) \mathrm{d}s
+E_{A}(X,\gamma).
\end{equation*}    
\end{corollary}

\section{Counting integral orbits on pairs of ternary quadratic forms with fixed determinants} \label{sec-countpairsternaryfixeddet}

Fix $a,d \in \Z \smallsetminus \{0\}$, and let $W_{a,d} \subset W$ denote the subvariety of pairs $(A,B)$ satisfying  
$\det A =a/4$ and $\det B = -d/4$. 
Recall from \S\ref{sec-stat} that we define two heights on $U_{a,d}(\R)$, the balanced height and the weighted height. We lift these heights to $W_{a,b}(\R)$ by setting $\rmH_\circ(A,B) \defeq \rmH_\circ(\Res(A,B))$ for each $\circ \in \{\on{bal},\on{wei}\}$.  
For a real number $X>0$ and a subset $T\subset U_{a,d}(\R)$ or $T\subset W_{a,d}(\R)$, we denote the set of elements in $T$ with balanced height (resp., weighted height) less than $X$ by $T^\bx_X$ (resp., $T^\wei_X$). 

In this section, our goal is to determine asymptotics for the number of non-distinguished $\SL_3(\Z)$-orbits on $W_{a,d}(\Z)$ having bounded balanced (resp., weighted) height. To state the result precisely, we first need to define two notions, namely splitting types over $\R$ and local masses.

\medskip
\noindent \textbf{Splitting types over $\R$.}
The analogue of an $\SL_3(\Z_p)$-invariant weight function at the infinite place is given by the notion of a  splitting type.
For a set $T\subset U(\R)$, let $T^{\pm} \defeq \{f\in T : \pm\on{Disc}(f)>0\}$. 
\begin{defn}
We say that a pair $(A,B)\in W(\R)$ with resolvent of nonzero discriminant has {\it splitting type} $(1111)$, $(112)$, or $(22)$ if the zero loci in $\P^2(\R)$ of the quadratic forms $A$ and $B$  intersect in four real points, two real points and one set of complex conjugate points, or two sets of complex conjugate points respectively.
\end{defn}
\noindent Note that the splitting type of $(A,B) \in W(\R)$ remains unchanged under the action of $\on{SL}_3(\R)$. From \cite[Lemmas 16,\,17]{BHSpreprint}, we obtain the following useful description of the splitting types of the orbits of $\on{SL}_3(\R)$ on $W(\R)$ having resolvent $f^{\pm} \in U(\R)^{\pm}$:
\begin{lemma}
We have the following two points:
\begin{itemize}
\item[{\rm (1)}] The set of elements in $W(\R)$ having resolvent $f^+$ consists of three $\SL_3(\R)$-orbits, one with splitting type $(1111)$ and three with splitting type $(22)$. The stabiliser in $\SL_3(\R)$ of each element in any of these orbits has size $\sigma^+=4$.
\item[{\rm (2)}] The set of elements in $W(\R)$ having resolvent $f^-$ consists of a single $\SL_3(\R)$-orbit with splitting type $(112)$. The stabiliser in $\SL_3(\R)$ of each element in this orbit has size $\sigma^-=2$.
\end{itemize}
\end{lemma}
Moreover, from \cite[Remark 3.18]{BHSpreprint}, we may canonically distinguish between the three $(22)$-orbits. The first consists of pairs $(A,B)\in W(\R)$, where $A$ is anisotropic; the second consists of pairs $(A,B)$, where $A$ is isotropic and $B$ takes positive values on the zeroes of $A$; and the third consists of pairs $(A,B)$, where $A$ is isotropic and $B$ takes negative value on the zeroes of $A$. We denote these refinements of the splitting type $(22)$ by $(22\#)$, $(22+)$, and $(22-)$, respectively. For any set $T\subset W(\R)$ and splitting type $(i)$, denote the set of elements in $T$ with splitting type $(i)$ by $T^{(i)}$.

\medskip
\noindent \textbf{Local masses.} Given $A,B\in V(\Z)$ of nonzero determinant, write $S_{AB} \defeq S_A\times S_B$. Next, we introduce the notion of local masses of pairs $(A,B)$ over binary cubic forms.

\begin{defn}
For each $f\in U_{a,d}(\Z)$, each $(A,B)\in W_{a,d}(\Z)$, each splitting type $(i)$ at infinity, and each $\SL_3(\Z)$-invariant function $\phi \colon W_{a,d}(\Z) \to [0,1]$ defined by congruence conditions, define \mbox{the following masses:}
\begin{equation*}
\begin{array}{rcl}
\on{Mass}_\infty\big(f;A,B;(i)\big)&:=&\displaystyle
\sum_{\substack{w\in\frac{S_{AB}(\R)^{(i)} \cap \on{Res}^{-1}(f)}{\SL_3(\R)}}}\frac{1}{\#\Stab_{\SL_3(\R)}(w)};
\\[.2in]
\on{Mass}_p\big(f;A,B;\phi\big)&:=&\displaystyle
\sum_{\substack{w\in\frac{S_{AB}(\Z_p)\cap \on{Res}^{-1}(f)}{\SL_3(\Z_p)}}}\frac{\phi_p(w)}{\#\Stab_{\SL_3(\Z_p)}(w)}.
\end{array}
\end{equation*}
\end{defn}
\noindent Given $k \in \Z \smallsetminus \{0\}$, let $\mc{G}_k$ denote the set of genera of quadratic forms in $V_k(\Z)$. 
Note that the masses defined above depend only on the genera of $A$ and $B$.

\medskip
We are now in position to state the main result of this section, which is as follows:
\begin{theorem}\label{th:sec4main}
Let $\circ \in \{\on{bal},\on{wei}\}$, let $(i)$ be a splitting type, and let $\phi \colon W_{a,d}(\Z) \to [0,1]$ be an $\SL_3(\Z)$-invariant function defined by finitely many congruence conditions. Then we have
\begin{equation*}
\sum_{\substack{w\in \frac{W(\Z)^{(i),\circ}_X}{\SL_3(\Z)}\\w\text{ \rm non-dist.}}}\phi(w)\sim
4\sum_{\substack{A\in \mc{G}_a\\B\in \mc{G}_{-d}}}\int_{f\in U_{a,d}(\R)^{\pm,\circ}_X}\on{Mass}_\infty\big(f;A,B;(i)\big)\mathrm{d}f\cdot\prod_p\int_{f\in U_{a,d}(\Z_p)}\on{Mass}_p(f;A,B;\phi)\mathrm{d}f .\footnote{Given two functions $h_1(X)$ and $h_2(X)$, we write $h_1 \sim h_2$ if the limit $\lim_{X \to \infty} h_1(X)/h_2(X) = 1$.} %,
\end{equation*}
\end{theorem}

 We now briefly outline the proof of Theorem~\ref{th:sec4main} in the setting of balanced height; the proof in the weighted height case is similar. Proving Theorem~\ref{th:sec4main} amounts to counting lattice points $(A,B) \in W_{a,d}(\Z)$ in fundamental sets for the action of $\on{SL}_3(\Z)$ on $W(\R)$. In \S\ref{sec-reductave}, we define the relevant fundamental sets and apply Bhargava's averaging method to express the count as an integral over $\gamma \in \on{SL}_3(\Z) \backslash \on{SL}_3(\R)$ of the number of non-distinguished points in a certain region $\mc{R}_\gamma \subset \gamma W_{a,d}(\R)^{(i)}$. In \S\ref{sec-cuspcut}, we prove that the contribution to this integral from the cuspidal regions is negligible, and in \S\ref{sec-dister}, we prove that, on average over $\gamma$, the region $\mc{R}_\gamma$ contains negligibly many distinguished and $\Delta$-distinguished points. Using the results of \S\ref{sec-countonsymmetricvarieties}, we express the number of integer elements in $\mc{R}_\gamma$ in terms of the volume of $\mc{R}_\gamma$. We compute this volume via a Jacobian change-of-variables argument in \S\ref{sec-jacobianchvals}. Finally, we complete the proof of Theorem~\ref{th:sec4main} in \S\ref{sec-congr}.

\subsection{Reduction theory and averaging} \label{sec-reductave}

There are two approaches to counting $\on{SL}_3(\Z)$-orbits on $W_{a,d}(\Z)$. We could do so directly, or we could realise the set of $\on{SL}_3(\Z)$-orbits on $W_{a,d}(\Z)$ as the union over $A \in V_a(\Z)$ of the set of $\on{SO}_A(\Z)$-orbits on $\{A\} \times V_{-d}(\Z)$. We shall use both approaches; the first is amenable to working with the balanced height, whereas the second is amenable to working with the weighted height.

Fix $A\in V_a(\Z)$. 
The following proposition follows immediately from \cite[Constructions 3.20 and 4.8]{BHSpreprint} and gives fundamental sets for the action of $\on{SL}_3(\R)$ (resp., $\on{SO}_A(\R)$) on subsets of $W(\R)$ (resp., $\{A\} \times V(\R)$) having given splitting type:

\begin{proposition}
 There exists continuous sections
\begin{equation*}
\begin{array}{rcl}
s^{(i)}_\bx:U(\R)^+&\to&W(\R)^{(i)} \;\;{\rm for}\; i\in\{0,2\#,2+,2-\},
\\[.1in]
s^{(1)}_\bx:U(\R)^-&\to&W(\R)^{(1)}, 
\\[.1in]
s_A^{(i)}:U_a(\R)^+&\to& \big(\{A\} \times V(\R)\big)^{(i)}
\;\;\;{\rm for}\;
i\in\{0,2+,2-\} \mbox{ when $A$ is isotropic},
\\[.1in]
s^{(2\#)}_A:U_a(\R)^+&\to& \big(\{A\} \times V(\R)\big)^{(2\#)}
\mbox{    \,when $A$ is anisotropic},
\\[.1in]
s^{(1)}_A:U_a(\R)^-&\to& \big(\{A\} \times V(\R)\big)^{(1)}
\;\;\;\mbox{when $A$ is isotropic},
\end{array}
\end{equation*}
such that the following properties hold for all $f$: $\Res(s_\bx(f))=f$ and $\Res(s_A(f))=f$; the coefficients of $s_\bx(f)$ are $\ll \rmH_\bx(f)^{1/3}$; and the coefficients of $s_A(f)$ are $\ll \rmH_\wei(f)$.
\end{proposition}

\noindent We denote the image of $s_\bx^{(i)}$ by $R_\bx^{(i)}$ and the image of $s_A^{(i)}$ by $R_A^{(i)}$. 

\medskip

Next, we construct fundamental domains for the action of $\on{SL}_3(\Z)$ on $\on{SL}_3(\R)$ (resp., $\on{SO}_A(\Z)$ on $\on{SO}_A(\R)$). Consider the Iwasawa decomposition $\SL_3(\R)=N_3T_3K_3$, where $N_3$ is the subgroup of lower-triangular unipotent matrices, $T_3$ is the subgroup of diagonal matrices with nonnegative entries, and $K_3 = \SO_3(\R)$. We write elements in $T_3$ as $(s_1,s_2)$, where $(s_1,s_2)$ corresponds to the diagonal matrix with diagonal entries $s_1^{-2}s_2^{-1}$, $s_1s_2^{-1}$, and $s_1s_2^2$. 
By a well-known result of Minkowski (see~\cite{MR1580668}), we may choose a fundamental domain $\FF_3$ for the action of $\SL_3(\Z)$ on $\SL_3(\R)$ that lies in a Siegel domain $N_3'T_3'K_3$, where $N_3' \subset N_3$ is a bounded subset and $T_3' = \{(s_1, s_2) \in T_3' : s_1, s_2 > c\}$ for some constant $c > 0$. Then for $i\in \{0,1,2\#,2+,2-\}$, the multiset $\FF_3\cdot R_\bx^{(i)}$ is generically a $\sigma_i$-fold cover of a fundamental set for the action of $\SL_3(\Z)$ on $W(\R)^{(i)}$, where $\sigma_0=\sigma_{2+}=\sigma_{2-}=\sigma_{2\#}=\sigma^+ = 4$, and $\sigma_1 = \sigma^- = 2$.

Following \cite[Construction 4.13]{BHSpreprint}, we choose a fundamental domain $\FF_A$ for the action of $\SO_A(\Z)$ on $\SO_A(\R)$. When $A$ is anisotropic over $\Q$, this domain $\FF_A$ may be chosen to be compact. On the other hand, when $A$ is isotropic, then $A$ is $\SL_3(\Q)$-equivalent to the element $A_a \in S(\Z)$, where $A_a$ is antidiagonal with entries $1/2$, $-a$, and $1/2$. That is, there exists $g_A\in\SL_3(\Q)$ with $g_A A g_A^T=A_a$. We thus obtain the following isomorphisms:
\begin{equation*}
\begin{array}{rclrcl}
m_A\colon \{A\} \times V(\R)&\to& \{A_a\} \times V(\R),\quad &(A,B)&\mapsto & (A_a,g_A Bg_A^T);
\\[.1in]
m_A'\colon \SO_A(\R)&\to&\SO_{A_a}(\R),\quad &g&\mapsto& g_A \cdot g\cdot g_A^{-1}.
\end{array}
\end{equation*}
Consider the Iwasawa decomposition $\on{SO}_{A_a}(\R) = N_aT_aK_a$, where $N_a$ is the subgroup of lower-triangular unipotent matrices, $T_a$ is the subgroup of diagonal matrices with nonnegative entries, and $K_a$ is a maximal compact subgroup. We write elements of $T_a$ as $t$, where $t$ corresponds to the diagonal matrix diagonal entries $t^{-2}$, $1$, and $t^2$. Choose a suitable Siegel domain $\mathscr{S}_a=N_a'T_a'K_a$ for the action of $\SO_{A_a}(\Z)$ on $\SO_{A_a}(\R)$, where $N_a' \subset N_a$ is a bounded subset and $T_a' = \{t \in T_a : t > c\}$ for some constant $c > 0$. Let $\mathscr{S}_A \defeq {m_A'}^{-1}(\mathscr{S}_a)$.
By \cite[Example 2.5]{MR0148666}, there exists a fundamental domain $\FF_A$ for the action of $\SO_A(\Z)$ on $\SO_A(\R)$ contained within a finite union of translates $\bigcup_i g_i\mathscr{S}_A$, where each $g_i$ belongs to $\SO_A(\Q)$.
That is, we may write
\begin{equation}\label{eq:FcA}
\FF_A\subset \bigcup_{i}g_ig_A^{-1}\mathscr{S}_ag_A=g_A^{-1}h_i\mathscr{S}_ag_A,
\end{equation}
where $h_i \defeq m_A'(g_i)$.
Then for $i\in\{0,1,2\#,2+,2-\}$, the multiset $\FF_A\cdot R_A^{(i)}$ is generically a $\sigma_i$-fold cover of a fundamental set for the action of $\SO_A(\Z)$ on $(\{A\} \times V(\R))^{(i)}$.

\medskip

Now, for any $\SL_3(\Z)$-invariant subset $L \subset W_{a,d}(\Z)$, define
\begin{equation*}
N_{a,d}^{\bx}(L,X) \defeq \#\bigl(\SL_3(\Z)\backslash L^{\nd,\bx}_X\bigr),\quad\quad
N_{a,d}^{\wei}(L,X) \defeq \#\bigl(\SL_3(\Z)\backslash L^{\nd,\wei}_X\bigr),
\end{equation*}
where $L^\nd \subset L$ denotes the subset of non-distinguished elements, with resolvent having Galois group $S_3$.
Let $G_A\subset \SO_A(\R)$ be a nonempty open bounded subset with measure-$0$ boundary, and let $G_0\subset \SL_3(\Z)$ be a subset satisfying the same conditions. 
Upon combining the above discussion with Bhargava's averaging method (see, e.g.,~\cite{MR2183288,MR2745272}), we obtain the following result:
\begin{proposition}\label{prop:avg}
Fix $i\in\{0,1,2\#,2+,2-\}$. Let $L\subset W_{a,d}(\Z)^{(i)}$ be an $\SL_3(\Z)$-invariant subset, and let $L_A \defeq L\cap \{A\} \times V(\Z)$. Then we have
\begin{equation*}
\begin{array}{rcl}
N_{a,d}^{\bx}(L,X)&=&\displaystyle
\frac{1}{\sigma_{(i)}\Vol(G_0)}\int_{\gamma\in\FF_3}\#\Bigl(L^\nd \cap \gamma G_0\bigl(R_\bx^{(i)}\bigr)^\bx_X\Bigl)\mathrm{d}\gamma;
\\[.2in]
N_{a,d}^{\wei}(L,X)&=&\displaystyle
\sum_{A\in\frac{V_a(\Z)}{\SL_3(\Z)}}
\frac{1}{\sigma_{(i)}\Vol(G_A)}\int_{\gamma\in\FF_A}\#\Bigl(L_A^\nd \cap \gamma G_A\bigl(R_A^{(i)}\bigr)^\wei_{X}\Bigl)\mathrm{d}\gamma.
\end{array}
\end{equation*}
\end{proposition}

\subsection{Cutting off the cusp} \label{sec-cuspcut}

Now take $A$ to be isotropic. For $\delta>0$, let $\FF_{3}^{(\delta)} \defeq \{n(s_1,s_2)k \in \mc{F}_3 : \max\{s_1, s_2\} > X^\delta\}$, and let $\FF_{A}^{(\delta)} \defeq \{ntk \in \mc{F}_A : t > X^\delta\}$. 
We think of $\FF_3^{(\delta)}$ and $\FF_A^{(\delta)}$ as the \emph{cuspidal regions} of the fundamental domains 
$\FF_{3}$ and $\FF_{A}$, respectively.
In this subsection, we prove that the contribution to Proposition \ref{prop:avg} from the cuspidal regions is negligible. 
\begin{theorem}\label{th:cuspcut}
We have
\begin{equation*}
\begin{array}{rcl}
\displaystyle
\int_{\gamma\in\FF_3^{(\delta)}}\#\Bigl(W_{a,d}(\Z)^\nd \cap \gamma G_0\bigl(R_\bx^{(i)}\bigr)^\bx_{X}\Bigl)\mathrm{d}\gamma &\ll&
X^{2-\delta_1},
\\[.2in]
\displaystyle
\int_{\gamma\in\FF_A^{(\delta)}}\#\Bigl(W_{a,d}(\Z)^\nd \cap \gamma G_A\bigl(R_A^{(i)}\bigr)^\wei_{X}\Bigl)\mathrm{d}\gamma &\ll&
X^{3-\delta_2},
\end{array}
\end{equation*}
for some $\delta_1,\delta_2 > 0$.
\end{theorem}

As it happens, the asymptotic derived in \S\ref{sec-countonsymmetricvarieties} using dynamical techniques is insufficient to prove Theorem~\ref{th:cuspcut}, as bounding the contribution from the cuspidal regions requires us to have a more explicit understanding of the error dependence on the skewing factor $\gamma$. Instead of making this dependence more explicit, we take a different approach: using a completely elementary argument, we obtain a tight upper bound on the number of ternary quadratic forms in skewed boxes having fixed determinant $k/4$ for $k \in \Z \smallsetminus \{0\}$. To this end, write the matrix $A\in V_k(\Z)$ as
$$A = \begin{pmatrix} a & b/2 & d/2 \\ b/2 & c & e/2 \\ d/2 & e/2 & f \end{pmatrix}$$
Let $\cB\subset V(\R)$ be a bounded subset. We write the equation $\det A=k/4$ as
\begin{equation} \label{eq-det-k}
-k=-4\det(A)= f(b^2-4ac)+(ae^2-bed+cd^2).
\end{equation}
Given a triple $(a,b,c) \in \Z^3$, let $Q(x,y) \defeq Q_{a,b,c}(x,y)$ denote the quadratic form $ax^2-bxy+cy^2$, which has discriminant $\Delta \defeq b^2 - 4ac$. Then our equation~\eqref{eq-det-k} takes the cleaner form
\begin{equation*}
-k=f\Delta+Q(e,d).
\end{equation*}
Now, for $s_1,s_2 \gg 1$, put 
\begin{equation*}
\begin{array}{lcl}
\displaystyle N(s_1,s_2;Y)& \defeq & \#\bigl( V_k(\Z) \cap (s_1,s_2)Y\cB\bigr);
\\[.1in]\displaystyle 
N^*(s_1,s_2;Y)& \defeq &
\#\bigl\{A\in V_k(\Z) \cap (s_1,s_2)Y\cB : a\neq 0\bigr\}.
\end{array}
\end{equation*}
and subdivide these counts as follows:
\begin{equation*}
\begin{array}{rcl}
N(s_1,s_2;Y)&=&N_{\Delta\neq 0}(s_1,s_2;Y)+N_{\Delta= 0}(s_1,s_2;Y),\\[.1in]
N^*(s_1,s_2;Y)&=&N^*_{\Delta\neq 0}(s_1,s_2;Y)+N^*_{\Delta= 0}(s_1,s_2;Y),
\end{array}
\end{equation*}
where $N_{\Delta\neq 0}(s_1,s_2;Y)$ is the contribution to $N(s_1,s_2;Y)$ from ternary quadratic forms $A$ with $\Delta\neq 0$, and $N_{\Delta= 0}(s_1,s_2;Y)$ is the contribution from forms $A$ with $\Delta= 0$, with the analogous convention for the counts $N^*$. 
Then we have the following result:

\begin{theorem}\label{th:skbound}
For any $Y>1$, we have
\begin{equation*}
\begin{array}{lclll}
\displaystyle N_{\Delta\neq 0}(s_1,s_2;Y)&\ll_\varepsilon& s_1^3Y^{3+\varepsilon};\quad\quad
N_{\Delta= 0}(s_1,s_2;Y)&\ll_\varepsilon& s_2^3Y^{3+\varepsilon}+s_1^4s_2^5Y^{2+\varepsilon};
\\[.1in]\displaystyle 
N^*_{\Delta\neq 0}(s_1,s_2;Y)&\ll_\varepsilon& Y^{3+\varepsilon};\quad\quad\,\,\,\,\,\,
N^*_{\Delta= 0}(s_1,s_2;Y)&\ll_\varepsilon& s_2^3Y^{3+\varepsilon},
\end{array}
\end{equation*}
where the implied constants depend on the region $\cB$.
\end{theorem}
\begin{proof}[Proof of Theorem~\ref{th:skbound}]
We denote the range of each matrix entry $\alpha$ in the set $(s_1,s_2)Y\cB$ by $R_\alpha$. Note that we have
\begin{equation*}
\begin{array}{rcl}
&R_a\ll s_1^{-4}s_2^{-2}Y;\quad
R_b\ll s_1^{-1}s_2^{-2}Y;\quad
R_c\ll s_1^{2}s_2^{-2}Y;&
\\[.05in]
&R_d\ll s_1^{-1}s_2Y;\quad
R_e\ll s_1^{2}s_2Y;&
\\[.05in]
&R_f\ll s_1^{2}s_2^{4}Y.&
\end{array}
\end{equation*}
We now split into four cases depending on whether or not $a$ and $\Delta$ are zero:

\medskip 
\noindent {\bf Case I: $a\neq 0$ and $\Delta\neq 0$.}
We start by bounding $N^*_{\Delta\neq 0}(s_1,s_2;Y)$.
For fixed $a,b,c \in \Z$ with $a,\Delta\neq 0$, let $N_{a,b,c}(k)$ denote the number of matrices $A \in V_k(\Z) \cap (s_1,s_2)X\cB$ where the $11$-, $12$-, and $22$-matrix entries of $A$ are $a$, $b/2$, and $c$, respectively.
Then we have
\begin{equation*}
\begin{array}{rcl}
N_{a,b,c}(k)&\ll&\#\bigl\{(d,e)\in\Z^2:|d|\ll R_d,\,|e|\ll R_e,\,Q(e,d)\equiv -k\pmod\Delta\bigr\}
\\[.1in]&\ll&
\#\bigl\{(d,e)\in\Z^2:|d|\ll R_d,\,|e|\ll R_e,\,4aQ(e,d)\equiv -4ak\pmod\Delta\bigr\}
\\[.1in]&\ll&
\#\bigl\{(d,e)\in\Z^2:|d|\ll R_d,\,|e|\ll R_e,\,(2ae-bd)^2\equiv -4ak\pmod\Delta\bigr\},
\end{array}
\end{equation*}
since $4aQ(x,y)\equiv (2ax-by)^2\pmod\Delta$.
This forces $(2ae-bd)^2$ to lie within the set $T$ of square roots of $-4ak$ modulo $\Delta$. Write $\Delta=s_\Delta^2q_\Delta$, where $q_\Delta$ is squarefree. Then $\#T\ll Y^\varepsilon s_\Delta$.
For any $\ell\in T$, two solutions $(e_1,d_1)$ and $(e_2,d_2)$ to $2ax-by\equiv\ell \pmod \Delta$ gives a solution $(e'=e_1-e_2,d'=d_1-d_2)$ to $2ax-by\equiv 0\pmod\Delta$; note that we have $|d'|\ll R_d$ and $|e'|\ll R_e$. It follows that
\begin{equation*}
N_{a,b,c}(k)\ll Y^\varepsilon s_\Delta \#\bigl\{(d',e')\in\Z^2:|d'|\ll R_d,\,|e'|\ll R_e,\, 2ae'-bd'\equiv 0\pmod{\Delta}\bigr\}.
\end{equation*}
We now fibre over $g$, $b \neq 0$, and $\Delta$, where $g = \gcd(b,\Delta)$. With these fixed, we also fibre over $\sigma$, the possible values of $2ae'-bd'$, which must be a multiple of $\Delta$. This gives a contribution of $\ll$
\begin{equation*}
\begin{array}{ll}
&\displaystyle\sum_{\substack{g\geq 1}}\sum_{\substack{|b|\ll R_b,|\Delta|\ll R_b^2\\ \gcd(b,\Delta)=g}}\sum_{\substack{a,c\\b^2-ac=\Delta}}
\sum_{\substack{|\sigma|\ll R_aR_e\\\Delta\mid\sigma}}
Y^\varepsilon s_{\Delta}\#\{(d',e')\in\Z^2:|d'|\ll R_d,\,|e'|\ll R_e,\, 2ae'-bd'=\sigma\} \ll
\\[.25in]
 &\displaystyle
\sum_{\substack{\\g\geq 1}}\sum_{\substack{|b|\ll R_b,|\Delta|\ll R_b^2\\ \gcd(b,\Delta)=g}}\sum_{\substack{a,c\\b^2-ac=\Delta}}
\sum_{\substack{|\sigma|\ll R_aR_e\\\Delta\mid\sigma}}
Y^\varepsilon s_{\Delta}\#\{(d'',e'')\in\Z^2:|d''|\ll R_d,\,|e''|\ll R_e,\, 2ae''-bd''=0\},
\end{array}
\end{equation*}
where the second estimate follows from the same trick of replacing two solutions $(d_1',e_1')$ and $(d_2',e_2')$ with $(d''=d_1'-d_2',e''=e_1'-e_2')$. The number of such solutions $(d'',e'')$ once $a,b \in \Z \smallsetminus \{0\}$ have been fixed is $\ll R_e(a,b)/b\leq gR_e/b$. We thus obtain a contribution of $\ll$
\begin{equation*}
Y^\varepsilon R_e\sum_{\substack{g\geq 1}}\sum_{\substack{|b|\ll R_b,|\Delta|\ll R_b^2\\ \gcd(b,\Delta)=g}}\sum_{\substack{\sigma\neq 0\\|\sigma|\ll R_aR_e\\\Delta\mid\sigma}}
\frac{s_{\Delta} g}{b}\ll Y^\varepsilon R_aR_e^2=Y^{3+\varepsilon},
\end{equation*}
as necessary. Meanwhile, the contribution from $b=0$ is $\ll R_aR_dR_e \ll s_1^{-3}Y^3$, which is also sufficiently small. This concludes the first case.

\medskip 
\noindent {\bf Case II: $a\neq 0$ and $\Delta= 0$.} In this case, since $\Delta=0$, we may write our quadratic form $Q$ as $Q(x,y)=\alpha(\beta x+\gamma y)^2$, with $\alpha,\beta,\gamma \in \Z$ such that $\alpha\beta\neq 0$. Since $\det A =k/4$ imposes no condition on $f$, we cannot improve on the $\ll R_f$ choices for $f$. It remains to estimate the number of possible $(\alpha,\beta,\gamma,d,e)$ satisfying
\begin{equation*}
\alpha(\beta e+\gamma d)^2=-k.
\end{equation*}
This restricts the value of $\alpha$, which must divide $k$, to $\ll |k|^\varepsilon$ possibilities. For each possible $\alpha$, we have $\beta e+\gamma d=\pm \sqrt{-k/\alpha}$, and we have the bounds $|\beta|\ll R_a^{1/2}$ and $|\gamma|\ll R_c^{1/2}$. Fibreing over $g$, $\beta \neq 0$, and $\gamma$, where $g = \gcd(\beta,\gamma)$ and using the trick from Case I, we find that
\begin{equation*}
\begin{array}{rcl}
N^*_{\Delta=0}(s_1,s_2;Y)&\ll&
\displaystyle R_f\sum_{g\geq 1}\sum_{\substack{|\beta|\ll R_a^{1/2},|\gamma|\ll R_c^{1/2} \\ \gcd(\beta,\gamma) = g}}
\#\{|d'|\ll R_d,\,|e'|\ll R_e,\, \beta e'+\gamma d'=0\}
\\[.25in]
&\ll&
\displaystyle  R_c^{1/2}R_f

\sum_{g\geq 1}\frac{1}{g}\sum_{|\beta|\ll R_a^{1/2}} gR_d/\beta\ll Y^\varepsilon R_c^{1/2}R_dR_f=s_1^2s_2^4Y^{5/2+\varepsilon}.
\end{array}
\end{equation*}
Since $R_a\gg 1$ (for otherwise, $a\neq 0$ would be impossible), we may multiply our bound by $R_a^{1/2}$ to obtain a bound of $s_2^3Y^{3+\varepsilon}$, as claimed.

\medskip 
\noindent {\bf Case III: $a= 0$ and $\Delta\neq 0$.} In this case, we have $\Delta=b^2\neq 0$, and so $s_\Delta=|b|$. Our equation is
\begin{equation*}
-k= b^2f+ (cd^2-bde).
\end{equation*}
If $c=0$, we fibre over $b$, $d$, and $e$; this fixes $f$, and the total number of possibilities is $O(Y^3)$, which is sufficiently small. On the other hand, for fixed $b,c \neq 0$, we obtain a contribution of
\begin{equation*}
\begin{array}{rcl}
N_{0,b,c}(k)&\ll&\#\bigl\{(d,e)\in\Z^2:|d|\ll R_d,\,|e|\ll R_e,\,(2cd-be)^2\equiv -4kc\pmod{b^2}\bigr\}
\\[.1in]&\ll&
\displaystyle |b|\frac{R_bR_e}{b^2}\#\bigl\{(d,e)\in\Z^2:|d|\ll R_d,\,|e|\ll R_e,\, 
2cd'-be'= 0\bigr\},
\end{array}
\end{equation*}
where the second estimate follows because $-4kc$ has $\ll |b|$ square roots modulo $b^2$, which gives a total of $\ll |b| R_bR_e/b^2$ possible values for $2cd-be$; then, for each such value, we apply our previous trick of replacing solutions $(d_1,e_1)$ and $(d_2,e_2)$ with $(d_1-d_2,e_1,e_2)$ to obtain a solution to $cd'-be'=0$. 
Fibreing over $g$, $b$, and $c$, where $g = \gcd(b,c)$, we obtain a contribution of $\ll$
\begin{equation*}
R_bR_e\sum_{g\geq 1}\sum_{\substack{|b|\ll R_b,|c|\ll R_c\\ \gcd(b,c) = g}}\frac{1}{|b|}\frac{gR_e}{|c|}\ll Y^{\varepsilon}R_bR_e^2\ll s_1^3Y^{3+\varepsilon},
\end{equation*}
as necessary.

\medskip 
\noindent {\bf Case IV: $a=0$ and $\Delta= 0$.}
This is equivalent to stipulating that $a=b=0$, and our equation becomes $-k=cd^2$, which determines $c$ and $d$ up to $|k|^\varepsilon$ choices. The coefficients $e$ and $f$ can take any values, leading to a total of $\ll R_eR_f=s_1^4s_2^5Y^{2}$ choices.
\end{proof}

\noindent With Theorem~\ref{th:skbound} in hand, we are now in position to prove the main result of the subsection.

\medskip

\noindent {\bf Proof of Theorem \ref{th:cuspcut}:} We begin with the first equation and write
\begin{equation*}
\begin{array}{rcl}
\displaystyle
\int_{\gamma\in\FF_3^{(\delta)}}\#\Bigl(W_{a,d}(\Z)^\nd \cap \gamma G_0\bigl(R_\bx^{(i)}\bigr)^\bx_{X} \Bigl)\mathrm{d}\gamma 
&\ll&
\displaystyle\int_{\substack{s = (s_1, s_2)\\ \max\{s_1,s_2\}\gg X^{\delta}\\s_1^4s_2^2\ll X^{1/3}}}
\#\Big(W_{a,d}(\Z)^{\nd} \cap s X^{1/3} S\Big)\frac{\mathrm{d}^\times s}{s_1^6s_2^6},
\end{array}
\end{equation*}
where $S=\cB\times\cB$ for some bounded set $\cB\subset V(\R)$, and where we can impose the condition $s_1^4s_2^2 \ll X^{1/3}$ because if $s_1^4s_2^2>CX^{1/3}$ for some sufficiently large $C$, then $|a_{11}|,|b_{11}|<1$ for every $(A,B)\in s X^{1/3}S$, implying that $W_{a,d}(\Z)^{\nd} \cap s X^{1/3}S = \varnothing$. The required bound now follows from an application of Theorem~\ref{th:skbound}. The second equation in Theorem \ref{th:cuspcut} \mbox{follows in identical fashion. \hspace*{\fill} \qed}

\subsection{Bounding the number of distinguished and $\Delta$-distinguished points in the main body} \label{sec-dister}

The complement of the cuspidal region of a fundamental domain is known colloquially as the \emph{main body}. We have the following lemma proving that the counts of distinguished and $\Delta$-distinguished points in the main bodies of $\mc{F}_3$ and $\mc{F}_A$ are negligible.
\begin{lemma}\label{lem:mbdist}
We have
\begin{equation*}
\begin{array}{rcl}
\displaystyle
\int_{\gamma\in\FF_3\smallsetminus \FF_3^{(\delta)}}\#\Bigl((W_{a,d}(\Z) \smallsetminus W_{a,d}(\Z)^\nd) \cap \gamma G_0\bigl(R_\bx^{(i)}\bigr)^\bx_{X}\Bigl)d\gamma &=&
o(X^{2});
\\[.2in]
\displaystyle
\int_{\gamma\in\FF_\cA\smallsetminus \FF_\cA^{(\delta)}}\#\Bigl((W_{a,d}(\Z) \smallsetminus W_{a,d}(\Z)^\nd) \cap \gamma G_\cA\bigl(R_\cA^{(i)}\bigr)^\wei_{X}\Bigl)d\gamma &=&
o(X^{3}).
\end{array}
\end{equation*}
\end{lemma}
\begin{proof}
Let $p\nmid 2ad$ be a prime. Then as shown in Proposition~\ref{prop-maxcalcs}, a positive proportion (bounded away from $0$, independent of $p$) of forms in $U_{a,d}(\F_p)$ have three distinct roots in $\P^1(\F_p)$. For each such form $f$, it follows from, e.g.,~\cite[Theorem~5.7]{MR3187931} that there exist $4$ distinct $\SL_3(\F_p)$-orbits on $W_{a,d}(\F_p)$ with resolvent $f$. At least two of these orbits are neither distinguished nor $\Delta$-distinguished. Hence, a positive proportion (bounded away from $0$, independent of $p$) of elements in $W_{a,d}(\F_p)$ are neither distinguished nor $\Delta$-distinguished.

Next, any $f \in U_{a,d}(\Z)$ with mod-$p$ reduction having an irreducible quadratic factor has Galois group $S_3$. As shown in Proposition~\ref{prop-maxcalcs}, a positive proportion (bounded away from $0$, independent of $p$) of forms in $U_{a,d}(\F_p)$ have an irreducible quadratic factor. For each such form $f$, it follows from, e.g.,~\cite[Theorem~5.7]{MR3187931} that there exist $2$ distinct $\SL_3(\F_p)$-orbits on $W_{a,d}(\F_p)$ with resolvent $f$. Hence, a positive proportion (bounded away from $0$, independent of $p$) of elements in $W_{a,d}(\F_p)$ cannot have resolvent equal to the mod-$p$ reduction of a form $f \in U_{a,d}(\Z)$ with Galois group $C_3$. %\todo{\small added this paragraph; is it correct?}

The lemma now follows from the above two paragraphs via an argument identical to that given in \cite[Lemma 13]{MR2745272}, with the only difference being that we apply the counting results from \S\ref{sec-countonsymmetricvarieties} instead of geometry-of-numbers methods.
\end{proof}

\subsection{A Jacobian change-of-variables} \label{sec-jacobianchvals}

In this subsection, we derive a Jacobian change-of-variables formula relating measures on $U_{a,d}$ and $W_{a,d}$. This formula is a key ingredient in the proof of the main result of this section, namely Theorem~\ref{th:sec4main}; see \S\ref{sec-congr}.

\subsubsection*{Choice of measures}
Regard $\mathbb{A}^1$, $U$, $V$, and $W$ as schemes over $\Q$. We now choose a volume form on each of them. First, on $\mathbb{A}^1 = \on{Spec}\Q[t]$ we take our volume form to be $\mathrm{d}t$. Second, on $U = \on{Spec} \Q[a,b,c,d]$ we take our volume form to be $\mathrm{d}f \defeq \mathrm{d}a \wedge \mathrm{d}b \wedge \mathrm{d}c \wedge \mathrm{d}d$. For fixed $a,d \in \Z \smallsetminus \{0\}$, we abuse notation by writing $\mathrm{d}f=db\wedge dc$ to denote the restriction of the form $\mathrm{d}f$ on $U$ to the subscheme $U_{a,d}$. Third, on $V = \on{Spec}\Q[\{B_{ij} : 1 \leq i \leq j \leq 3\}]$ we take our volume form to be $\mathrm{d}v \defeq \bigwedge_{1 \leq i \le j \leq 3} \mathrm{d}B_{ij}$. Fourth, on $W = V \times_{\Q} V$ we take our volume form to be $\mathrm{d}w \defeq \mathrm{d}v_1 \wedge \mathrm{d}v_2$, where $v_1$ (resp., $v_2$) denotes the coordinate on the left (resp., right) factor. One readily verifies that $\mathrm{d}v$ and $\mathrm{d}w$ are $\on{SL}_3$-invariant. 

Next, let $A \in V(\Z)$ be of nonzero determinant. Recall from \S\ref{sec-countonsymmetricvarieties} that we chose volume forms $\mathrm{d}g$ on $G = \on{SL}_3$ and $\mathrm{d}h$ on $H = \on{SO}_A$ defined over $\Q$, and that the measures $\mathrm{d}g$ and $\mathrm{d}h$ give rise to an $\on{SL}_3$-invariant quotient measure $\mathrm{d}s$ on $\mc{S}_A$. For the sake of clarity in this subsection, we denote these measures $\mathrm{d}g$, $\mathrm{d}h$, and $\mathrm{d}s$ by $\omega_{\on{SL}_3}$, $\omega_{\on{SO}_A}$, and $\omega_{\mc{S}_A}$, respectively. We now give a useful explicit description of the measure $\omega_{\mc{S}_A}$. Let $\on{det} \colon V \rightarrow \mathbb{A}^1$ denote the determinant map. Choose a differential $5$-form $\omega$ on $V$, defined over $\Q$, such that 
\begin{equation} \label{eq-vasomega}
\mathrm{d}v = \omega \wedge \on{det}^*\mathrm{d}t.
\end{equation}
Note that such a form $\omega$ exists; as an example, we could take $\omega$ explicitly as follows. Letting $\on{adj} \colon V \to V$ denote the adjugate map, we have 
$$\frac{\partial \det B}{\partial B_{ij}} = \begin{cases} 2 (\on{adj}B)_{ji} , & \text{if $i\neq j$,} \\  (\on{adj}B)_{ii} , & \text{if $i=j$,} \end{cases}$$
so $\on{det}^* \mathrm{d}t = \sum_i \on{adj}(B)_{ii} \mathrm{d}B_{ii} + \sum_{i < j} 2 \on{adj}(B_{ij}) \mathrm{d}B_{ij}$, and we can take \[\omega = \frac{1}{3\det B} \left(\sum_{i \le j} \pm B_{ij} \mathrm{d}B_{11}\wedge \ldots \wedge \widehat{\mathrm{d}B_{ij}} \wedge \ldots \wedge \mathrm{d}B_{nn}\right)\] 
where the $\pm$ signs are chosen such that $\mathrm{d}v = \mathrm{d}B_{ij} \wedge (\pm \mathrm{d}B_{11}\wedge \ldots \wedge \widehat{\mathrm{d}B_{ij}} \wedge \ldots \wedge \mathrm{d}B_{nn})$. With such a choice of $\omega$, we can now take $\omega_{\mc{S}_A}$ to be the pullback of $\omega$ via the inclusion $\mc{S}_A \subset V$. It follows from the proof of~\cite[Lemma~A6]{MR2155083} that $\omega_{\mc{S}_A}$ is independent of the choice of $\omega$. It is easy to verify that $\omega_{\mc{S}_A}$ is left-invariant under the action of $\on{SL}_3$: indeed, given any field $K/\Q$ and any $g \in \on{SL}_3(K)$, we have the following equalities of differential forms over $K$:
\[\mathrm{d}v = g^* \mathrm{d}v = g^*\omega \wedge (\on{det} \circ\, g)^*(\mathrm{d}t) =  g^*\omega \wedge \on{det}^* \mathrm{d}t,\]
from which it follows that $g^*\omega_{\mc{S}_A} = \omega_{\mc{S}_A}$. Now, let $\iota,\,\pi$ be as in~\eqref{eq-exactgroups}. Then the form $\iota_* \omega_{\on{SO_A}} \wedge \pi^* \omega_{\mc{S}_A}$ is nowhere-vanishing top-degree differential form on $\on{SL}_3$ defined over $\Q$, and so we may rescale $\omega_{\on{SL}_3}$ to arrange that $\omega_{\on{SL}_3} = \iota_* \omega_{\on{SO_A}} \wedge \pi^* \omega_{\mc{S}_A}$.

\subsubsection*{Statement of the formula}

Fix $A,B \in V(\Z)$ of determinants $a,d$, respectively, and write $\mc{S}_{AB} \defeq \mc{S}_A \times \mc{S}_B$. We prove the following Jacobian change-of-variables formulas relating the Haar measure on $\mc{S}_{AB}$ with the product measure on $U_{a,d} \times \on{SL}_3$, as well as the Haar measure on $\mc{S}_B$ with the product measure on $U_{a,d} \times \on{SO}_A$:

\begin{proposition} \label{prop-thejac}
Let $K$ be $\R$ or $\Z_p$ for a prime $p$, and let $|\cdot|$ denote the usual absolute value on $K$. Then there exists a nonzero rational constant $\mc{J}$ $($resp., $\mc{J}_A${}$)$, independent of $K$, such that for any measurable function $\phi$ on $\mc{S}_{AB}(K)$ $($resp., $\mc{S}_{B}(K)${}$)$, we have
\begin{align*}
    & \int_{(A',B') \in \mc{S}_{AB}(K)} \phi(A',B')\omega_{\mc{S}_A}(A')\omega_{\mc{S}_B}(B') =  |\mc{J}|\int_{\substack{f \in U_{a,d}(K) \\ \on{Disc}(f) \neq 0}} \sum_{\substack{w \in \frac{\mc{S}_{AB}(K)}{\on{SL}_3(K)} \\ \on{Res}(w) = f}} \int_{g \in \on{SL}_3(K)} \phi(g \cdot w)\omega_{\on{SL}_3}(g)\mathrm{d}f;
    \\
       & \int_{B'  \in  \mc{S}_{B}(K)} \phi(A,B')\omega_{\mc{S}_{B}}(B') =  |\mc{J}_A|\int_{\substack{f \in U_{a,d}(K) \\ \on{Disc}(f) \neq 0}} \sum_{ \substack{B' \in \frac{ \mc{S}_B(K)}{\on{SO}_A(K)} \\ \on{Res}(A,B') = f}} \int_{g \in \on{SO}_A(K)} \phi(g \cdot (A,B'))\omega_{\on{SO}_A}(g)\mathrm{d}f.
\end{align*}
\end{proposition}
\begin{proof}
We only prove the first claim of the proposition; the proof of the second claim is identical. By imitating the proof of~\cite[Proposition~3.12]{MR3272925}, we see that it suffices to take $K = \R$ and prove the following: given a piecewise-analytic section $\sigma \colon U(\R) \to W(\R)$, we have
\begin{align} \label{eq-intersect}
    \int_{\on{SL}_3(\R) \cdot \sigma(U_{a,d}(\R))} \phi(w)\mathrm{d}w & = |\mc{J}|\int_{\substack{f \in U_{a,d}(\R) \\ \on{Disc}(f) \neq 0}}  \int_{g \in \on{SL}_3(K)} \phi(g \cdot \sigma(f))\omega_{\on{SL}_3}(g)\mathrm{d}f.
\end{align}
(Note that $\on{SL}_3(\R) \cdot \sigma(U_{a,d}(\R)) \subset W_{a,d}(\R) = \mc{S}_{AB}(\R)$, so the integral on the left-hand side of~\eqref{eq-intersect} is indeed over a subset of $\mc{S}_{AB}(\R)$.) To prove~\eqref{eq-intersect}, consider the following commutative diagram: 
\begin{equation}
\begin{tikzcd} \label{eq-diag}
U \times \SL_3 \arrow[r, "\sigma"]                             & W                     \\
U_{a,d} \times \SL_3 \arrow[u, "\iota_1"] \arrow[r,swap, "\sigma"] & \mc{S}_{AB} \arrow[u, "\iota_2"']
\end{tikzcd}
\end{equation}
 By~\cite[Proposition~3.41]{BHSpreprint}, we have that $\sigma^* \mathrm{d}w = |\J| \times \mathrm{d}f \wedge \omega_{\SL_3}$ for some $\J \in \Q^\times$. As in~\eqref{eq-vasomega}, we can write $\mathrm{d}w = \mathrm{d}v_1 \wedge \mathrm{d}v_2 = (\omega_1 \wedge \det^* \mathrm{d}t_1) \wedge (\omega_2 \wedge \det^* \mathrm{d}t_2)$. Thus, we obtain the following equality of measures on $U \times \on{SL}_3$:
 \begin{equation} \label{eq-withad}
 |\J| \times \mathrm{d}f \wedge \omega_{\SL_3} = \sigma^* \mathrm{d}w = \sigma^*(\omega_1 \wedge \omega_2) \wedge \sigma^*(\on{det}^*\mathrm{d}t_1 \wedge \on{det}^*\mathrm{d}t_2) = \sigma^*(\omega_1 \wedge \omega_2) \wedge (\mathrm{d}a \wedge \mathrm{d}d).
 \end{equation}
By ``cancelling out'' $\mathrm{d}a \wedge \mathrm{d}d$ on both sides of~\eqref{eq-withad} and pulling back via $\iota_1$, we find that
  \begin{equation} \label{eq-withadcan}
 |\J| \times \mathrm{d}f \wedge \omega_{\SL_3} =  \iota_1^*\sigma^*(\omega_1 \wedge \omega_2),
 \end{equation}
 as measures on $U_{a,d} \times \on{SL}_3$.
Then, using the commutativity of the diagram~\eqref{eq-diag} along with the fact that $\iota_2^*(\omega_1 \wedge \omega_2) = \omega_{\mc{S}_{A}} \wedge \omega_{\mc{S}_{B}}$, we deduce that
\begin{equation} \label{eq-theonewithiota}
\iota_1^* \sigma^* (\omega_1 \wedge \omega_2) = \sigma^*\iota_2^*(\omega_1 \wedge \omega_2) = \sigma^*(\omega_{\mc{S}_{A}}\wedge\omega_{\mc{S}_{B}})
\end{equation}
Combining~\eqref{eq-withadcan} and~\eqref{eq-theonewithiota}, we conclude that $|\J| \times  \mathrm{d}f \wedge \omega_{\SL_3} =  \sigma^*(\omega_{\mc{S}_{A}} \wedge \omega_{\mc{S}_{B}})$ as measures on $U_{a,d} \times \on{SL}_3$, from which~\eqref{eq-intersect} is an immediate consequence.
\end{proof}

\subsection{Proof of Theorem \ref{th:sec4main}} \label{sec-congr}

Denote by $N_{a,d}^\circ(W(\Z)^{(i)};X;\phi)$ the left-hand side of Theorem \ref{th:sec4main}. We begin by considering balanced height ($\circ = \on{bal}$).  From Proposition \ref{prop:avg}, Theorem \ref{th:cuspcut}, and Lemma \ref{lem:mbdist}, we have
\begin{equation*}
\begin{array}{rcl}
N_{a,d}^{\bx}\big(W(\Z)^{(i)};X;\phi\big)&=&
\displaystyle
\frac{1}{\sigma_{(i)}\Vol(G_0)}\int_{\gamma\in\FF_3}\Bigl(\sum_{w\in W_{a,d}(\Z)^{(i),\nd} \cap \gamma G_0(R_\bx^{(i)})^\bx_X} \phi(w)\Bigr)\mathrm{d}\gamma
\\[.2in]&=&
\displaystyle
\frac{1}{\sigma_{(i)}\Vol(G_0)}\int_{\gamma\in\FF_3\smallsetminus \FF_3^{(\delta)}}\Bigl(\sum_{w\in W_{a,d}(\Z)^{(i)} \cap \gamma G_0(R_\bx^{(i)})^\bx_X} \phi(w)\Bigr)\mathrm{d}\gamma+o(X^2),
\end{array}
\end{equation*}
The function $\phi\colon W_{a,d}(\Z)=V_a(\Z) \times V_{-d}(\Z) \to\R$ is a finite sum of product functions $f_1 \times f_2$, where $f_1 \colon V_a(\Z)\to\R$ and $f_2 \colon V_{-d}(\Z) \to \R$. Similarly, up to sets of arbitrarily small measure, we may tile $W_{a,d}(\R) \cap G_0(R^{(i)}_{\on{bal}})_X^{\on{bal}}$ with product sets in $W_{a,d}(\R)=V_a(\R) \times V_{-d}(\R)$ that cover it (resp., are covered by it). Hence, we may sum $\phi$ on elements $(A,B)\in W_{a,d}(\Z) \cap \gamma G_0\bigl(R_\bx^{(i)}\bigr)^\bx_X$ by independently summing functions over $A$ and $B$ in $\gamma$-translates of homogeneously expanding sets intersected with $V_a(\Z)$ and $V_{-d}(\Z)$, respectively. Applying Proposition \ref{prop-singcong}, we conclude that, up to an error of $E(X,\gamma) \defeq O(X^{2-\delta}h(\gamma)^M)$ for some $\delta, M > 0$, we have
\begin{equation*}
\sum_{w \in W_{a,d}(\Z)\cap \gamma G_0(R_\bx^{(i)})^\bx_X} \phi(w)=\sum_{A\in\frac{V_a(\Z)}{\SL_3(\Z)}}\sum_{B\in\frac{V_{-d}(\Z)}{\SL_3(\Z)}}\frac{\tau_{A,\infty}\tau_{B,\infty}}{\tau_{3,\infty}^2}\Vol\bigl(G_0(R^{(i)})^\bx_X\cap S_{AB}(\R)\bigr)\wt{\nu}_{A,B}(\phi),
\end{equation*}
where $\wt{\nu}_{A,B}(\phi)$ is defined as
\begin{equation*}
\wt{\nu}_{A,B}(\phi)\defeq
\Vol\big(S_{AB}(\widehat{\Z})\big)^{-1}
\nu_{A,B}(\phi);\quad\quad
\nu_{A,B}(\phi) \defeq 
\prod_p\int_{w\in S_{AB}(\Z_p)}\phi_p(w)\mathrm{d}w.
\end{equation*}
It is clear that $\delta$ can be chosen small enough so that the integral of $E(X,\gamma)$ over $\FF_3\smallsetminus \FF_3^{(\delta)}$ is $o(X^2)$. Note that $S_A(\Z_p)$ and $S_A(\R)$ only depend on the genus of $A$ and that $\Vol(S_A(\Z_p))=\Vol(\SL_3(\Z_p))/\Vol(\SO_A(\Z_p))$, with the analogous statements holding for $B$ as well. Thus, we can group forms $A$ and $B$ according to genus. Upon doing so and applying the Tamagawa number identities
\begin{equation}\label{eq:tam}
\tau_3 \defeq \prod_v\tau_{3,v}=1,\quad\quad
\sum_{A''\in\mc{G}(A')}\tau_{A''}=2,\quad\quad \tau_{A''} \defeq \prod_v \tau_{A'',v},
\end{equation}
where $\mc{G}(A')$ denotes the genus of $A'$ for any $A' \in V(\Z)$ of nonzero determinant, we find that
\begin{equation}\label{eq:4mainprt1}
N_{a,d}^\bx(W(\Z)^{(i)};X;\phi)\sim
\frac{4\tau_{3,\infty}}{\sigma_{(i)}\Vol(G_0)}\sum_{A\in \frac{V_a(\Z)}{\on{SL}_3(\Z)}} \sum_{B\in \frac{V_{-d}(\Z)}{\on{SL}_3(\Z)}}
\Vol\bigl(G_0(R^{(i)})^\bx_X\cap S_{AB}(\R)\bigr)\nu_{A,B}(\phi).
\end{equation}
Finally, we use the first Jacobian change-of-variables formula in Proposition \ref{prop-thejac} to compute
\begin{equation}\label{eq:4mainprt2}
\begin{array}{rcl}
\Vol\bigl(G_0(R^{(i)})^\bx_X\cap S_{AB}(\R)\bigr)&=&
\displaystyle|\J|\Vol(G_0)\int_{f\in U_{a,d}(\R)^{\pm,\bx}_X}\#\bigg(\frac{S_{AB}(\R)\cap\Res^{-1}(f)}{\SL_3(\R)}\bigg)\mathrm{d}f;
\\[.2in]
\displaystyle\int_{w\in S_{AB}(\Z_p)}\phi(w)dw
&=&\displaystyle
|\J|_p\tau_{3,p}\int_{f\in U_{a,d}(\Z_p)}\Bigl(\sum_{w\in\frac{S_{AB}(\Z_p)\cap\Res^{-1}(f)}{\SL_3(\Z_p)}}\frac{\phi(w)}{\#\Stab_{\SL_3(\Z_p)}(w)}
\Bigr)\mathrm{d}f.
\end{array}
\end{equation}
Theorem \ref{th:sec4main} for the case $\circ=\bx$ now follows by combining \eqref{eq:4mainprt1} and \eqref{eq:4mainprt2}.

\medskip

We now turn to the case $\circ = \on{wei}$. Since the proof is similar, we will be brief. We begin by writing
\begin{equation*}
N_{a,d}^{\wei}\big(W(\Z)^{(i)};X;\phi\big)\sim
\displaystyle\sum_{A\in\frac{V_a(\Z)}{\SL_3(\Z)}}
\frac{1}{\sigma_{(i)}\Vol(G_A)}\int_{\gamma\in\FF_A\smallsetminus\FF_A^{(\delta)}}\Bigl(\sum_{(A,B)\in  W_{a,d}(\Z)^{(i)} \cap \gamma G_A(R_\wei^{(i)})^\wei_X} \phi(A,B)\Bigr)\mathrm{d}\gamma.
\end{equation*}
As before, we apply Proposition \ref{prop-singcong} to conclude that, up to a negligible error, we have
\begin{equation*}
\sum_{(A,B)\in W_{a,d}(\Z)^{(i)} \cap \gamma G_A(R_\wei^{(i)})^\wei_X} \phi(A,B)=\sum_{B\in\frac{V_{-d}(\Z)}{\SL_3(\Z)}}\frac{\tau_{B,\infty}}{\tau_{3,\infty}}\Vol\bigl(G_A(R^{(i)})^\wei_X\cap S_{B}(\R)\bigr)\wt{\nu}_{A|B}(\phi),
\end{equation*}
where we abuse notation by identifying the set $G_A(R^{(i)})^\wei_X$ with its image under the projection map $W \to V$ defined by $(A,B) \mapsto B$, and where
\begin{equation*}
\wt{\nu}_{A|B}(\phi)\defeq
\Vol(S_{B}(\widehat{\Z}))^{-1}
\nu_{A|B}(\phi);\quad\quad
\nu_{A|B}(\phi) \defeq 
\prod_p\int_{B'\in S_{B}(\Z_p)}\phi_p(A,B')\mathrm{d}B'.
\end{equation*}
It is clear that $\nu_{\cA|B}$ (and hence $\widetilde{\nu}_{\cA|B}$ depends only on $\cA$ and the $\SL_3(\Z)$-class of $B$. Moreover, since $\phi$ is $\SL_3(\Z)$-invariant, it also follows that $\nu_{A|B}$ is unchanged when $A$ is replaced with an $\SL_3(\Z)$-translate. Therefore, we may again group $A$ and $B$ according to genus and write
\begin{equation}\label{eq:4mainprt3}
N_{a,d}^\wei\big(W(\Z)^{(i)};X;\phi\big)\sim
\sum_{A\in\overline{\mc{G}}_a}\frac{2\tau_{A,\infty}}{\sigma_{(i)}\Vol(G_A)}\sum_{B\in\overline{\mc{G}}_{-d}}
\Vol\bigl(G_A(R^{(i)})^\wei_X\cap S_{B}(\R)\bigr)\nu_{A|B}(\phi).
\end{equation}
Using the second Jacobian change-of-variables formula in Proposition \ref{prop-thejac}, we obtain
\begin{equation}\label{eq:4mainprt4}
\begin{array}{rcl}
\Vol\bigl(G_A(R^{(i)})^\bx_X\cap S_{B}(\R)\bigr)&=&
\displaystyle|\J_A|\Vol(G_A)\int_{f\in U_{a,d}(\R)^{\pm,\bx}_X}\#\frac{(\{A\}\times S_{B}(\R)^{(i)})\cap\Res^{-1}(f)}{\SO_A(\R)}\mathrm{d}f;
\\[.2in]
\displaystyle\int_{B'\in S_{B}(\Z_p)}\phi(A,B')dB'
&=&\displaystyle
|\J_A|_p\tau_{A,p}\int_{f\in U_{a,d}(\Z_p)}\Bigl(\sum_{(A,B')\in\frac{\{A\}\times S_{B}(\Z_p)\cap\Res^{-1}(f)}{\SO_A(\Z_p)}}\frac{\phi(w)}{\#\Stab_{\SO_A(\Z_p)}(w)}
\Bigr)\mathrm{d}f.
\end{array}
\end{equation}
Theorem \ref{th:sec4main} for the case $\circ = \on{wei}$ now follows by combining \eqref{eq:4mainprt3} and \eqref{eq:4mainprt4}. \hspace*{\fill} \qed

\section{Uniformity estimates and a square-free sieve}

%\subsection{The number of binary cubic forms in acceptable families}

Fix $a,d \in \Z \smallsetminus \{0\}$. For each prime $p$, let $\Sigma_p\subset U_{a,d}(\Z_p)$ be an open subset with measure-$0$ boundary such that $\Sigma_p \supset \{f\in U_{a,d}(\Z_p) : p^2\nmid\on{Disc}(f)\}$. Let $\Sigma\subset U_{a,d}(\Z)$ be the family corresponding to the collection $\{\Sigma_p\}_p$. In this section, we prove the following result.
\begin{theorem}\label{th:cubiccount}
Let $\circ \in \{\bx,\wei\}$. Then we have
\begin{equation*}
\#\Sigma^{\pm,\circ}_X\sim\Vol\bigl(U_{a,d}(\R)^{\pm,\circ}_X\bigr)\prod_p\Vol(\Sigma_p).
\end{equation*}
\end{theorem}

A straightforward application of geometry-of-numbers methods (for example, Davenport's lemma \cite{MR43821}) shows that Theorem~\ref{th:cubiccount} is true when $\Sigma$ is defined by finitely many congruence conditions (equivalently, $\Sigma_p=U_{a,d}(\Z_p)$ for sufficiently large $p$). The theorem then follows for any acceptable $\Sigma$ from an application of a standard sieve combined with the following tail estimates.

\begin{theorem}\label{th:unif}
For positive real numbers $M$ and $X$, we have
\begin{equation*}
\begin{array}{rcl}
\displaystyle\sum_{p>M}\#\bigl\{f\in U_{a,d}(\Z)^{\pm,\bx}_X:p^2\mid\on{Disc}(f)\bigr\}&\ll_{a,d,\varepsilon}&
\displaystyle\frac{X^{2+\varepsilon}}{M}+o(X^{2});
\\[.2in]
\displaystyle\sum_{p>M}\#\bigl\{f\in U_{a,d}(\Z)^{\pm,\wei}_X:p^2\mid\on{Disc}(f)\bigr\}&\ll_{a,d,\varepsilon}&
\displaystyle\frac{X^{3+\varepsilon}}{M}+o(X^{3}).
\end{array}
\end{equation*}
\end{theorem}
\begin{proof}
If $p^2\mid\on{Disc}(f)$ for $f\in U_{a,d}(\Z)$, then either $f$ has a triple root modulo $p$ or $R_f$ is non-maximal at $p$. Let $\W^{(1)}_p$ be the set of forms in $U_{a,d}(\Z)$ with nonzero discriminant having a triple root modulo $p$, and let $\W^{(2)}_p$ be the set of forms $f\in U_{a,d}(\Z)$ with nonzero discriminant such that $R_f$ is non-maximal at $p$. A direct application of \cite[Theorem 3.3]{geosieve} yields the estimates
\begin{equation*}
\sum_{p>M}
\#\bigl(\W_p^{(1)}\bigr)^{\bx}_X\ll_{a,d,\varepsilon}\frac{X^{2+\varepsilon}}{M}+X;\quad\quad
\sum_{p>M}
\#\bigl(\W_p^{(1)}\bigr)^{\wei}_X\ll_{a,d,\varepsilon}\frac{X^{3+\varepsilon}}{M}+X^{2}.
\end{equation*}
To obtain the analogous estimates for the case of $\W_p^{(2)}$, we break up the primes $p>M$ into two ranges, the ``small range'' and the ``large range,'' and we use different methods to handle each range.

\medskip 

\noindent {\bf The small range:}
Let $p\nmid ad$ be a fixed prime. We use the characterisation of elements $f\in \W_p^{(2)}$ from \cite[Lemma 3.5]{STpreprint}: namely, for a form $f\in U_{a,d}(\Z)$ with nonzero discriminant, the ring $R_f$ is non-maximal at $p$ if and only if there exists $\bar{r}\in\F_p$ such that for every lift $r\in\Z$ of $\bar{r}$, we have $p^2\mid f(r,1)$ and $p\mid f'(r,1)$, where $f'$ denotes the derivative of $f$ with respect to $x$. If we write forms $f\in \W_p^{(2)}$ as $f(x,y)=ax^3+bx^2y+cxy^2+dy^3$, then the value of $\bar{r}\in\F_p$ determines the reductions of $br+c$ and $2br+c$ modulo $p$, and hence determines the values of $b$ and $c$ modulo $p$. Moreover, once $b$ and $\bar{r}$ are determined, the reduction of $c$ modulo $p^2$ is determined. Therefore, we have
\begin{equation}\label{eq:smallrangep}
\#\bigl(\W_p^{(2)}\bigr)^{\bx}_X\ll p\Bigl(\frac{X}{p}+1\Bigr)\Bigl(\frac{X}{p^2}+1\Bigr);\quad\quad
\#\bigl(\W_p^{(2)}\bigr)^{\wei}_X\ll
p\Bigl(\frac{X}{p}+1\Bigr)\Bigl(\frac{X^2}{p^2}+1\Bigr).
\end{equation}

\medskip

\noindent {\bf The large range:} For large primes $p$, we follow the strategy in \cite{MR4277110} of fibering by $K_f\defeq \Q\otimes_{\Z} R_f$. First note that for primes $p\nmid a$, replacing $f(x,y)=ax^3+bx^2y+cxy^2+dy^3\in U_{a,d}(\Z)$ by $x^3+bx^2y+acxy^2+a^2dy^3$ does not change the condition $p^2\mid\on{Disc}(f)$, and it modifies both heights of $f$ by a factor of $O(1)$. Hence, we may assume without loss of generality that $a=1$. 

Consider the bijection
\begin{equation} \label{eq-1dbij}
U_{1,d}(\Z)^\irr \longleftrightarrow
\{(K,\alpha)\},
\end{equation}
between irreducible elements in $U_{1,d}(\Z)$ and the set of pairs $(K,\alpha)$, where $K$ is a cubic field and $\alpha\in\cO_K \smallsetminus \Z$ has norm $d$, given by sending $(K,\alpha)$ to the homogenisation of the minimal polynomial of $\alpha$. 
We have the following lemma concerning this bijection:
\begin{lemma} \label{lem-onbij}
Let $K$ be a cubic field. Then the number of elements $\alpha\in\cO_K\smallsetminus \Z$ corresponding to $f\in U_{1,d}(\Z)$ with $($either balanced or weighted$)$ height less than $X$ is $\ll_d \log X/\log|\Disc(K)|$.
\end{lemma}
\begin{proof}[Proof of Lemma~\ref{lem-onbij}]
Denote the length of $\alpha\in\cO_K$, with respect to the standard embedding $\cO_K\to \R\otimes_{\Q} K$, by $\ell(\alpha)$. It is easy to see that if $f\in U_{1,d}(\Z)^\irr$ corresponds to $(K,\alpha)$ under the bijection~\eqref{eq-1dbij}, then we have $\rmH_\bx(f)\geq \rmH_\wei(f)\asymp \ell(\alpha)$. Thus, for the purpose of proving the lemma, it is necessary for us to bound the number of elements $\alpha\in\cO_K$ with norm $d$ and length $\ell(\alpha)\ll X$. Let $S_K(d)$ denote the set of elements in $\cO_K$ with norm $d$. Then the group $\cO_K^\times$ acts on $S_K(d)$, and $\#(\cO_K^\times\backslash S_K(d))$ is bounded above by the number of principal ideals of norm $d$, which is $\leq 3^{\omega(d)}\ll_d 1$.

Recall the standard logarithmic embedding map $L\colon \cO_K^\times\to\R^r$, where $r$ denotes the rank of $\mc{O}_K^\times$. 
Every unit $\alpha\in\cO_K^\times\smallsetminus\{\pm 1\}$ (indeed, every element $\alpha\in\cO_K\smallsetminus\Z$) must have length $\ell(\alpha)\gg |\Disc(K)|^{1/6}$, since the covolume of $\Z[\alpha]$ in $\R \otimes_{\Q} K$ is at least $\sqrt{|\Disc(K)|}$. Thus, the length of $L(\alpha)$ must be $\gg \log |\Disc(K)|$ for every $\alpha \in \mc{O}_K^\times \smallsetminus \{\pm1\}$. This implies that every $\cO_K^\times$-orbit in $S_K(d)$ has $\ll \log X/\log|\Disc(K)|$ elements in it corresponding to polynomials with height less than $X$, as necessary.
\end{proof}

Now observe that for $p>M$, any $f\in (\W_p^{(2)})^{\bx}_X$ (resp., any $f\in (\W_p^{(2)})^{\wei}_X$) corresponds to 
a pair $(K,\alpha)$ with $|\Disc(K)|\ll X^4/M^2$ (resp., $|\Disc(K)|\ll X^6/M^2$). Moreover, the Davenport--Heilbronn theorem (see, e.g.,~\cite{MR3090184}) implies that the number of such cubic fields is $\ll X^4/M^2$ (resp., $\ll X^6/M^2$). Hence, we have
\begin{equation}\label{eq:largerangep}
\sum_{p>M}\#\bigl(\W_p^{(2)}\bigr)^{\bx}_X\ll \frac{X^4}{M^2};\quad\quad
\sum_{p>M}\#\bigl(\W_p^{(2)}\bigr)^{\wei}_X\ll
\frac{X^6}{M^2}.
\end{equation}
Theorem \ref{th:unif} now follows for the balanced height by using \eqref{eq:smallrangep} when $p\leq X\log\log X$ and \eqref{eq:largerangep} when $p>X\log\log X$, and for the weighted height by using \eqref{eq:smallrangep} when $p\leq X^{3/2+\delta}$ and \eqref{eq:largerangep} when $p>X^{3/2+\delta}$ for some suitably small $\delta > 0$.
\end{proof}

Next, we have the following immediate consequence of Theorem \ref{th:sec4main}.

\begin{theorem}\label{th:sec4maininf}
Let $\circ \in \{\on{bal},\on{wei}\}$, let $(i)$ be a splitting type, and let $\phi \colon W_{a,d}(\Z) \to [0,1]$ be an $\SL_3(\Z)$-invariant function defined by $($possibly$)$ infinitely many congruence conditions. Then we have
\begin{equation*}
\begin{array}{rcl}
\displaystyle\sum_{\substack{w\in \frac{W(\Z)^{(i),\circ}_X}{\SL_3(\Z)}\\w\text{ \rm non-dist.}}}\phi(w)&\leq&
\displaystyle 4\sum_{\substack{A\in \mc{G}_a\\B\in \mc{G}_{-d}}}\int_{f\in U_{a,d}(\R)^{\pm,\circ}_X}\on{Mass}_\infty\big(f;A,B;(i)\big)\mathrm{d}f\cdot\prod_p\int_{f\in U_{a,d}(\Z_p)}\on{Mass}_p(f;A,B;\phi)\mathrm{d}f . 
\\[.2in]&&+\displaystyle
O(\Vol(U_{a,d}(\R)_X^{\pm,\circ})).
\end{array}
\end{equation*}
\end{theorem}

\noindent The above result follows by writing $\phi$ as a product $\phi=\prod_p \phi_p$ of functions $\phi_p:W_{a,d}(\Z_p)\to [0,1]$, applying Theorem \ref{th:sec4main} on the function $\prod_{p<M}\phi_p$, and then letting $M$ go to infinity.

\begin{remark} \label{rmk:unif estimate}
Suppose that for infinitely many primes $p$, the function $\phi_p$ is $1$ on the set of elements $w\in W_{a,d}(\Z_p)$ with $p^2\nmid\Disc(w)\defeq \Disc(\Res(w))$. Then a tail estimate (analogous to Theorem \ref{th:unif} on the number of elements
$w\in \frac{W(\Z)^{(i),\circ}_X}{\SL_3(\Z)}$
with $p^2\mid\Disc(w)$) would allow us to replace the ``$\leq$'' in Theorem \ref{th:sec4maininf} with an equality. Though we expect such an bound to be true, we are not yet able to prove it. This is the reason we have upper bounds for our main results Theorems \ref{first main}--\ref{third main2}, and why we expect that this upper bound is also the true lower bound.
\end{remark}

\section{Average $2$-torsion in the class groups of cubic fields}\label{sec-proofofmainresults}

In this section, we prove our main results on class groups of cubic fields in thin families, namely Theorems \ref{first main}, \ref{first main2}, \ref{second main real}, \ref{second main complex}, and \ref{main stable}. We start in \S\ref{sec-masses} by evaluating the relevant local masses, and we complete the proofs in \S\ref{sec-thatsit}.

\subsection{Evaluating the local masses}  \label{sec-masses}

The purpose of this section is to compute the local masses $\on{Mass}_\infty(f;A,B;(i))$ and $\on{Mass}_p(f;A,B;\phi)$ that appear in Theorem~\ref{th:sec4main} in the case when $f$ is maximal, $\phi=1$, and $(i)=(1111)$ or $(i)=(112)$. To facilitate this computation, we start by recalling from~\cite{BHSpreprint} some notation and results concerning ternary quadratic forms over $\R$ and $\Z_p$.

\begin{defn}[\protect{Adapted from \cite[Definition~4.23]{BHSpreprint}}] \label{def-kappas}
Let $v$ be a place of $\Q$, and let $A \in V(\Z_v)$. We define a quantity $\kappa_v(A) \in \{0,\pm1\}$ as follows:
\begin{itemize}
\item[{\rm (a)}] Suppose $v = \infty$. If $A$ is isotropic over $\R$, then $\kappa_\infty(A) = 1$; otherwise, $\kappa_\infty(A) = -1$.
\item[{\rm (b)}] Suppose $v = p$ is a prime, and further suppose that the mod-$p$ reduction $\ol{A}$ of $A$ is either smooth or the union of two distinct lines.\footnote{To be clear, $\kappa_p(A)$ is only defined for those $A$ such that $\ol{A}$ has one of these two reduction types.} If $\ol{A}$ is a union of two distinct lines, and if these lines are defined over $\mathbb{F}_p$, then $\kappa_p(A) = 1$; if the lines are conjugates over $\mathbb{F}_{p^2}$, then $\kappa_p(A) = -1$. Otherwise, $\kappa_p(A) = 0$.
\end{itemize}
\end{defn}

  We now recall two lemmas from~\cite{BHSpreprint} concerning the quantity $\kappa_v$ introduced in Definition~\ref{def-kappas}. The first lemma describes the orbits of the action of $\on{SL}_3(\Z_v)$ on the set of ternary quadratic forms having specified determinant and value of $\kappa_v$:

\begin{lemma}[\protect{\cite[Lemma~4.24]{BHSpreprint}}] \label{lem-kappa1}
Let $v$ be a place of $\Q$, let $a \in \Z \smallsetminus \{0\}$, and let $\varepsilon \in \{0,\pm 1\}$. Then the set $\{A \in V_a(\Z_v) : \kappa_v(A) = \varepsilon\}$ is nonempty if and only if $v = \infty$ or $v = p$ is a prime and $\nu_p(a) = \varepsilon = 0$ or $\nu_p(a) \geq |\varepsilon| > 0$; moreover, if this set is nonempty, it consists of a single $\on{SL}_3(\Z_v)$-orbit.
\end{lemma}

The second lemma demonstrates that if $(A,B)$ is a pair of ternary quadratic forms over $\Z_p$ whose resolvent is maximal, then the quantities $\kappa_p(A)$ and $\kappa_p(B)$ are defined:

\begin{lemma}[\protect{\cite[Lemma~4.24]{BHSpreprint}}] \label{lem-kappa2}
Let $f \in U(\Z_p)$ be maximal, and let $(A,B) \in W(\Z_p)$ with $\on{Res}(A,B) = f$. Then each of $\ol{A}$ and $\ol{B}$ is either smooth or the union of two distinct lines.
\end{lemma}

We are now in position to compute the desired local masses. For each $f \in U(\Z_v)$, maximal if $v = p$ is a prime, and each $\varepsilon_1, \varepsilon_2 \in \{0,\pm 1\}$, define the following mass:
\begin{align*}
\on{Mass}_\infty^{\varepsilon_1, \varepsilon_2}(f; (i)) & \defeq \sum_{\substack{(A,B) \in \frac{W(\R)^{(i)} \cap \on{Res}^{-1}(f)}{\on{SL}_3(\R)} \\ \kappa_\infty(A) = \varepsilon_1,\, \kappa_\infty(B) = \varepsilon_2}} \frac{1}{\#\on{Stab}_{\on{SL}_3(\R)}(A,B)}; \\
\on{Mass}_p^{\varepsilon_1, \varepsilon_2}(f) & \defeq \sum_{\substack{(A,B) \in \frac{W(\Z_p) \cap \on{Res}^{-1}(f)}{\on{SL}_3(\Z_p)} \\ \kappa_p(A) = \varepsilon_1,\, \kappa_p(B) = \varepsilon_2}} \frac{1}{\#\on{Stab}_{\on{SL}_3(\Z_p)}(A,B)}.
\end{align*}
As a consequence of Lemmas~\ref{lem-kappa1} and~\ref{lem-kappa2}, the masses $\on{Mass}_\infty(f; A; B; (i))$ with $i \in \{(1111),(112)\}$ and $\on{Mass}_p(f;A,B;\phi)$ are respectively of the form $\on{Mass}_\infty^{\varepsilon_1, \varepsilon_2}(f; (i))$ and $\on{Mass}_p^{\varepsilon_1, \varepsilon_2}(f)$ for some $\varepsilon_1, \varepsilon_2 \in \{0,\pm1\}$. Thus, to apply Theorem~\ref{th:sec4main}, it suffices to evaluate these new masses for each choice of $\varepsilon_1,\varepsilon_2 \in \{0,\pm1\}$, $(i) \in \{(1111), (112)\}$, and $f \in U(\Z_v)$, maximal if $v = p$ is a prime. We achieve this in the next pair of propositions; we omit the proof of the first as it is evident.

\begin{proposition}
Let $f \in U(\R)^+$. Then we have
\begin{equation*}
\on{Mass}_\infty^{1,1}\big(f;(1111)\big)=\frac14;\quad \on{Mass}_\infty^{-1,\pm1}\big(f;(1111)\big)=\on{Mass}_\infty^{\pm1,-1}\big(f;(1111)\big)=0.
\end{equation*}
Let $f \in U(\R)^-$. Then we have
\begin{equation*}
\on{Mass}_\infty^{1,1}\big(f;(112)\big)=\frac12;\quad \on{Mass}_\infty^{-1,\pm1}\big(f;(112)\big)=\on{Mass}_\infty^{\pm1,-1}\big(f;(112)\big)=0.
\end{equation*}
\end{proposition}

\begin{proposition} \label{prop-allthemasses}
Let $f(x,y) = ax^3 + bx^2y + cxy^2 + dy^3 \in U_{a,d}(\Z_p)_{\max}$. Then we have
\begin{equation} \label{eq-allthemasses}
\begin{array}{ll}
\on{Mass}_p^{0,0}(f) = 1, & \text{if $p \nmid a$ and $p \nmid d$,} \\[0.1cm] 
\on{Mass}_p^{\pm1,0}(f) = \frac{1}{2}, & \text{if $p \mid a$ and $p \nmid d$, and $f$ is not sufficiently-ramified,}\\[0.1cm]
\on{Mass}_p^{1,0}(f) = 1, & \text{if $p \mid a$ and $p \nmid d$, and $f$ is left-sufficiently-ramified,}\\[0.1cm]
\on{Mass}_p^{0,\pm1}(f) = \frac{1}{2}, & \text{if $p \nmid a$ and $p \mid d$, and $f$ is not sufficiently-ramified,}\\[0.1cm]
\on{Mass}_p^{0,1}(f) = 1, & \text{if $p \nmid a$ and $p \mid d$, and $f$ is right-sufficiently-ramified,}\\[0.1cm]
\on{Mass}_p^{\pm1,\pm1}(f) = \frac{1}{4}, & \text{if $p \mid a$ and $p \mid d$, and $f$ is not sufficiently-ramified,} \\[0.1cm]
\on{Mass}_p^{\pm 1,1}(f) = \frac{1}{2},& \text{if $p \mid a$ and $p \mid d$, and $f$ is left-sufficiently-ramified,}\\[0.1cm]
\on{Mass}_p^{1,\pm 1}(f) = \frac{1}{2},& \text{if $p \mid a$ and $p \mid d$, and $f$ is right-sufficiently-ramified.}\\[0.1cm]
\end{array}
\end{equation}
If $\varepsilon_1,\, \varepsilon_2,\, f$ are such that the symbol $\on{Mass}_p^{\varepsilon_1,\varepsilon_2}(f)$ does not appear above, then $\on{Mass}_p^{\varepsilon_1,\varepsilon_2}(f) =0$.
\end{proposition}
\begin{proof}
When $p \nmid a$ and $p \nmid d$, we have $\kappa_p(A) = \kappa_p(B) = 0$ for every $(A,B) \in \on{Res}^{-1}(f)$. Thus, in this case we have
\begin{equation} \label{eq-00mass}
\on{Mass}_p^{0,0}(f) = \sum_{(A,B) \in \frac{W(\Z_p) \cap \on{Res}^{-1}(f)}{\on{SL}_3(\Z_p)}} \frac{1}{\#\on{Stab}_{\on{SL}_3(\Z_p)}(A,B)},
\end{equation}
and the sum on the right-hand side of~\eqref{eq-00mass} was shown in~\cite[Corollary~3.40]{BHSpreprint} to be as claimed in~\eqref{eq-allthemasses}.

When $p \mid a$ and $p \nmid d$, we have $\kappa_p(A) \in \{\pm 1\}$ and $\kappa_p(B) = 0$. Thus, in this case we have
\begin{equation} \label{eq-10mass}
\on{Mass}_p^{\pm 1,0}(f) = \sum_{\substack{(A,B) \in \frac{W(\Z_p) \cap \on{Res}^{-1}(f)}{\on{SL}_3(\Z_p)} \\ \kappa_p(A) = \pm 1}} \frac{1}{\#\on{Stab}_{\on{SL}_3(\Z_p)}(A,B)},
\end{equation}
and the sum on the right-hand side of~\eqref{eq-10mass} was shown in~\cite[Corollary~4.26]{BHSpreprint} to be as claimed in~\eqref{eq-allthemasses}. The case where $p \nmid a$ and $p \mid d$ follows by symmetry.

It remains to handle the case where $p \mid a$ and $p \mid d$, where we have $\kappa_p(A), \kappa_p(B) \in \{\pm 1\}$. In this case it was shown in~\cite[Corollary~4.26]{BHSpreprint} that
\begin{equation} \label{eq-totmass}
\on{Mass}_p^{\pm 1, 1}(f) + \on{Mass}_p^{\pm 1, -1}(f) = \begin{cases} \frac{1}{2}, & \text{if $f$ is not right-sufficiently-ramified,} \\ \frac{1}{2}(1 \pm 1), & \text{if $f$ is right-sufficiently-ramified.} \end{cases}
\end{equation}
Given~\eqref{eq-totmass}, it suffices to determine the quantity
$$\on{Orb}_p^{\varepsilon_1, \varepsilon_2}(f) \defeq \#\left(\frac{\{(A,B) \in \on{Res}^{-1}(f) : \kappa_p(A) = \varepsilon_1, \kappa_p(B) = \varepsilon_2\}}{\on{SL}_3(\Z_p)}\right)$$
for all $\varepsilon_1, \varepsilon_2 \in \{\pm 1\}$, because $\#\on{Stab}_{\on{SL}_3(\Z_p)}(A,B)$ is constant over all $(A,B) \in \on{Res}^{-1}(f)$. The following result from~\cite{BHSpreprint} reduces this local orbit-counting problem into the even simpler problem of counting quadratic extensions of cubic $\Q_p$-algebras:

\begin{theorem}[\protect{\cite[Theorem~2.22 and proof of Corollary~4.26]{BHSpreprint}}] \label{thm-zpparam}
Let $p$ be a prime, and let $f(x,y) = ax^3 + bx^2y + cxy^2 + dy^3 \in U_{a,d}(\Z_p)_{\max}$ be a binary cubic form. 
\begin{itemize}
\item[{\rm (a)}] The set of $\on{SL}_3(\Z_p)$-orbits on $W(\Z_p)$ with resolvent $f$ is in bijection with the set of unramified quadratic extensions $L/K_f$ such that norm to $\Q_p$ of the discriminant $\Delta(L/K_f)$ is a square.
\item[{\rm (b)}] Suppose $p \mid a$, and let $f$ factor over $\Q_p$ as $f = f_1 \times f_2$, where $f_1$ is a power of $y$ modulo $p$, and where $f_2$ is indivisible by $y$ modulo $p$. Take $(A,B) \in W(\Z_p)$ with $\on{Res}(A,B) = f$, let $L$ be the unramified quadratic extension of $K_f$ corresponding to the $\on{SL}_3(\Z_p)$-orbit of $(A,B)$ under the bijection in part ${\rm (a)}$; write $L = L_1 \times L_2$ corresponding to the decomposition $K_f = K_{f_1} \times K_{f_2}$. Then $\kappa_p(A) = 1$ $($resp., $\kappa_p(A) = -1${}$)$ if and only if $L_1/K_{f_1}$ is split $($resp., inert$)$. The corresponding statement holds with $a$ replaced by $d$, $y$ replaced by $x$, and $A$ replaced by $B$.
\end{itemize}
\end{theorem}

By Theorem~\ref{thm-zpparam}, determining the values of $\on{Orb}_p^{\varepsilon_1,\varepsilon_2}(f)$ amounts to counting the unramified quadratic extensions of $K_f$ with specified splitting type whose discriminants have square norm to $\Q_p$. To this end, let $f$ factor over $\Q_p$ as $f = f_\ell \times f_c \times f_r$, where $f_\ell$ is a power of $y$ modulo $p$, $f_c$ is coprime to $xy$ modulo $p$, and $f_r$ is a power of $x$ modulo $p$. We split into cases depending on whether $f$ is sufficiently-ramified:
\begin{itemize}
\item First suppose $f$ is not sufficiently-ramified. Then $f_\ell$, $f_c$, and $f_r$ are linear, so $K_{f_\ell}$, $K_{f_c}$, and $K_{f_r}$ are isomorphic to $\Q_p$. Thus, every unramified quadratic extension $L/K_f$ can be written as $L \simeq L_\ell \times L_c \times L_r$, where $L_\bullet$ is an unramified quadratic extension of $K_\bullet \simeq \Q_p$ for each $\bullet \in \{\ell,c,r\}$. By part (b) of Theorem~\ref{thm-zpparam}, the values of $\varepsilon_1$ and $\varepsilon_2$ determine whether $L_\ell$ and $L_r$ are split or inert. Moreover, since the discriminant of an unramified quadratic extension has square norm to $\Q_p$ if and only if it is split, the splitting type of $L_c$ is determined by the splitting types of $L_\ell$ and $L_r$. Thus, $\on{Orb}_p^{\varepsilon_1,\varepsilon_2}(f) = 1$ for each $\varepsilon_1,\varepsilon_2 \in \{\pm 1\}$.
\item Next suppose $f$ is right-sufficiently-ramified. Then $f_\ell$ is linear, $f_c$ is a constant, and $f_r$ is quadratic, so $K_{f_\ell} \simeq \Q_p$, $K_{f_c} \simeq 0$, and $K_{f_r}$ is a ramified quadratic extension of $\Q_p$. Thus, every unramified quadratic extension $L/K_f$ can be written as $L \simeq L_\ell \times L_r$, where $L_\ell$ is an unramified quadratic extension of $\Q_p$ and $L_r$ is an unramified quadratic extension of $K_{f_r}$. Moreover, since the discriminant of an unramified quadratic extension of $\Q_p$ has square norm to $\Q_p$ if and only if it is split, and since the discriminant of an unramified quadratic extension of a ramified quadratic extension of $\Q_p$ always has square norm to $\Q_p$, the extension $L_\ell/K_{f_\ell}$ must be split. Thus, $\on{Orb}_p^{1,\pm 1}(f) = 1$ and $\on{Orb}_p^{-1,\pm 1}(f) = 0$.

The case where $f$ is left-sufficiently-ramified follows by symmetry.
\end{itemize}
This completes the proof of Proposition~\ref{prop-allthemasses}.
\end{proof}

\subsection{Proofs of Theorems~\ref{first main},~\ref{first main2},~\ref{second main real}, \ref{second main complex}, and \ref{main stable}} \label{sec-thatsit}

We start by proving Theorems~\ref{second main real} and~\ref{second main complex}. As in the theorem statements, fix $\varepsilon \in \{\pm\}$ and $\Sigma = U_{a,d}(\Z)_{\max}^\varepsilon \,\cap\, \bigcap_p \Sigma_p$ to be an acceptable family. 

Fix $\circ \in \{\on{bal},\on{wei}\}$. Suppose $\varepsilon  = +$, and for each prime $p$, let $\phi_p$ be the characteristic function of the set $\on{Res}^{-1}(\Sigma_p) \cap W_{a,d}(\Z_p)$. Let $\mathfrak{d}_\Sigma$ denote the density of $f \in \Sigma$ above which the $\Delta$-distinguished orbit contains an integral representative, and let $\mathfrak{d}_\infty^{(0)}$ be as in Proposition~\ref{prop-realvols}. Then combining Theorems~\ref{thm-dual2torsparam},~\ref{thm-howmanydist},~\ref{th:sec4main}, and~\ref{th:cubiccount} yields that the average $2$-torsion in the class group of the maximal cubic rings $R_f$ cut out by forms $f \in \Sigma$ is bounded as follows:
\begin{align}
    \limsup_{X \to \infty} \underset{f \in \Sigma_X^\circ}{\on{Avg}}\,\on{Cl}(R_f)[2] & \leq 1 + \mathfrak{d}_\infty^{(0)}\mathfrak{d}_\Sigma +  \sum_{\substack{A \in \mc{G}_a,\, B \in \mc{G}_{-d} \\ (A,B) \in W(\R)^{(0)}}} \prod_p\frac{1}{ \on{Vol}(\Sigma_p)}\int_{f \in \Sigma_p} \on{Mass}_p^{\kappa_p(A),\kappa_p(B)}(f) \mathrm{d}f,
    \label{eq-firsttry}
    \end{align}
    where the ``$1$'' term on the right-hand side constitutes the contribution from the distinguished orbit. We obtain an inequality rather than an equality because we are imposing infinitely many congruence conditions; proving a tail estimate for $W(\Z)$ akin to Theorem~\ref{th:unif} would allow one to prove the reverse inequality too.

    For $m \in \Z \smallsetminus\{0\}$, let $T_m$ be the set of primes dividing $m$, and let $T_m^{\on{odd}} \defeq \{p \in T_m : \nu_p(m) \equiv 1 \pmod 2\}$ and $T_m^{\on{even}} \defeq T_m \smallsetminus T_m^{\on{odd}}$. To evaluate the above sum, we use the following lemma, which states that, for each $m \in \Z \smallsetminus \{0\}$, the elements of $\mc{G}_m$ are parametrised by the $\kappa_v$-values defined in Definition~\ref{def-kappas}:
    \begin{lemma}[\protect{\cite[Proof of Theorem 7, p.~36]{BHSpreprint}}] \label{lem-kappagenerate}
The map sending $M \in \mc{G}_m$ to the tuple $(\kappa_v(M))_v \in \prod_v \{0,\pm1\}$ defines a bijection from the set $\mc{G}_m$ to the set of tuples $(\varepsilon_v)_v \in \prod_v \{0,\pm1\}$ over places $v$ of $\Q$, where $\varepsilon_p = 0$ if and only if $p \nmid m$ and $\varepsilon_\infty \times \prod_{p \in T_m^{\on{odd}}} \varepsilon_p = 1$.
\end{lemma}
\noindent 
Applying Lemma~\ref{lem-kappagenerate} to the right-hand side of~\eqref{eq-firsttry}, we find that
\begin{align}
\limsup_{X \to \infty} \underset{f \in \Sigma_X^\circ}{\on{Avg}}\,\on{Cl}(R_f)[2]  & \leq 
    1 +  \mathfrak{d}_\infty^{(0)}\mathfrak{d}_\Sigma + \sum_{\substack{(\varepsilon_{1,p},\varepsilon_{2,p})_p \in \prod_p \{0,\pm1\}^2 \\ \prod_{p \in T_a^{\on{odd}}} \varepsilon_{1,p}=\prod_{p \in T_d^{\on{odd}}} \varepsilon_{2,p}  = 1 \\ \prod_{p \in T_a^{\on{even}}} |\varepsilon_{1,p}|=\prod_{p \in T_d^{\on{even}}} |\varepsilon_{2,p}| = 1}} \prod_p\frac{1}{ \on{Vol}(\Sigma_p)}\int_{f \in \Sigma_p} \on{Mass}_p^{\varepsilon_{1,p},\varepsilon_{2,p}}(f)\mathrm{d}f. \label{eq-withgenus}
\end{align}
Using Proposition~\ref{prop-allthemasses}, we find that the factor at $p$ in~\eqref{eq-withgenus} is given as follows:
\begin{equation} \label{eq-massinrhos}
   \frac{1}{ \on{Vol}(\Sigma_p)}\int_{f \in \Sigma_p}  \on{Mass}_p^{\varepsilon_{1,p},\varepsilon_{2,p}}(f)\mathrm{d}f = \begin{cases} \displaystyle 1 , & \text{if $p \not\in (T_a \cup T_d)$,} \\[0.2cm]
   \displaystyle \frac{1}{2} + \varepsilon_{1,p}\frac{\rhoright_\Sigma(p)}{2}, & \text{if $p \in T_a \smallsetminus T_d$,}  \\[0.2cm]
    \displaystyle \frac{1}{2} + \varepsilon_{2,p}\frac{\rholeft_\Sigma(p)}{2}, & \text{if $p \in T_d \smallsetminus T_a$,}  \\[0.2cm]
      \displaystyle \frac{1}{4} + \varepsilon_{1,p}\frac{\rhoright_\Sigma(p)}{4}+ \varepsilon_{2,p}\frac{\rholeft_\Sigma(p)}{4}, & \text{if $p \in T_a \cap T_d$.}
   \end{cases}
\end{equation}
Substituting the result of~\eqref{eq-massinrhos} into~\eqref{eq-withgenus}, we deduce that
\begin{align}
    & \limsup_{X \to \infty} \underset{f \in \Sigma_X^\circ}{\on{Avg}}\,\on{Cl}(R_f)[2]  \leq 1 + \mathfrak{d}_\infty^{(0)}\mathfrak{d}_\Sigma +  \sum_{\substack{(\varepsilon_{1,p})_p\in \{\pm1\}^{T_a^{\on{odd}}}, (\varepsilon_{2,p})_p \in \{\pm 1\}^{T_d^{\on{odd}}} \\ \prod_{p \in T_a^{\on{odd}}} \varepsilon_{1,p}= \prod_{p \in T_d^{\on{odd}}} \varepsilon_{2,p} = 1}} \prod_{p \in T_a^{\on{odd}} \smallsetminus T_d^{\on{odd}}} \bigg(\frac{1}{2} + \varepsilon_{1,p}\frac{\rhoright_\Sigma(p)}{2}\bigg) \times \label{eq-longprod} \\
   & \qquad\qquad\qquad\qquad\qquad \prod_{p \in T_d^{\on{odd}} \smallsetminus T_a^{\on{odd}}} \bigg(\frac{1}{2} + \varepsilon_{2,p}\frac{\rholeft_\Sigma(p)}{2}\bigg) \prod_{p \in T_a^{\on{odd}} \cap T_d^{\on{odd}}} \bigg(\frac{1}{4} + \varepsilon_{1,p}\frac{\rhoright_\Sigma(p)}{4}+ \varepsilon_{2,p}\frac{\rholeft_\Sigma(p)}{4}\bigg) \nonumber
\end{align}
To simplify the right-hand side of~\eqref{eq-longprod}, we make use of the following identity:
\begin{lemma}[\protect{\cite[Lemma~4.17]{hankerhoproduct}}] \label{lem-hankeprod}
Suppose $\mathbb{T}$ is a nonempty finite set and that $X_i$ and $Y_i$ are indeterminates for each $i \in \mathbb{T}$. Then for any $c \in \{\pm 1\}$, we have the polynomial identity
$$\sum_{\substack{(\varepsilon_i)_i \in \{\pm 1\}^{\mathbb{T}} \\ \prod_{i \in \mathbb{T}} \varepsilon_i = c}}\prod_{i \in \mathbb{T}}(X_i + \varepsilon_iY_i) = 2^{\#\mathbb{T}-1}\times \bigg(\prod_{i \in \mathbb{T}} X_i + c \prod_{i \in \mathbb{T}} Y_i\bigg).$$
\end{lemma}
\noindent A straightforward application of Lemma~\ref{lem-hankeprod} to~\eqref{eq-longprod} yields that
\begin{equation*}
    \limsup_{X \to \infty} \underset{f \in \Sigma_X^\circ}{\on{Avg}}\,\on{Cl}(R_f)[2]  \leq \frac{5}{4} + \frac{\rhoright_\Sigma + \rholeft_\Sigma + \chi_{a,d}\, \rhoright_\Sigma \rholeft_\Sigma}{4} + \mathfrak{d}_\infty^{(0)}\mathfrak{d}_\Sigma  ,
\end{equation*}
which is precisely as claimed in Theorem~\ref{second main real}.

\medskip

Next, suppose $\varepsilon = -$, and let $\phi$ be as above. Let $\mathfrak{d}_\Sigma$ denote the density of $f \in \Sigma$ above which the $\Delta$-distinguished orbit contains an integral representative. Then combining Theorems~\ref{thm-dual2torsparam},~\ref{thm-howmanydist},~\ref{th:sec4main}, and~\ref{th:cubiccount}, along with Proposition~\ref{prop-coincidenot}, Lemmas~\ref{lem-kappagenerate} and~\ref{lem-hankeprod} and the formula~\eqref{eq-massinrhos}, yields that the average $2$-torsion in the class group of the maximal cubic rings $R_f$ cut out by forms $f \in \Sigma$ is bounded as follows:
\begin{align}
    \limsup_{X \to \infty} \underset{f \in \Sigma_X^\circ}{\on{Avg}}\,\on{Cl}(R_f)[2]   & 
    \leq 1 + \mathfrak{d}_\Sigma + 2 \times\sum_{\substack{A \in \mc{G}_a,\, B \in \mc{G}_{-d} \\ (A,B) \in W(\R)^{(1)}}} \prod_p\frac{1}{ \on{Vol}(\Sigma_p)}\int_{f \in \Sigma_p} \on{Mass}_p^{\kappa_p(A),\kappa_p(B)}(f) \mathrm{d}f\nonumber \\
    & = 1 + \mathfrak{d}_\Sigma + 2\times \sum_{\substack{(\varepsilon_{1,p},\varepsilon_{2,p})_p \in \prod_p \{0,\pm1\}^2 \\ \prod_{p \in T_a^{\on{odd}}} \varepsilon_{1,p}=\prod_{p \in T_d^{\on{odd}}} \varepsilon_{2,p}  = 1 \\ \prod_{p \in T_a^{\on{even}}} |\varepsilon_{1,p}|=\prod_{p \in T_d^{\on{even}}} |\varepsilon_{2,p}| = 1}} \prod_p\frac{1}{ \on{Vol}(\Sigma_p)}\int_{f \in \Sigma_p} \on{Mass}_p^{\varepsilon_{1,p},\varepsilon_{2,p}}(f)\mathrm{d}f \nonumber \\
    & = 1 +\mathfrak{d}_\Sigma + 2 \times \sum_{\substack{(\varepsilon_{1,p})_p\in \{\pm1\}^{T_a^{\on{odd}}}, (\varepsilon_{2,p})_p \in \{\pm 1\}^{T_d^{\on{odd}}} \\ \prod_{p \in T_a^{\on{odd}}} \varepsilon_{1,p}= \prod_{p \in T_d^{\on{odd}}} \varepsilon_{2,p} = 1}} \prod_{p \in T_a^{\on{odd}} \smallsetminus T_d^{\on{odd}}} \bigg(\frac{1}{2} + \varepsilon_{1,p}\frac{\rhoright_\Sigma(p)}{2}\bigg) \times \nonumber \\
   & \qquad\qquad \prod_{p \in T_d^{\on{odd}} \smallsetminus T_a^{\on{odd}}} \bigg(\frac{1}{2} + \varepsilon_{2,p}\frac{\rholeft_\Sigma(p)}{2}\bigg) \prod_{p \in T_a^{\on{odd}} \cap T_d^{\on{odd}}} \bigg(\frac{1}{4} + \varepsilon_{1,p}\frac{\rhoright_\Sigma(p)}{4}+ \varepsilon_{2,p}\frac{\rholeft_\Sigma(p)}{4}\bigg)  \nonumber \\
   & = \frac{3}{2}  + \frac{\rhoright_\Sigma + \rholeft_\Sigma + \chi_{a,d}\, \rhoright_\Sigma \rholeft_\Sigma}{2}+ \mathfrak{d}_\Sigma, \nonumber
\end{align}
which is precisely as claimed in Theorem~\ref{second main complex}. In the above, the ``$1$'' (resp., ``$\mathfrak{d}_\Sigma$'') term on the right-hand side constitutes the contribution from the distinguished (resp., $\Delta$-distinguished) orbit. 

\medskip

Theorem \ref{main stable} follows from Theorems~\ref{second main real} and~\ref{second main complex} by setting/redefining each of $\rho_\Sigma$, $\lambda_\Sigma$, $\chi_{a,d}$, $\mathfrak{d}_\Sigma$, and $\mathfrak{d}_\infty^{(0)}$ to be $0$ or $1$ according to the parameters defining the stable family of interest. 
As for Theorems~\ref{first main} and~\ref{first main2}, these follow by applying Theorems~\ref{second main real} and~\ref{second main complex} to the case where $\Sigma_p = U_{a,d}(\Z_p)_{\max}$ for each prime $p$. The relevant local densities are computed in Appendix~\ref{appendix-localcomputationsat2}; see Proposition~\ref{prop-maxcalcs} for the sufficient-ramification densities and Corollary~\ref{cor-deltadist} for the density $\mathfrak{d}_\Sigma$. 

\section{Average size of the $2$-Selmer groups of elliptic curves} \label{sec-selproofs}

In this section, we prove our main results on Selmer groups of elliptic curves in thin families, namely Theorems \ref{third main} and \ref{third main2}. We start in \S\ref{sec-selmasses} by evaluating the relevant local masses, and we complete the proofs in \S\ref{sec-selthatsit}.

\subsection{Evaluating the local masses}  \label{sec-selmasses}

Let $ d\in \Z \smallsetminus \{0\}$ such that $d \equiv 1 \pmod 8$. The purpose of this section is to compute the local masses $\on{Mass}_\infty(4f; A; B; (i))$ and $\on{Mass}_p(4f; A; B; \phi)$ that appear in Theorem~\ref{th:sec4main} in the case when $f$ is maximal, $\phi$ is the characteristic function of locally soluble elements in $W_{1,d}^\vee(\Z)$, and $(i) = (1111)$, $(i) = (2+)$, or $(i) = (112)$.

For a place $v$ of $\Q$ and a form $A \in V_d(\Z_v)$, let $\wt{\kappa}_v(A) \defeq \kappa_v(A)$ if $v \neq 2$, $\wt{\kappa}_v(A) = 1$ if $v = 2$ and $A$ is isotropic, and $\wt{\kappa}_v(A) = -1$ if $v = 2$ and $A$ is anisotropic. 
For each $f \in U(\Z_v)$, maximal if $v = p$ is a prime, and each $\varepsilon_1, \varepsilon_2 \in \{0,\pm 1\}$, define the following masses:
\begin{align*}
\wt{\on{Mass}}_\infty^{\varepsilon_1, \varepsilon_2}(f; (i)) & \defeq \sum_{\substack{(A,B) \in \frac{W^\vee(\R)^{(i)} \cap \on{Res}^{-1}(4f)}{\on{SL}_3(\R)} \\ \wt{\kappa}_\infty(A) = \varepsilon_1,\, \wt{\kappa}_\infty(B) = \varepsilon_2 \\ (A,B) \text{ soluble}}} \frac{1}{\#\on{Stab}_{\on{SL}_3(\R)}(A,B)}; \\
\wt{\on{Mass}}_p^{\varepsilon_1, \varepsilon_2}(f) & \defeq \sum_{\substack{(A,B) \in \frac{W^\vee(\Z_p) \cap \on{Res}^{-1}(4f)}{\on{SL}_3(\Z_p)} \\ \wt{\kappa}_p(A) = \varepsilon_1,\, \wt{\kappa}_p(B) = \varepsilon_2 \\ (A,B) \text{ soluble}}} \frac{1}{\#\on{Stab}_{\on{SL}_3(\Z_p)}(A,B)}.
\end{align*}
Then the masses $\on{Mass}_\infty(4f; A; B; (i))$ with $i \in \{(1111),(22+),(112)\}$ and $\on{Mass}_p(4f;A,B;\phi)$ are respectively of the form $\wt{\on{Mass}}_\infty^{\varepsilon_1, \varepsilon_2}(f; (i))$ and $\wt{\on{Mass}}_p^{\varepsilon_1, \varepsilon_2}(f)$ for some $\varepsilon_1, \varepsilon_2 \in \{0,\pm1\}$. Thus, to apply Theorem~\ref{th:sec4main}, it suffices to evaluate these new masses for each choice of $\varepsilon_1, \varepsilon_2 \in \{0, \pm 1\}$, $(i)\in \{(1111),(22+),(112)\}$, and $f \in U(\Z_v)$, maximal if $v = p$ is a prime. We achieve this in the next triple of propositions.

We start by computing the mass at infinity, for which we require the following intermediate result.

\begin{lemma}\label{lem:Rorbits}
Let $f \in U_{1,d}(\R)^+$, and denote the three real roots of $f(x,1)$ by $r_1<r_2<r_3$. Then the following four pairs of matrices are representatives for the action of $\SL_3(\R)$ on pairs $(A,B)\in W^\vee(\R)$ with $\det(xA-B)=f(x)$:
\begin{equation*}
\begin{array}{rcl}
(A_1,B_1)=\left(\begin{pmatrix} -1 & & \\ & \hphantom{-}1 & \\ & & -1 \end{pmatrix},
\begin{pmatrix} -r_1 & & \\ & \hphantom{-}r_2 & \\ & & -r_3 \end{pmatrix}
\right),&\hspace*{-15pt}&
(A_2,B_2)=\left(\begin{pmatrix} \hphantom{}1 & & \\ & -1 & \\ & & -1 \end{pmatrix},
\begin{pmatrix} \hphantom{}r_1 & & \\ & -r_2 & \\ & & -r_3 \end{pmatrix}
\right),
\\[.3in]
(A_3,B_3)=\left(\begin{pmatrix} -1 & & \\ & -1 & \\ & & \hphantom{-}1 \end{pmatrix},
\begin{pmatrix} -r_1 & & \\ & -r_2 & \\ & & \hphantom{-}r_3 \end{pmatrix}
\right),&\hspace*{-15pt}&
(A_4,B_4)=\left(\begin{pmatrix} \hphantom{}1 & & \\ & \hphantom{-}1 & \\ & & \hphantom{-}1 \end{pmatrix},
\begin{pmatrix} \hphantom{}r_1 & & \\ & \hphantom{-}r_2 & \\ & & \hphantom{-}r_3 \end{pmatrix}
\right).
\end{array}
\end{equation*}
Moreover, the splitting types of $(A_1,B_1)$ is $(1111)$, the splitting types of $(A_2,B_2)$ and $(A_3,B_3)$ are $(22\pm)$, and the splitting type of $(A_4,B_4)$ is $(22\#)$.
\end{lemma}
\begin{proof}
From a parametrisation result of Bhargava~\cite[Theorem~4]{MR2081442}, it follows that the set of $\SL_3(\R)$-orbits on such pairs $(A,B)$ are in bijection with classes in $(\R^{\times}/\R^{\times2})^3$ having norm $1 \in \Q^\times/\Q^{\times2}$. There are four such classes, namely the classes of $(1,1,1)$, $(1,-1,-1)$, $(-1,1,-1)$, and $(-1,-1,1)$, and the four pairs $(A_i,B_i)$ correspond respectively to these classes. It is evident that $A_1, B_1$ have common zeros, and so the splitting type of $(A_1, B_1)$ must be $(1111)$. Since $A_4$ is anisotropic, the splitting type of $(A_4,B_4)$ must be $(22\#)$. The other two pairs must then have splitting type $(22\pm)$.
\end{proof}

Lemma~\ref{lem:Rorbits} immediately yields the following values for the infinite mass:
\begin{proposition} \label{lem-rorbits}
Let $f \in U(\R)^+$. Then we have
\begin{equation*}
\wt{\on{Mass}}_\infty^{1,1}\big(f;(1111)\big)=\frac14;\quad \wt{\on{Mass}}_\infty^{-1,\pm1}\big(f;(1111)\big)=\wt{\on{Mass}}_\infty^{\pm1,-1}\big(f;(1111)\big)=0.
\end{equation*}
Let notation be as in Lemma \ref{lem:Rorbits}. Then we have 
\begin{equation*}
\begin{array}{rcl}
\wt{\on{Mass}}_\infty^{-1,\pm1}\big(f;(22+)\big)
&=&0;
\\[.1in]
\wt{\on{Mass}}_\infty^{1,1}\big(f;(22+)\big)&=& \begin{cases}
\frac14, & \text{ if $r_2 > 0$ or $r_3 < 0$}; \\
0, &\text{else};
\end{cases}
\\[.25in]
\wt{\on{Mass}}_\infty^{1,-1}(f;(22+))&=&
\frac14-\wt{\on{Mass}}_\infty^{1,1}\big(f;(22+)\big).
\end{array}
\end{equation*}
Let $f \in U(\R)^-$. Then we have
\begin{equation*}
\wt{\on{Mass}}_\infty^{1,1}\big(f;(112)\big)=\frac12;\quad \wt{\on{Mass}}_\infty^{-1,\pm1}\big(f;(112)\big)=\wt{\on{Mass}}_\infty^{\pm1,-1}\big(f;(112)\big)=0.
\end{equation*}
\end{proposition}
%\begin{proof}
%The claim in the first displayed equation is immediate. 
%Next assume that $A$ is isotropic over $\R$, and let $(A,B) \in \on{Res}^{-1}(f) \cap W(\R)^{(22+)}$. It is clear that we have $\on{Mass}_\infty^{+,\pm}(f;(22+))=1/4$ and $\on{Mass}_\infty^{+,\mp}(f;(22+))=0$ when $\wt{\kappa}_\infty(B)=\pm 1$. To finish the proof, we divide into cases depending on the signs of the roots $r_1 < r_2 < r_3$ of $f(x,1)$ and apply Lemma \ref{lem:Rorbits}. When all the roots $r_i$ have the same sign, both $B_2$ and $B_3$ are isotropic, implying that $\on{Mass}_\infty^{+,+}(f;(22+))=1/4$. 
%
%Next, suppose that $r_3$ is the only positive root. Then $B_2$ is isotropic and $B_3$ is negative-definite. Since $B_3$ takes only negative values (in particular, it takes negative values on the zero set of $A_3$), it follows that the splitting type of $(A_3,B_3)$ is $(22-)$ and hence also that the splitting type of $(A_2,B_2)$ is $(22+)$. Since $B_2$ is isotropic, the claim follows. Finally, when $r_1$ is the only negative root, $B_2$ is positive-definite, hence anisotropic, and the splitting type of $(A_2,B_2)$ is $(22+)$. Once again, the claim follows.
%\end{proof}

We now compute the mass at odd primes $p$:

\begin{proposition}\label{lem:porbits}
    Let $p > 2$ be a prime, let $f \in U_{1,d}(\Z_p)_{\max}$, and let $\varepsilon_1, \varepsilon_2 \in \{0,\pm 1\}$. Then we have $$\wt{\on{Mass}}_p^{\varepsilon_1, \varepsilon_2}(f) = \on{Mass}_p^{\varepsilon_1, \varepsilon_2}(f).$$
\end{proposition}
\begin{proof}
    Since $p$ is odd, every pair $(A,B) \in \on{Res}^{-1}(f) \cap W^\vee(\Z_p)$ is soluble over $\Q_p$ by~\cite[Corollary~8.4]{MR3156850}.
\end{proof}

We finish by computing the mass at $2$. For the sake of simplicity, we restrict our forms $f$ to lie in the set $\Sigma_2$ defined by 
\begin{equation} \label{eq-thisissigma2}
\Sigma_2 \defeq \{f \in U_{1,d}(\Z_2) : f(x,1) \equiv  x^3 + x + 1 \pmod 8\} \subset U_{1,d}(\Z_2)_{\max}.
\end{equation}
Since $x^3 + x + 1$ is irreducible over $\F_2$, it follows from~\cite[\S6.2]{MR3156850} that there exist exactly $2$ $\Q_2$-soluble orbits of $\on{SL}_3(\Q_2)$ on $W^\vee(\Q_2)$ having any given resolvent $f \in \Sigma_2$. In the following pair of lemmas, we identify and characterise these $2$ $\Q_2$-soluble orbits:

\begin{lemma} \label{lem-soldist}
    Let $f \in \Sigma_2$. The two $\Q_2$-soluble orbits of $\on{SL}_3(\Q_2)$ on $W^\vee(\Q_2)$ having resolvent $4f$ are given by the distinguished and $\Delta$-distinguished orbits.
\end{lemma}
\begin{proof}
    Recall from \S\ref{sec-dist} that the binary quartic form $g$ associated to a representative of the distinguished orbit is necessarily reducible. Thus, the equation $z^2 = g(x,y)$ is $\Q_2$-soluble, implying that the distinguished orbit is $\Q_2$-soluble.

    Next, take $b = 8b'$, $c = 8c' + 1$, and $d = 8d' + 1$ for $b',c',d' \in \Z_2$. Then, from~\eqref{eq-binquartdel}, the binary quartic form $g$ associated to a representative of the $\Delta$-distinguished orbit is given by
    \begin{equation} \label{eq-binquart2}
    g(x,y) = d_k x^4 - \frac{b}{2} x^2y^2 + d_m xy^3 + \left(\frac{b^2}{16d_k} - \frac{c}{4d_k}\right)y^4
    \end{equation}
    It is evident that the equation $z^2 = g(x,y)$ is $\Q_2$-soluble: indeed, $(x,y) = (1,0)$ is a solution because $d_k \equiv 1 \pmod 8$ is a quadratic residue. Thus, the $\Delta$-distinguished orbit is $\Q_2$-soluble. It now remains to show that the distinguished and $\Delta$-distinguished orbits are distinct, for which it suffices to show that $g(x,y)$ has no primitive roots modulo $8$ (implying that $g(x,y)$ has no linear factors over $\Q_2$). Consider the binary quartic form $d_kg(x,2y)$, which has coefficients in $\Z_2$. Modulo $8$, we have $d_kg(x,2y) \equiv x^4-4y^4$, which has no primitive roots, as desired.
\end{proof}

\begin{lemma} \label{lem-soliso}
    Let $(A,B) \in W_{1,d}(\Z_2)$ be soluble with resolvent $f \in \Sigma_2$. Then $B$ is isotropic; i.e., $\wt{\kappa}_2(B) = 1$.
\end{lemma}
\begin{proof}
    This holds by definition for both the distinguished and $\Delta$-distinguished orbits. The lemma then follows from Lemma~\ref{lem-soldist}.
\end{proof}

Now, combining Lemmas~\ref{lem-soldist} and~\ref{lem-soliso} yields that the $2$-adic mass is as follows: 

\begin{proposition} \label{lem-2orbits}
    Let $f \in \Sigma_2$. Then we have 
    $$\wt{\on{Mass}}_2^{1,1}(f) = 2; \quad \wt{\on{Mass}}_2^{-1,\pm1}(f)=\wt{\on{Mass}}_2^{\pm1,-1}(f)=0.$$
\end{proposition}

\subsection{Proofs of Theorems~\ref{third main} and \ref{third main2}} \label{sec-selthatsit}

We are now ready to deduce Theorems \ref{third main} and \ref{third main2}. As in the theorem statements, fix $\varepsilon \in \{\pm\}$ and $\Sigma = U_{1,d}(\Z)_{\max}^\varepsilon \cap \bigcap_p \Sigma_p$ to be an acceptable family with $\Sigma_2$ as in~\eqref{eq-thisissigma2}.  By considering $W^\vee(\Z)$ as a ($\SL_3(\Z)$-invariant) subset of 
$W(\Z)$, we may apply Theorem \ref{th:sec4main} to count $\on{SL}_3(\Z)$-orbits on $W^\vee(\Z)$. 

Fix $\circ \in \{\on{bal},\on{wei}\}$. 
Suppose $\varepsilon  = +$, and for each prime $p$, let $\phi_p$ be the characteristic function of the set $\{(A,B) \in \on{Res}^{-1}(4\Sigma_p) \cap W_{1,d}^{\vee}(\Z_p) : (A,B) \text{ is $\Q_p$-soluble}\}$. 
   Given $k \in \Z \smallsetminus \{0\}$, let $\wt{\mc{G}}_k$ denote the set of genera of integer-matrix quadratic forms $V_{4k}(\Z)$. 
 Now, combining Theorems~\ref{thm-2selparam},~\ref{thm-howmanydist},~\ref{th:sec4main}, and~\ref{th:cubiccount}, along with Proposition~\ref{lem-rorbits}, yields that the contribution to the average size of $\on{Sel}_2(E_f)$ over forms $f \in \Sigma$ from those $2$-Selmer elements that correspond to pairs $(A,B) \in W(\R)^{(0)}$ that are not $\Delta$-distinguished over $\Z$ is bounded above by
\begin{align}
  & 1 +  \sum_{\substack{B \in \wt{\mc{G}}_{-d} \\ (A,B) \in W(\R)^{(0)}}} \prod_p\frac{1}{ \on{Vol}(\Sigma_p)}\int_{f \in \Sigma_p} \wt{\on{Mass}}_p^{\wt{\kappa}_p(A),\wt{\kappa}_p(B)}(f) \mathrm{d}f. \label{eq-sepitout}
  \end{align}
As stated in Theorem~\ref{thm-2selparam}, any pair $(A,B)$ corresponding to a $2$-Selmer element can be translated under the action of $\on{SL}_3(\Z)$ to be in the form $(\cA, B)$, where $\cA$ is the anti-diagonal matrix with anti-diagonal entries given by $1$, $-1$, and $1$. Thus, we may take $A = \cA$ and $\wt{\kappa}_p(A) = 1$ for all $p$.
  
Instead of Lemma~\ref{lem-kappagenerate}, we use the following analogue for integer-matrix quadratic forms, which follows from~\cite[Chapter 15, \S7.7]{MR1662447}:
    \begin{lemma} \label{lem-kappagenerate2}
Let $m \equiv 1 \pmod 8$. The map sending $M \in \wt{\mc{G}}_m$ to the tuple $(\wt{\kappa}_v(M))_v \in \prod_v \{0,\pm1\}$ defines a bijection from the set $\wt{\mc{G}}_m$ to the set of tuples $(\varepsilon_v)_v \in \prod_v \{0,\pm1\}$ over places $v$ of $\Q$, where $\varepsilon_p = 0$ if and only if $p \nmid m$ and $\varepsilon_\infty \times \prod_{p \in T_m^{\on{odd}}} \varepsilon_p = 1$.
\end{lemma}
\noindent Applying Lemma~\ref{lem-kappagenerate2}, Propositions~\ref{lem:porbits} and~\ref{lem-2orbits}, and the formula~\eqref{eq-massinrhos} to~\eqref{eq-sepitout}, we find that the aforementioned contribution is bounded above by
  \begin{align}
& 1 +  \sum_{\substack{(\varepsilon_{2,p})_p \in \prod_p \{0,\pm1\} \\ \prod_{p \in T_d^{\on{odd}}} \varepsilon_{2,p}  = \prod_{p \in T_d^{\on{even}}} |\varepsilon_{2,p}| = 1}} \prod_p\frac{1}{ \on{Vol}(\Sigma_p)}\int_{f \in \Sigma_p} \wt{\on{Mass}}_p^{1,\varepsilon_{2,p}}(f)\mathrm{d}f = \nonumber \\
    &\qquad\qquad\qquad\qquad\qquad\qquad\qquad 1 +  2\times \sum_{\substack{(\varepsilon_{2,p})_p \in \{\pm 1\}^{T_d^{\on{odd}}} \\ \prod_{p \in T_d^{\on{odd}}} \varepsilon_{2,p} = 1}}  \prod_{\substack{p \in T_d^{\on{odd}}}} \bigg(\frac{1}{2} + \varepsilon_{2,p}\frac{\rholeft_\Sigma(p)}{2}\bigg) = 2 + \rholeft_\Sigma. \label{eq-0contr}
\end{align}
Now, with notation as in Lemma~\ref{lem:Rorbits}, let 
\begin{equation} \label{eq-defpid}
\Pi_{d,\on{H}_{\circ}} \defeq \limsup_{X \to \infty} \frac{\on{Vol}(\{f \in U_{1,d}(\R)_X^{\circ,+} : r_2 > 0 \text{ or } r_3 < 0\})}{\on{Vol}(U_{1,d}(\R)_X^{\circ,+})}.
\end{equation}
Note that $\Pi_{d,\on{H}_{\circ}} = 1$ if $d > 0$.  The contribution to the average size of $\on{Sel}_2(E_f)$ over forms $f \in \Sigma$ from those $2$-Selmer elements that correspond to pairs $(A,B) \in W(\R)^{(2+)}$ that are not $\Delta$-distinguished over $\Z$ is bounded above by
\begin{align}
  & 
   \Pi_{d,\on{H}_{\circ}} \times \sum_{\substack{(\varepsilon_{2,p})_p \in \prod_p \{0,\pm1\} \\ \prod_{p \in T_d^{\on{odd}}} \varepsilon_{2,p}  = \prod_{p \in T_d^{\on{even}}} |\varepsilon_{2,p}| = 1}} \prod_p\frac{1}{ \on{Vol}(\Sigma_p)}\int_{f \in \Sigma_p} \wt{\on{Mass}}_p^{\varepsilon_{1,p},\varepsilon_{2,p}}(f)\mathrm{d}f + \nonumber \\
&\qquad\qquad\qquad\qquad  (1-\Pi_{d,\on{H}_{\circ}}) \times \sum_{\substack{(\varepsilon_{2,p})_p \in \prod_p \{0,\pm1\} \\ -\prod_{p \in T_d^{\on{odd}}} \varepsilon_{2,p}  = \prod_{p \in T_d^{\on{even}}} |\varepsilon_{2,p}| = 1}} \prod_p\frac{1}{ \on{Vol}(\Sigma_p)}\int_{f \in \Sigma_p} \wt{\on{Mass}}_p^{\varepsilon_{1,p},\varepsilon_{2,p}}(f)\mathrm{d}f   = \nonumber \\
    &   2\Pi_{d,\on{H}_{\circ}}\times \sum_{\substack{(\varepsilon_{2,p})_p \in \{\pm 1\}^{T_d^{\on{odd}}} \\ \prod_{p \in T_d^{\on{odd}}} \varepsilon_{2,p} = 1}}  \prod_{\substack{p \in T_d^{\on{odd}}}} \bigg(\frac{1}{2} + \varepsilon_{2,p}\frac{\rholeft_\Sigma(p)}{2}\bigg) +\nonumber \\
    &\qquad\qquad\qquad\qquad  2(1-\Pi_{d,\on{H}_{\circ}})\times \sum_{\substack{(\varepsilon_{2,p})_p \in \{\pm 1\}^{T_d^{\on{odd}}} \\ \prod_{p \in T_d^{\on{odd}}} \varepsilon_{2,p} = -1}}  \prod_{\substack{p \in T_d^{\on{odd}}}} \bigg(\frac{1}{2} + \varepsilon_{2,p}\frac{\rholeft_\Sigma(p)}{2}\bigg) = \nonumber \\
&
\Pi_{d,\rmH_\circ}(1+\rholeft_\Sigma)+(1-\Pi_{d,\rmH_\circ})(1 - \rholeft_\Sigma)= 1 + \big(2\Pi_{d,\on{H}_{\circ}} -1\big) \rholeft_\Sigma.
 \label{eq-2contr}
\end{align}
Now, let $\mathfrak{d}_{\Sigma}'$ denote the density of $f \in \Sigma$ above which the $\Delta$-distinguished orbit contains a locally soluble integral representative. Then it follows from Theorem~\ref{th:deltadistex}, Corollary~\ref{cor-rsol}, and Lemma~\ref{lem-soldist}, along with the fact that when $d < 0$, the condition $r_2 > 0$ or $r_3 < 0$ is equivalent to $b < 0$ and $c > 0$, that 
\begin{equation} \label{eq-deltprim}
\mathfrak{d}_\Sigma' = \Pi_{d,\on{H}_{\circ}} \rholeft_\Sigma.
\end{equation}
Note in particular that $\Pi_{d,\rmH_{\circ}}$ depends only on the sign of $d$ and the choice of height, and we have
\begin{equation} \label{eq-newpid}
    \Pi_{d,\rmH_{\circ}} = \begin{cases} 1, & \text{if $d > 0$,} \\ 1/4,& \text{if $\circ = \on{bal}$ and $d < 0$,} \\ 1/26,& \text{if $\circ = \on{wei}$ and $d < 0$.} \end{cases}
\end{equation}
Combining~\eqref{eq-0contr},~\eqref{eq-2contr}, and~\eqref{eq-deltprim}, we conclude that
\begin{equation*}
    \limsup_{X \to \infty} \underset{f \in \Sigma_X^\circ}{\on{Avg}}\,\on{Sel}_2(E_f)  \leq 3 + 3\Pi_{d,\on{H}_{\circ}}\rholeft_\Sigma 
\end{equation*}
which is precisely as claimed in Theorem~\ref{third main}.

\medskip

Next, suppose $\varepsilon = -$, and let $\phi$ be as above. Then combining Theorems~\ref{thm-2selparam},~\ref{thm-howmanydist},~\ref{th:sec4main}, and~\ref{th:cubiccount}, along with Lemma~\ref{lem-kappagenerate2}, Propositions~\ref{lem:porbits},~\ref{lem-rorbits}, and~\ref{lem-2orbits}, and the formula~\eqref{eq-massinrhos}, yields that the average size of $\on{Sel}_2(E_f)$ over forms $f \in \Sigma$ is bounded as follows:
\begin{align}
    \limsup_{X \to \infty} \underset{f \in \Sigma_X^\circ}{\on{Avg}}\,\on{Sel}_2(E_f)   & 
    \leq 1 + \mathfrak{d}_\Sigma' + 2 \times\sum_{\substack{B \in \wt{\mc{G}}_{-d} \\ (\cA,B) \in W(\R)^{(1)}}} \prod_p\frac{1}{ \on{Vol}(\Sigma_p)}\int_{f \in \Sigma_p} \wt{\on{Mass}}_p^{\wt{\kappa}_p(\cA),\wt{\kappa}_p(B)}(f) \mathrm{d}f\nonumber \\
    & = 1 + \mathfrak{d}_\Sigma' + 2\times \sum_{\substack{(\varepsilon_{2,p})_p \in \prod_p \{0,\pm1\} \\ \prod_{p \in T_d^{\on{odd}}} \varepsilon_{2,p}  = \prod_{p \in T_d^{\on{even}}} |\varepsilon_{2,p}| = 1}} \prod_p\frac{1}{ \on{Vol}(\Sigma_p)}\int_{f \in \Sigma_p} \wt{\on{Mass}}_p^{1,\varepsilon_{2,p}}(f)\mathrm{d}f \nonumber \\
    & = 1 +\mathfrak{d}_\Sigma' + 4 \times \sum_{\substack{(\varepsilon_{2,p})_p \in \{\pm 1\}^{T_d^{\on{odd}} } \\ \prod_{p \in T_d^{\on{odd}}} \varepsilon_{2,p} = 1}} \prod_{\substack{p \in T_d^{\on{odd}}}} \bigg(\frac{1}{2} + \varepsilon_{2,p}\frac{\rholeft_\Sigma(p)}{2}\bigg) \nonumber \\
   & = 3 + 3\rholeft_\Sigma, \nonumber
\end{align}
which is precisely as claimed in Theorem~\ref{third main2}.

\appendix

\section{Computations of local densities} \label{appendix-localcomputationsat2}

This section is devoted to computing the local densities required to deduce Theorems~\ref{first main} and~\ref{first main2} from Theorems~\ref{second main real} and~\ref{second main complex}. In \S\ref{sec-maxsuff}, we determine the local density of maximal binary cubic forms, as well as the local density of sufficiently-ramified forms among the maximal ones. In \S\ref{sec-distdensities}, we compute the density of forms in an acceptable family above which there is an \mbox{integral $\Delta$-distinguished orbit.}

\subsection{Densities of maximal and sufficiently ramified cubic forms} \label{sec-maxsuff}

A binary cubic form $f \in U(\Z)$ is maximal if and only if it is maximal over $\Z_p$ for each prime $p$. The following proposition gives the probability that $f \in U_{a,d}(\Z_p)$ is maximal, as well as the probability that $f \in U_{a,d}(\Z_p)_{\max}$ is left- or right-sufficiently-ramified at $p$.

\begin{proposition} \label{prop-maxcalcs}
    The $p$-adic density of $U_{a,d}(\Z_p)_{\max}$ in $U_{a,d}(\Z_p)$ is given by
\begin{equation} \label{eq-maxdens}
\begin{array}{ll}
1 - 3(p^{-2} - p^{-3})\chi_p(d), & \text{if $p \not\in T_a$ and $p \not\in T_d$ and $p \equiv 1\,\,\,(\on{mod} 3)$,} \\[0.1cm] 
1 - p^{-2} + p^{-3}, & \text{if $p \not\in T_a$ and $p \not\in T_d$ and $p \not\equiv 1\,\,\,(\on{mod} 3)$,} \\[0.1cm]
1 - p^{-2} + p^{-3}, & \text{if $p \mid\mid a$ and $p \not\in T_d$, or $p \not\in T_a$ and $p \mid\mid d$,}\\[0.1cm]
(1-p^{-1})(1-p^{-2}), & \text{if $\nu_p(a) \geq 2$ and $p \not\in T_d$, or $p \not\in T_a$ and $\nu_p(d) \geq 2$,} \\[0.1cm]
1 - p^{-2}, & \text{if $p \mid\mid a$ and $p \mid\mid d$,} \\[0.1cm]
1 - p^{-1}, & \text{if $p \mid\mid a$ and $\nu_p(d) \geq 2$, or $\nu_p(a) \geq 2$ and $p \mid\mid d$,} \\[0.1cm]
(1 - p^{-1})^2, & \text{if $\nu_p(a) \geq 2$ and $\nu_p(d) \geq 2$,}
\end{array}
\end{equation}
    where $\chi_p(d) = 1$ if $d \in (\Z/p\Z)^{\times3}$ and $\chi_p(d) = 0$ otherwise. Moreover, the $p$-adic density of left-sufficiently-ramified elements in $U_{a,d}(\Z_p)_{\max}$ is given by
    \begin{equation} \label{eq-sqdens}
\begin{array}{ll}
(p^2 - p+1)(p^3 - p + 1)^{-1}, & \text{if $p \mid\mid a$ and $p \not\in T_d$,}\\[0.1cm]
(p-1)^2(p^3 - p + 1)^{-1}, & \text{if $\nu_p(a) \geq 2$ and $p \not\in T_d$,} \\[0.1cm]
(p+1)^{-1}, & \text{if $p \mid\mid a$ and $p \mid\mid d$,} \\[0.1cm]
p^{-1}, & \text{if $\nu_p(a) \geq 2$ and $p \mid\mid d$,} \\[0.1cm]
0, & \text{else.}
\end{array}
\end{equation}
The right-sufficient-ramification densities may be obtained from~\eqref{eq-sqdens} by switching $a$ and $d$.
\end{proposition}
\begin{proof}
A binary cubic form $f \in U(\Z/p^2\Z)$ is said to satisfy \emph{Dedekind's criterion} if neither of the following two conditions holds:
\begin{itemize}
\item There exists $\gamma \in \on{SL}_2(\Z_p)$ such that $\gamma \cdot f$ has $x^3$-coefficient divisible by $p^2$ and $x^2y$-coefficient divisible by $p$; or
\item Every coefficient of $f$ is divisible by $p$ $($i.e., $f$ is imprimitive$)$.
\end{itemize}
It is well-known (see, e.g.,~\cite[Lemma~4.33]{BHSpreprint}) that a binary cubic form $f \in U(\Z_p)$ is maximal if and only if it satisfies Dedekind's criterion. Thus, it suffices to count, for each $a,d \in \Z$, the number of binary cubic forms $f\in U_{a,d}(\Z/p^2\Z)$ having each possible splitting type\footnote{The \emph{splitting type} of a binary cubic form $f \in U_{a,d}(\F_p)$ is the tuple $\big(e_1^{f_1}\cdots e_k^{f_k}\big)$, where the $e_i$ are the degrees of the irreducible factors of $f$, and the $f_i$ are their corresponding multiplicities. This tuple is unordered unless it is $(11^2)$ or $(1^21)$, in which case the $1^2$ appears on the left if and only if the double root of $f$ is given by $(x,y) = (1,0)$.} that satisfy Dedekind's criterion. We divide into cases depending on whether $p \mid a$ and/or $p \mid d$:

\smallskip
\noindent \emph{Case 1}: $p \nmid a$ and $p \nmid d$. It is easy to verify that the following table gives the number of forms in $U_{a,d}(\Z/p^2\Z)$ having each possible splitting type that satisfy Dedekind's criterion:

\begin{center}
\begin{tabular}{ |c|c|c| } 
 \hline
 splitting type & $\#$ of forms \\ \hline 
 $(111)$ & $\begin{cases} p^2\big((p-1)(p-4)/6 + \chi_p(d)\big), & \text{if $p \equiv 1\pmod 3$} \\ p^2(p-2)(p-3)/6, & \text{if $p \not\equiv 1 \pmod 3$} \end{cases}$ \\
 $(12)$ &  $p^3(p-1)/2$ \\
  $(3$) & $\begin{cases} p^2\big((p^2+p+1)/3 - \chi_p(d)\big), & \text{if $p \equiv 1 \pmod3$}\\ p^3(p+1)/3, & \text{if $p \not\equiv 1 \pmod3$} \end{cases}$ \\
 $(11^2)$ & $\begin{cases} p(p-1)\big(p-1-3\chi_p(d)\big), & \text{if $p \equiv 1 \pmod 3$} \\ p(p-1)(p-2), & \text{if $p \not\equiv 1 \pmod 3$} \end{cases}$ \\
 $(1^3)$ & $p(p-1)$ \\
 \hline
\end{tabular}
\end{center}

  \smallskip
\noindent \emph{Case 2}: $p \mid a$ and $p \nmid d$. It is easy to verify that the following table gives the number of forms in $U_{a,d}(\Z/p^2\Z)$ having each possible splitting type that satisfy Dedekind's criterion:

\begin{center}
\begin{tabular}{ |c|c|c| } 
 \hline
 splitting type & $\#$ of forms \\ \hline 
 $(111)$ & $p^2(p-1)(p-2)/2$ \\
 $(12)$ &  $p^3(p-1)/2$ \\
 $(1^21)$ & $\begin{cases} p^2(p-1), & \text{if $\nu_p(a) = 1$,} \\ 0, & \text{else} \end{cases}$ \\
 $(11^2)$ & $p(p-1)^2$ \\
 $(1^3)$ & $\begin{cases} p^2, & \text{if $\nu_p(a) = 1$,} \\ 0, & \text{else} \end{cases}$ \\
 \hline
\end{tabular}
\end{center}

  \smallskip
\noindent \emph{Case 3}: $p \nmid a$ and $p \mid d$. On account of the transformation $f(x,y) \mapsto f(y,x)$, the results for this case are the same as those for Case 2.

\smallskip
\noindent \emph{Case 4}: $p \mid a$ and $p \mid d$. It is easy to verify that the following table gives the number of forms in $U_{a,d}(\Z/p^2\Z)$ satisfying Dedekind's criterion and having each possible splitting type:

\begin{center}
\begin{tabular}{ |c|c|c| } 
 \hline
 splitting type & $\#$ of forms \\ \hline 
 $(111)$ & $p^2(p-1)^2$ \\
 $(1^21)$ & $\begin{cases} p^2(p-1), & \text{if $\nu_p(a) = 1$,} \\ 0, & \text{else} \end{cases}$ \\
 $(11^2)$ & $\begin{cases} p^2(p-1), & \text{if $\nu_p(d) = 1$,} \\ 0, & \text{else} \end{cases}$ \\
 \hline
\end{tabular}
\end{center}
The densities claimed in~\eqref{eq-maxdens} and~\eqref{eq-sqdens} follow immediately from the three tables presented above. This completes the proof of Proposition~\ref{prop-maxcalcs}.
\end{proof}

\subsection{The density of $\Delta$-distinguished orbits} \label{sec-distdensities}

We start by giving the proof of Theorem~\ref{th:deltadistex}, which characterises in terms of congruence conditions when the $\Delta$-distinguished orbit exists over a maximal binary cubic form. The proof of Theorem~\ref{th:deltadistex2} is completely analogous, so we omit it.

\begin{proof}[Proof of Theorem~\ref{th:deltadistex}]
The uniqueness claim follows from the general fact that every $\on{SL}_3(\Q)$-equivalence class on $W(\Z)$ with maximal resolvent contains exactly one $\on{SL}_3(\Z)$-orbit.

 Let $(A,B)$ be a representative in a $\Delta$-distinguished orbit with $\Res(A,B)=f$ satisfying the conditions $a_{11}a_{22} - a_{12}^2 = b_{11}b_{22} - b_{12}^2 = 0$. We claim that, using the action of $\on{SL}_3(\Z)$, we can replace $(A,B)$ with a pair of the following form:
    %By Definition~\ref{def-thatsdelt}, there exists a $\P^1_K\subset \P^2K$ such that the restriction of $A$ and $B$ to this $\P^1_K$ give reducible binary quadratic forms. (Note that neither of these binary forms can be $0$, otherwise we would have $\Delta(A,B)=0$.) 
    %By replacing $(A,B)$ with an $\SL_3(\Z)$-translate, we may assume that this $\P^1$ is the subspace $(*,*,0)$,
    %Using the action of the block-diagonal subgroup $\{g \in \on{GL}_2(\Z) \times \on{GL}_1(\Z) : \det g = 1\} \subset \on{SL}_3(\Z)$, we can further arrange that $a_{12}=a_{22}=b_{11}=b_{12}=0$ (this corresponds to moving the points of bitangency to $0$ and $\infty$ on the bitangent line). %Since $a_{11}$ and $b_{22}$ are nonzero, we may use another $\SL_3(K)$-transformation to further assume that $a_{13}=b_{23}=0$. 
    %Then we have $a/4=\det(A)=-a_{11}a_{23}^2$ and $b/4=\det(B)=-b_{22}b_{13}^2$. Thus, with another $\SL_3(K)$ transformation, we can ensure that $a_{11}=-a$, $b_{22}=-d$, and $a_{23}=b_{13}=1/2$. Finally, the coefficients $a_{33}$ and $b_{33}$ are uniquely determined by the values of $b$ and $c$. The lemma follows.
%Consequently, we can assume $(A,B)$ to be of the following form:
\begin{equation}  \label{eq-secondkostant}
(A,B) = \left(\begin{pmatrix}
-a_{11} & 0 & a_{13}/2\\
0 & 0 & a_{23}/2\\
a_{13}/2 & a_{23}/2 & a_{33}
\end{pmatrix},
\begin{pmatrix}
0 & 0 & b_{13}/2 \\
0 & b_{22} & b_{23}/2 \\
b_{13}/2 & b_{23}/2 & b_{33}
\end{pmatrix}\right)
\end{equation}
Indeed, using the action of the block-diagonal subgroup $\{g \in \on{GL}_2(\Z) \times \on{GL}_1(\Z) : \det g = 1\} \subset \on{SL}_3(\Z)$, we can certainly arrange that $B$ takes the shape given in~\eqref{eq-secondkostant}. Then, if $a_{11} \mid a_{12}$, we can use the action of the lower-triangular unipotent subgroup of $\on{SL}_3(\Z)$ to clear out $a_{12}$. Otherwise, there must exist a prime $p$ such that $p^2 \mid a_{11}$ and $p \mid a_{12}$. If we also have $p \mid a_{13}$, then taking the resolvent of $f = \on{Res}(A,B)$, we see that $\wt{f}$ fails Dedekind's criterion, as $p^2$ divides the $x^3$-coefficient and $p$ divides the $x^2y$-coefficient. We may thus assume $p \nmid a_{13}$ and consider the $\on{GL}_2(\Z)$-equivalent binary cubic form $f(a_{13}^{-1}x+a_{13}^{-1}b_{13}y,y)$, which also fails Dedekind's criterion, this time by examining the coefficients of $xy^2$ and $y^3$. Thus, we must have that $a_{11} \mid a_{12}$, and the claim follows. 
%We now prove existence. One can check using the action of the block-diagonal subgroup $\on{SL}_2 \times \on{SL}_1 \subset \on{SL}_3$ that every $\Delta$-distinguished orbit on $W(\Z)$ admits a representative of the following form:

Now, the resolvent of the pair~\eqref{eq-secondkostant} is given explicitly by
\begin{align} \label{eq-secondkostantresolvent}
& (a_{11}a_{23}^2)x^3 + (a_{13}^2b_{22} + 4a_{11}a_{33}b_{22} - 2a_{11}a_{23}b_{23})x^2y + \\
& \qquad\qquad\qquad\,\,\, (a_{11}b_{23}^2 - 4a_{11}b_{22}b_{33} -2a_{13}b_{13}b_{22}  )xy^2 + (b_{13}^2b_{22})y^3. \nonumber
\end{align}
Setting the form in~\eqref{eq-secondkostantresolvent} equal to $f$, we see that $\gcd(a_{11}, b_{22}) = 1$, for otherwise $f$ would fail to be primitive and thus fail Dedekind's criterion.

Suppose $p \mid b_{22}$ for some prime $p$. If $p \mid a_{23}$, then $p^2$ divides the coefficient of $x^3$ and $p$ divides the coefficient of $x^2y$, so by Dedekind's criterion, we must have $p \nmid a_{23}$. Thus, we have $\gcd(b_{22}, a_{23}) = 1$, and by symmetry, we also have $\gcd(a_{11},b_{13}) = 1$. Now suppose $p^2 \mid b_{22}$. Then the $\on{GL}_2(\Z_p)$-equivalent binary cubic form $f(a_{23}^{-1}x + a_{23}^{-1}b_{23}y, y)$ fails Dedekind's criterion, as the coefficients of both $xy^2$ and $y^3$ are divisible by $p^2$. Thus, $b_{22}$ is squarefree, and by symmetry, so is $a_{11}$. It follows that $a_{11} = a_k$, $|a_{23}| = a_m$, $b_{22} = d_k$, and $|b_{13}| = d_m$, and also that $\gcd(a_k,d) = \gcd(a, d_k) = 1$.

We have thus established that the $\Delta$-distinguished orbit with resolvent $f$, if it exists, contains a representative of the canonical form~\eqref{eq-secondkostant}, except possibly for the signs of $a_{23}$ and $b_{13}$. To address this, suppose that $a_{23} = -a_m$. Then by acting on $(A,B)$ with the diagonal matrix with diagonal entries $-1$, $1$, and $-1$, we can flip the sign of $a_{23}$ while preserving the sign of $b_{13}$. Similarly, suppose that $b_{13} = -d_m$. Then by acting on $(A,B)$ with the diagonal matrix with diagonal entries $1$, $-1$, and $-1$, we can flip the sign of $b_{13}$ while preserving the sign of $a_{23}$.

Now, suppose $a_{11}$, $|a_{23}|$, $b_{22}$, and $|b_{13}|$ are determined. We have the following identities:
\begin{align}
 a_{33} & = \frac{-a_{13}^2b_{22}+2a_{11}a_{23}b_{23}+b}{4a_{11}b_{22}}, \label{eq-id1}\\
 b_{33} & = \frac{a_{11}b_{23}^2-2a_{13}b_{13}b_{22}-c}{4a_{11}b_{22}}. \label{eq-id2}
\end{align}
We then split into cases according as $a_{11}$ and $b_{22}$ are odd or even:

\smallskip
\noindent \textbf{Case I: $2 \nmid a_{11}b_{22}$}:
In this case, the identities~\eqref{eq-id1} and~\eqref{eq-id2} may be interpreted as imposing three conditions on the numerators of their right-hand sides: (1) they must be divisible by $a_{11}$; (2) they must be divisible by $b_{22}$; and (3) they must be divisible by $4$. We consider each of these three conditions sequentially as follows:
\begin{enumerate}
    \item[(1)] For $a_{33}$ to be an integer, it follows from~\eqref{eq-id1} that $a_{13} \equiv \pm \sqrt{b_{22}^{-1}b} \pmod{a_{11}}$, which is possible if $c^2 - 4bd \equiv 0 \pmod{a_{11}}$. Then, for $b_{33}$ to be an integer, it follows from~\eqref{eq-id2} that $a_{11}b_{23}^2 - 2a_{13}b_{13}b_{22} - c \equiv \pm 2\sqrt{bd} - c \equiv 0 \pmod{a_{11}}$, which is equivalent to stipulating that $c^2 - 4bd \equiv 0 \pmod{a_{11}}$. Thus, there exists a choice of $a_{13}$ such that the numerators of the right-hand sides of~\eqref{eq-id1} and~\eqref{eq-id2} are both divisible by $a_{11}$ if and only if we have \mbox{$c^2 - 4bd \equiv 0 \pmod{a_{11}}$.}
\item[(2)] By symmetry with the analysis for condition (1), there exists a choice of $b_{23}$ such that the numerators of the right-hand sides of~\eqref{eq-id1} and~\eqref{eq-id2} are both divisible by $b_{22}$ if and only if \mbox{$b^2 - 4ac \equiv 0 \pmod{b_{22}}$.}
\item[(3)] For $a_{33}$ and $b_{33}$ to be integers, we must have 
\begin{equation} \label{eq-4modcase1}
-a_{13}^2b_{22} + 2a_{11}a_{23}b_{23}+b \equiv 0 \pmod{4}\quad \text{and} \quad a_{11}b_{23}^2 - 2a_{13}b_{13}b_{22} - c \equiv 0 \pmod 4.
\end{equation}
It thus suffices to determine for each choice of $(a,b,c,d)$ whether there exists a choice of $(a_{13}, b_{23})$ such that the congruences~\eqref{eq-4modcase1} are satisfied. This can be done by inspection, and the results are given by the first table in the theorem statement.
\end{enumerate}

\medskip
\noindent \textbf{Case II: $2 \nmid a_{11}$ but $2 \mid b_{22}$}: 
In this case, the identities~\eqref{eq-id1} and~\eqref{eq-id2} may be interpreted as imposing three conditions on the numerators of their right-hand sides: (1) they must be divisible by $a_{11}$; (2) they must be divisible by $b_{22}$; and (3) they must be divisible by $8$. The analyses for conditions (1) and (2) are exactly the same as in Case I, so we omit them. As for condition (3), for $a_{33}$ and $b_{33}$ to be integers, we must have
    \begin{equation} \label{eq-4modcase2}
-a_{13}^2b_{22} + 2a_{11}a_{23}b_{23}+b \equiv 0 \pmod{8}\quad \text{and} \quad a_{11}b_{23}^2 - 2a_{13}b_{13}b_{22} - c \equiv 0 \pmod 8.
\end{equation}
It thus suffices to determine for each choice of $(a,b,c,d)$ whether there exists a choice of $(a_{13}, b_{23})$ such that the congruences~\eqref{eq-4modcase2} are satisfied. This can be done by inspection, and the results are given by the second table in the theorem statement.

\medskip
\noindent \textbf{Case III: $2 \mid a_{11}$ but $2 \nmid b_{22}$}: The argument here is entirely analogous to that of Case II.
\end{proof}

Combining Proposition~\ref{prop-maxcalcs} with Theorem~\ref{th:deltadistex} yields the following result, which gives a formula for $\mathfrak{d}_\Sigma$ in the case where $\Sigma_2 = U_{a,d}(\Z_2)_{\max}$.

\begin{corollary} \label{cor-deltadist}
    Let $\Sigma \subset U_{a,d}(\Z)_{\max}$ be an acceptable family with $\Sigma_2 = U_{a,d}(\Z_2)_{\max}$. Then we have $$\mathfrak{d}_\Sigma = \chi_{a,d} \, \mathfrak{d}_2 \times \prod_{\substack{p \mid a_k \\ p > 2}} \rhoright_{\Sigma}(p) \prod_{\substack{p \mid d_k \\ p > 2}} \rholeft_{\Sigma}(p),$$
    where $\mathfrak{d}_2$ is given by
     \begin{equation*} 
\mathfrak{d}_2 = \begin{cases}
3/14, & \text{if $2 \nmid ad$,}\\[0.1cm]
1/6, & \text{if $2 \nmid a_kd_k$ and $2 \mid a_m$ or $d_m$ but not both,} \\[0.1cm]
1/4, & \text{if $2 \nmid a_kd_k$ and $2 \mid \gcd(a_m,d_m)$,} \\[0.1cm]
3/28, & \text{if $2 \nmid ad_m$ and $2 \mid d_k$, or $2 \nmid a_md$ and $2 \mid a_k$,} \\[0.1cm]
1/16, & \text{if $2 \nmid a$ and $8 \mid d$, or $8 \mid a$ and $2 \nmid d$,} \\[0.1cm]
0, & \text{else.}
\end{cases}
\end{equation*}
\end{corollary}

\newpage

\section{Data} \label{section:appendixdata}

 We include two tables of numerically observed average $2$-torsion in class groups of randomly generated fields belonging to each of the possible stable families for each $(a,d) \in \{1,2,3,4,5\}^2$. To be precise, for each such $(a,d)$ and each stable family of $(a,d)$-monogenised cubic fields, we generated a ``bucket'' containing 2000 random members $K$ of the stable family and computed the average value of $\#\Cl(K)[2]$ over the bucket. 

The Monte Carlo algorithm to fill the buckets was written in Python, and the one to compute the average of $\#\Cl(K)[2]$ was written in Magma. The computations were run in parallel on 20 Intel(R) Xeon(R) Gold 6136 CPU @ 3.00GHz cores of Princeton's Program in Applied and Computational Mathematics (PACM) server and took roughly 12 hours. 
The code used to produce these averages, as well as the raw data of the cubic fields, may be found on \href{https://github.com/ashvin-swaminathan/special-cubics}{GitHub}.

For the sake of brevity, we denote the stable families by $\pm B_{\delta_{\on{fin}},\delta_\infty}^{\lambda}$ for $F^{\pm}_{(a,d)}(\lambda;\delta_{\on{fin}},\delta_\infty)$. This matches the names of the files on \href{https://github.com/ashvin-swaminathan/special-cubics}{GitHub}.

\begin{table}[h]
\centering
\begin{tabularx}{\textwidth}{|c||*{12}{X|}}\hline
Family/Bucket & $+B_{0,0}^0$ & $+B_{0,1}^0$ & $-B_{0,1}^0$ & $+B_{0,0}^1$ & $+B_{0,1}^1$ & $-B_{0,1}^1$ & $+B_{0,0}^2$ & $+B_{0,1}^2$ & $+B_{1,0}^2$ & $-B_{0,1}^2$ & $+B_{1,1}^2$ & $-B_{1,1}^2$ \\
\hline\hline 
\textbf{Prediction} & 1.25 & 1.25 & 1.5 & 1.5 & 1.5 & 2 & 2 & 2 & 2 & 3 & 3 & 4 \\ \hline\hline

$(1,1)$ & \cellcolor{lightgray} & \cellcolor{lightgray} & \cellcolor{lightgray} & \cellcolor{lightgray} & \cellcolor{lightgray} & \cellcolor{lightgray} & 1.986 & 1.985 & 1.966 & 2.900 & 2.971 & 4.055 \\ \hline

$(1,2)$ & \cellcolor{lightgray} & \cellcolor{lightgray} & \cellcolor{lightgray} & 1.490 & 1.469 & 2.049 & 2.044 & 1.954 & 2.016 & 2.919 & 3.022 & 4.077 \\ \hline

$(1,3)$ & \cellcolor{lightgray} & \cellcolor{lightgray} & \cellcolor{lightgray} & 1.503 & 1.522 & 1.976 & 2.007 & 1.982 & 1.999 & 2.990 & 2.963 & 3.961 \\ \hline

$(1,4)$ & \cellcolor{lightgray} & \cellcolor{lightgray} & \cellcolor{lightgray} & \cellcolor{lightgray} & \cellcolor{lightgray} & \cellcolor{lightgray} & 2.053 & 2.019 & 2.034 & 2.968 & 2.966 & 3.950 \\ \hline

$(1,5)$ & \cellcolor{lightgray} & \cellcolor{lightgray} & \cellcolor{lightgray} & 1.509 & 1.500 & 2.012 & 2.015 & 1.993 & 1.947 & 2.988 & 2.962 & 3.928 \\ \hline

$(2,1)$ & \cellcolor{lightgray} & \cellcolor{lightgray} & \cellcolor{lightgray} & 1.510 & 1.464 & 1.985 & 2.011 & 2.008 & 2.034 & 3.028 & 2.986 & 4.040 \\ \hline

$(2,2)$ & 1.256 & 1.251 & 1.492 & 1.468 & 1.513 & 1.944 & \cellcolor{lightgray} & \cellcolor{lightgray} & \cellcolor{lightgray} & \cellcolor{lightgray} & \cellcolor{lightgray} & \cellcolor{lightgray} \\ \hline

$(2,3)$ & 1.246 & 1.244 & 1.509 & 1.500 & 1.507 & 2.034 & 2.011 & 1.981 & 2.036 & 2.980 & 2.980 & 4.017 \\ \hline

$(2,4)$ & \cellcolor{lightgray} & \cellcolor{lightgray} & \cellcolor{lightgray} & 1.485 & 1.487 & 1.953 & \cellcolor{lightgray} & \cellcolor{lightgray} & \cellcolor{lightgray} & \cellcolor{lightgray} & \cellcolor{lightgray} & \cellcolor{lightgray} \\ \hline

$(2,5)$ & 1.249 & 1.241 & 1.499 & 1.519 & 1.492 & 1.975 & 1.977 & 1.974 & 2.003 & 2.959 & 2.970 & 3.952 \\ \hline

$(3,1)$ & \cellcolor{lightgray} & \cellcolor{lightgray} & \cellcolor{lightgray} & 1.472 & 1.510 & 2.007 & 1.996 & 1.979 & 2.005 & 3.008 & 3.027 & 4.030 \\ \hline

$(3,2)$ & 1.248 & 1.244 & 1.490 & 1.484 & 1.469 & 2.005 & 1.977 & 2.024 & 2.036 & 2.987 & 2.972 & 4.033 \\ \hline

$(3,3)$ & 1.234 & 1.240 & 1.491 & 1.473 & 1.488 & 2.001 & \cellcolor{lightgray} & \cellcolor{lightgray} & \cellcolor{lightgray} & \cellcolor{lightgray} & \cellcolor{lightgray} & \cellcolor{lightgray} \\ \hline

$(3,4)$ & \cellcolor{lightgray} & \cellcolor{lightgray} & \cellcolor{lightgray} & 1.503 & 1.505 & 1.924 & 2.029 & 2.009 & 2.051 & 2.942 & 2.987 & 3.946 \\ \hline

$(3,5)$ & 1.244 & 1.254 & 1.482 & 1.506 & 1.490 & 2.007 & 2.120 & 1.936 & 1.949 & 3.031 & 2.924 & 4.056 \\ \hline

$(4,1)$ & \cellcolor{lightgray} & \cellcolor{lightgray} & \cellcolor{lightgray} & \cellcolor{lightgray} & \cellcolor{lightgray} & \cellcolor{lightgray} & 2.010 & 2.048 & 1.918 & 2.960 & 3.013 & 4.031 \\ \hline

$(4,2)$ & \cellcolor{lightgray} & \cellcolor{lightgray} & \cellcolor{lightgray} & 1.466 & 1.501 & 1.963 & \cellcolor{lightgray} & \cellcolor{lightgray} & \cellcolor{lightgray} & \cellcolor{lightgray} & \cellcolor{lightgray} & \cellcolor{lightgray} \\ \hline

$(4,3)$ & \cellcolor{lightgray} & \cellcolor{lightgray} & \cellcolor{lightgray} & 1.514 & 1.479 & 1.927 & 1.956 & 1.924 & 1.992 & 3.013 & 2.948 & 3.980 \\ \hline

$(4,4)$ & \cellcolor{lightgray} & \cellcolor{lightgray} & \cellcolor{lightgray} & \cellcolor{lightgray} & \cellcolor{lightgray} & \cellcolor{lightgray} & 1.954 & 1.967 & 2.032 & 3.050 & 2.949 & 4.014 \\ \hline

$(4,5)$ & \cellcolor{lightgray} & \cellcolor{lightgray} & \cellcolor{lightgray} & 1.508 & 1.459 & 1.990 & 2.105 & 1.985 & 1.993 & 2.936 & 2.896 & 4.002 \\ \hline

$(5,1)$ & \cellcolor{lightgray} & \cellcolor{lightgray} & \cellcolor{lightgray} & 1.486 & 1.477 & 2.010 & 2.031 & 1.980 & 2.068 & 2.968 & 2.973 & 3.976 \\ \hline

$(5,2)$ & 1.262 & 1.232 & 1.507 & 1.524 & 1.469 & 2.031 & 2.017 & 1.989 & 2.065 & 2.954 & 2.987 & 3.919 \\ \hline

$(5,3)$ & 1.234 & 1.249 & 1.524 & 1.488 & 1.498 & 2.006 & 2.015 & 1.961 & 1.988 & 3.048 & 2.940 & 3.925 \\ \hline

$(5,4)$ & \cellcolor{lightgray} & \cellcolor{lightgray} & \cellcolor{lightgray} & 1.501 & 1.490 & 1.988 & 2.030 & 2.036 & 1.995 & 2.936 & 3.012 & 4.002 \\ \hline

$(5,5)$ & 1.254 & 1.246 & 1.503 & 1.487 & 1.500 & 1.948 & \cellcolor{lightgray} & \cellcolor{lightgray} & \cellcolor{lightgray} & \cellcolor{lightgray} & \cellcolor{lightgray} & \cellcolor{lightgray} \\ \hline

\end{tabularx}
\caption{Predicted averages for stable subfamilies of $(a,d)$-monogenised fields, along with averages observed in random samples of fields from each subfamily.}
\label{table:app}
\end{table}

\clearpage

\bibliographystyle{abbrv}
\bibliography{bibfile}
\end{document}